%% file: main.tex
\documentclass[reqno]{amsart}
\usepackage{color}
\usepackage{comment}
\usepackage{graphicx}
\usepackage{dcolumn}
\usepackage{dsfont}
\usepackage{yfonts}

\usepackage{tikz} 
\usepackage{dcolumn}
\usetikzlibrary{arrows,positioning}
\tikzset{
    >=stealth',
    pil/.style={
           ->,
           thick,
           shorten <=2pt,
           shorten >=2pt,}
}

\RequirePackage[numbers]{natbib}
\RequirePackage[colorlinks=true, pdfstartview=FitV, linkcolor=blue,
  citecolor=blue, urlcolor=blue]{hyperref}
\RequirePackage{hypernat}
\usepackage{paralist}

\usepackage{graphicx}
\graphicspath{{./figures/}}
\usepackage{amsfonts}
\usepackage{amsmath}
\usepackage{amsthm}
\usepackage{amssymb}
\usepackage{amsbsy}
\usepackage{epsfig}
\usepackage{natbib, mathrsfs}
\usepackage{verbatim}
\usepackage[latin1]{inputenc}
\usepackage{mhequ}
\usepackage{algorithm}
\usepackage{algorithmic}

\usepackage{import}
\usepackage{enumitem} 
\usepackage{caption,subcaption}

\numberwithin{equation}{section}

\renewcommand{\leq}{\;\leqslant\;}                   
\renewcommand{\geq}{\;\geqslant\;}                   
\newcommand{\sumtwo}[2]{\sum_{\substack{#1 \\ #2}}} 
\newcommand{\sumthree}[3]{\sum_{\substack{#1 \\ #2 \\ #3}}} 

\newcommand{\1}{\mathds{1}}

\newcommand\weakconv{\stackrel{d}{\rightarrow}}

\newcommand\w{\omega}
\newcommand\dP{\mathbb{P}} 
\newcommand\dE{\mathbb{E}} 
\newcommand\pP{\bP} 
\newcommand\pE{\bE} 

\newcommand\ptp{\textgoth{Z}} 

\newcommand\red[1]{\textcolor{red}{#1}}

\newcommand{\R}{\mathbb R}
\newcommand{\N}{\mathbb N}
\newcommand{\bP}{{\bf P}}
\newcommand{\bE}{{\bf E}}

\newcommand{\Z}{{\ensuremath{\mathbb Z}} }

\definecolor{WowColor}{rgb}{.75,0,.75}
\definecolor{SubtleColor}{rgb}{0.9,0,0}

\newcounter{margincounter}

\newcounter{latercounter}

\newtheorem{theorem}{Theorem}[section]
\newtheorem{lemma}[theorem]{Lemma}
\newtheorem{proposition}[theorem]{Proposition}
\newtheorem{corollary}[theorem]{Corollary}
\newtheorem{remark}[theorem]{Remark}

\newtheorem*{question*}{Question}

\newtheorem*{remark*}{Remark}
\newtheorem*{idefinition*}{Definition}
\newtheorem*{example*}{Example}

\begin{document}
\title[Central limit theorems for subcritical directed polymers on $\Z^{2+1}$]{Central limit theorems for the ($2+1$)-dimensional directed polymer in the weak disorder limit} 
\author[S. Gabriel]{Simon Gabriel}
\address{S. Gabriel, Mathematics Institute, University of Warwick, Coventry CV4 7AL, United Kingdom } 
\email{simon.gabriel@warwick.ac.uk}
\keywords{directed polymer model, random environment, weak disorder, invariance principle, local limit theorem, functional central limit theorem.}
\subjclass[2010]{Primary: 82D60; Secondary: 60F17, 60K37, 82B44}
  
\begin{abstract}
In this article, we present an invariance principle for the paths of the directed random polymer in space dimension two in the subcritical intermediate disorder regime. More precisely, the distribution of diffusively rescaled polymer paths converges in probability to the law of Brownian motion when taking the weak disorder limit. So far analogous results have only been established for $d\neq 2$. 
Along the way, we prove a local limit theorem which allows us to factorise the point-to-point partition function of the directed polymer into a product of two point-to-plane partition functions.
\end{abstract}

\newpage

\maketitle


\setcounter{tocdepth}{1}
\tableofcontents

\import{sections/}{introduction.tex}

\import{sections/}{p2p.tex}

\import{sections/}{annealed.tex}

\import{sections/}{quenched.tex}

\import{sections/}{appendix.tex}

\bibliography{cites}
\bibliographystyle{alpha}

\end{document}

%% file: sections/introduction.tex
\section{Introduction and main results}
The directed polymer model describes a random walk
whose law is exponentially tilted by a random  environment. The strength of the environment is described by a non-negative parameter $\beta\geq 0$, corresponding to the inverse temperature, which we will refer to as disorder strength.
Individual models may vary but the most common definition is as follows:
consider the law $\pP_{N,x}$ of a nearest neighbour  random walk of length $N$ starting at $x\in\mathbb{Z}^d$. Furthermore, let $\w=(\w_{n,z})_{(n,z)\in \mathbb{N}\times \mathbb{Z}^d}$ be a family of random variables with law $\dP$ (independent of $\pP_{N,x}$).
For a fixed realisation of $\w$, the directed polymer measure of length $N$ and disorder strength $\beta\geq 0$ is then defined using the following change of measure
\begin{align*}
    \pP_{\beta, N,x}^{\w}(d S)
    :=
    \frac{1}{Z_{\beta,N}(0,x,N,\star)}
    \exp \left(
    \sum_{n=1}^N (\beta\, \w_{n,S_n}
    -\lambda(\beta))
    \right)
    \pP_{N,x}(d S),
\end{align*}
where $\lambda(\beta)$ is a positive constant, which we will fix in the subsequent section. The denominator
\begin{align*}
    Z_{\beta,N}(0,x,N,\star):=\pE_{N,x}\left[
    \exp \left(
    \sum_{n=1}^N(\beta\, \w_{n,S_n}
    -\lambda(\beta))
    \right)
    \right]
\end{align*}
is a (random) normalising constant, called the (point-to-plane) partition function, making $\pP_{\beta, N,x}^{\w}$ a probability measure. Here, $\star$ denotes the free boundary condition of the endpoint $S_N$ to take arbitrary values in $\Z^d$. The notation for the partition function might seem overloaded at this point but will become clear below. 

The original model goes back to the physics literature \cite{HH85} where directed random polymers were introduced to study the interface in two-dimensional Ising models with random interactions. 
Subsequently, the model was studied by the mathematical community \cite{IS88,Bo89} 
and attracted attention because of its application to stochastic partial differential equations (SPDEs), see for example \cite{BC95,BG97}.
But even on its own, the directed random polymer remains an interesting mathematical model. 
We refer to \cite{CSY04,Co17} and references therein for an overview of the polymer literature.\\

Henceforth we will choose $\pP_{N,x}$ to be the path measure of the symmetric simple random walk of length $N$ starting at $x\in \Z^d$. Without loss of generality we will assume $x=0$ and omit the space index, thus, simply writing $\pP_{N}=\pP_{N,0}$ and $\pE_N$ for the corresponding expectation. 
Moreover, let $\w=(\w_{n,z})_{(n,z)\in \mathbb{N}\times \mathbb{Z}^d}$ be a collection of i.i.d. real random variables satisfying 
\begin{align*}
    \dE[\w_{n,z}]=0, \qquad
    \dE[\w_{n,z}^2]=1, \qquad
    \lambda(\beta):=\log \dE[e^{\beta \w_{n,z}}]<\infty \;\;\forall \beta>0 \text{ small enough},
\end{align*}
which we will refer to as disorder. 
For technical reasons, we also require a concentration inequality for the law $\dP$. More precisely, we assume the existence of $\gamma>1$ and $C_1,C_2\in (0,\infty)$ such that for every $n\in\mathbb{N}$ and convex, $1$-Lipschitz $f:\R^n \mapsto \R$ we have 
\begin{align}\label{eq_conc_cond}
    \dP(|f(\w_1,\ldots,\w_n)-M_f|\geq t)\leq C_1 \exp \left(-\frac{t^{\gamma}}{C_2}\right),
\end{align}
where $(\w_i)_{1\leq i \leq n}$ is a subset of the family of random variables introduced above and $M_f$ is a median of $f(\w_1,\ldots,\w_n)$. Condition \eqref{eq_conc_cond} guarantees control on the negative tail of the environment and is for example satisfied whenever $\w$ is bounded or Gaussian; see \cite{Le01} for an even wider class of potential laws and more details. See also Remark~\ref{rem_weaker_neg_moments} for a discussion and possible approach to weaken this assumption.
\\

It was shown in \cite{Bo89}, that the partition function's limit, as $N$ diverges, is either positive or equal zero $\dP$-almost surely.
The monotonicity in $\beta$ of this behaviour was proven in \cite{CY06}: for arbitrary dimension $d$ there exists a $\beta_c=\beta_c(d)$ such that $\dP$-almost surely
\begin{align*}
    \lim_{N\to\infty} Z_{\beta,N}(0,0,N,\star)
    \begin{cases}
    >0 & \text{ if } \beta\in \{0\} \cup (0,\beta_c), \\
    =0 &  \text{ if } \beta> \beta_c.
    \end{cases}
\end{align*}
The subcritical phase is referred to as \emph{weak disorder regime}, whereas the supercritical phase is known as the \emph{strong disorder regime}. 
In particular, they established that $\beta_c=0$ whenever $d\leq 2$. 
In the interest of seeing a limit of  $Z_{\beta,N}(0,0,N,\star)$ with non-trivial fluctuations in $d\leq 2$, $\beta$ needs to be rescaled as a function of $N$ appropriately; see \cite{ALKQ14,CSZ17} and the account we give in Section~\ref{sec_background}.
Note that the concept of considering polymers when the disorder strength is scaled as a function of $N$ appeared already earlier in the physics literature \cite{Br00,CDR10}.

The case of dimension two is special: after rescaling $\beta=\beta_N\sim \widehat{\beta}\sqrt{\pi (\log N)^{-1}}$, for some $\widehat{\beta}\geq 0$, we see a phase transition in $\widehat{\beta}$ \cite{CSZ17}. More precisely, the limit $\lim_{N\to\infty} Z_{\beta_N,N}(0,0,N,\star)$ is strictly positive if $\widehat{\beta}\in (0,1)$ and vanishes if $\widehat{\beta}\geq 1$. 
Note that the $\log$-scaling of $\beta_N$, but not the corresponding phase transition in $\widehat{\beta}$, was already observed in \cite{BC98}.
In the one-dimensional case no such phase transition exists and the limiting partition function is strictly positive for $\beta=\beta_N\sim\widehat{\beta}N^{-1/4}$ and arbitrary $\widehat{\beta}>0$.

\subsection{Main result}

In this paper we consider directed random polymers in two space-dimensions. 
We are particularly interested in the asymptotic behaviour of paths under the measure $\pP^{\w}_{\beta,N}$ in the large $N$ limit. The first step is to determine the scaling exponent $\xi>0$ such that  $N^{-\xi}S_N$ under $\pP^{\w}_{\beta,N}$ has a non-trivial (random) limit, before determining the exact limiting distribution of the rescaled endpoint. Here, we say the polymer is diffusive if $\xi=\tfrac{1}{2}$, subdiffusive if $\xi<\tfrac{1}{2}$ and superdiffusive otherwise.
After establishing the endpoint distribution, the next natural step is to determine the limit law of the paths $N^{-\xi}(S_n)_{0\leq n\leq N}$. As we will see, this is not straightforward due to the random environment, cf. Remark~\ref{rem_flct_ivp}. \\

Consider the directed polymer measure $\pP_{\beta, N}^{\w}$ introduced above. The random measure is supported on the space $\{(S_n)_{n}\in (\mathbb{Z}^2)^{N+1} \}$, more precisely its support is given by the subset $\Omega_{0,N}$ of nearest neighbour paths starting at the origin. 
Because we use the space $C[0,1]:=C([0,1],\R^2)$ equipped with the supremum-norm as reference space for the paths, we introduce the mapping $\pi_N:\Omega_{0,N}\mapsto C[0,1]$ given by 
\begin{align}
X^{(N)}_t:=
    \left(\pi_N(S)\right)_t = \frac{1}{\sqrt{N}}\left(
 S_{\lfloor tN \rfloor} + (tN-\lfloor tN \rfloor)(S_{\lfloor tN \rfloor+1}-S_{\lfloor tN \rfloor})
\right),
\end{align}
which embeds discrete nearest-neighbour paths in the space of continuous functions by linearly interpolating between integer points and rescaling space-time diffusively. Furthermore, we equip $C[0,1]$ with the Borel $\sigma$-algebra $\mathcal{F}$ which implies measurability of the projection maps $\pi_N$.\\

The paper's main contribution is an invariance principle for rescaled polymer paths. Along the way we determine the limiting finite-dimensional distributions of the rescaled polymers.
We take the weak disorder limit, which is the large $N$ limit when scaling $\beta=\beta_N\sim \widehat{\beta}\sqrt{\pi (\log N)^{-1}}$ where $0<\widehat{\beta}<1$. More precisely, we consider
\begin{align}\label{eq_R_N}
\beta_N:=\frac{\widehat{\beta}}{\sqrt{R_N}}
, \quad \text{ where }\quad 
    R_N:= \sum_{n=1}^N \sum_{z\in\Z^2} \pP_N(S_n=z)^2= \tfrac{\log N}{\pi}+O(1)
\end{align}
denotes the \emph{replica overlap} of the simple random walk.

Our first result states that the finite-dimensional distributions of the quenched polymer measure $ \pi_N^{*}\pP_{\beta_N,N}^{\w}$ converge to the ones of Brownian motion with diffusion matrix $\tfrac{1}{\sqrt{2}}I_2$ in $\dP$-probability.

\begin{theorem}\label{theo_fdd_quenched}
Let $\widehat{\beta}\in (0,1)$ and $\beta_N$ be as in \eqref{eq_R_N}, then for any $0\leq t_1<\ldots<t_k\leq 1$ we have 
\begin{align*}
    \pi_{N}^{*}\pP_{\beta_{N},{N}}^{\w}\big((X_{t_1}, \ldots, X_{t_k})\in \cdot \big)
    \weakconv 
    \pP\big(
\tfrac{1}{\sqrt{2}}(W_{t_1}, \ldots, W_{t_k}\big)\in \cdot), \quad \text{ in }\dP\text{-probability},
\end{align*}
where $\pP$ denotes the Wiener measure on $C[0,1]$.
We wrote $\pi_N^{\ast}$ for the push-forward operation under $\pi_N$.
\end{theorem}

Convergence of finite-dimensional distributions and tightness of the disorder-averaged polymer measure (cf. Lemma~\ref{lem_annealed_tightness}) suffice to show that the limiting rescaled polymer paths have the same law as  Brownian motion.


\begin{theorem}\label{theo_main}
Let $\widehat{\beta}\in (0,1)$ and $\beta_N$ be as in \eqref{eq_R_N}.
Then
$$\pi_N^{*}\pP_{\beta_N, N}^{\w}
\weakconv \pP\big(\tfrac{1}{\sqrt{2}}W\in \cdot \,\big), \qquad \text{as }N\to\infty,
\quad \text{in }\dP\text{-probability},
$$ 
where $\pP$ denotes the Wiener measure on $C[0,1]$.
\end{theorem}

Despite the random polymer converging to a stochastic process (which is independent of the disorder) on a macroscopic scale, the disorder influences the behaviour of the polymer on small scales.
On the microscopic level, the disorder prevails and we can deduce a local limit theorem,
which allows to compare the microscopic polymer transition probabilities to the ones of Brownian motion, weighted by random multiplicative factors which depend on the rescaled transition space-time points. 

\begin{proposition}\label{prop_polymer_llt}
Let $\widehat{\beta}\in (0,1)$, $\beta_N$ as in \eqref{eq_R_N} and
\begin{align*}
    &(z_j)_{j=1}^k=(z_j(N))_{j=1}^k\in\Z^2 \text{ such that } \lim_{N\to\infty}\tfrac{z_j}{\sqrt{N}}=:x_j \text{ exists},\\
    (m_j)_{j=1}^k=(m_j(N))_{j=1}^k &\in\mathbb{N} \text{ such that } \lim_{N\to\infty} \tfrac{m_j}{N}=:t_j\in(0,1) \text{ exists with } 0< t_1<\cdots<t_k<1.
\end{align*}
Then, 
\begin{align}\label{eq_polymer_llt}
\Big(\frac{N}{2}\Big)^k\,
    \pP_{\beta_N,N}^{\w}(S_{m_1}=z_1,\ldots, S_{m_k}=z_k)
\weakconv
\prod_{j=1}^k
    :e^{Y^-(t_j,x_j)}:\,:e^{Y^+(t_j,x_j)}:
    \prod_{j=1}^k p_{\frac{1}{2}(t_j -t_{j-1})}(x_j-x_{j-1}),
\end{align}
where $Y_j^{\pm}$'s are i.i.d. centred Gaussians with variance $\log (1-\widehat{\beta}^2)^{-1}$. We used the shorthand notation $:e^{Y}:=e^{Y-\frac{1}{2}\dE[Y^2]}$ for the Wick exponential.
\end{proposition}

The local limit theorem above reinforces the picture that the considered subcritical intermediate disorder regime is indeed the region where effects of disorder start to emerge in a non-trivial way.

\begin{remark}
We excluded the case $t_k=1$ from Proposition~\ref{prop_polymer_llt}, since it would only give rise to a single factor $:e^{Y^-(1,x)}:$. However, the proof can be repeated almost verbatim to include this case.
\end{remark}



\subsection{Comparison to the literature}

We give a short overview on the literature of diffusivity of directed random polymers.
Consider the two dimensional case and $\beta_N$ scaled as in \eqref{eq_R_N}, then it was proven in \cite{CSZ17} that the diffusively rescaled field
$$
\{ Z_{\beta_N,N}(0,\lfloor\sqrt{N}x\rfloor,\lfloor tN\rfloor,\star)\, :\, t>0,\, x\in \R^2 \}
$$
converges to the solution of the stochastic heat equation with additive
space-time white noise, but to the author's best knowledge there are no results on diffusivity of the polymer paths in this case. 
Subcritical scalings $\beta_N^2 \ll R_N^{-1}$ in $d=2$ were considered by Feng \cite{Fe12}, who proved diffusivity of the polymer endpoint. 
However, under such subcritical scalings the partition function's variance vanishes in the large $N$ limit, which essentially brings us to the situation of setting $\widehat{\beta}=0$.
Theorem~\ref{theo_fdd_quenched} and \ref{theo_main}, on the other hand, consider a critical scaling under which the partition function converges to a non-trivial random variable and a transition (in $\beta$) of the polymer path behaviour is expected.
Our result not only covers the diffusivity of the polymer endpoint but fully determines the behaviour of the limiting polymer paths in the corresponding subcritical regime under diffusive scaling.

In dimension $d\geq 3$, diffusivity of the directed random polymer in the weak disorder regime was first proven to hold with probability one in \cite{IS88} for sufficiently small disorder strength. The works by Bolthausen \cite{Bo89} and Kifer \cite{Ki97} simplified and extended the original result further. The first invariance principles were deduced in \cite{AZ96,SZ96} where, for $\beta>0$ in the $L^2$-phase, almost sure convergence to the law of Brownian motion with dimension-dependent diffusion matrix was achieved. 
Later, it was extended to the full weak disorder regime by Comets and Yoshida \cite{CY06} in the sense of a functional central limit theorem which holds in probability. 
Theorem~\ref{theo_main} can be viewed as the analogous result in $d=2$.


Recently, Junk \cite{Ju21,Ju21b} gave an alternative proof of determining the limit of the polymer endpoint distribution (for bounded bond disorder) in $d\geq 3$ by introducing a comparison principle for partition functions of distinct parameters $\beta$. This allows them to transform the polymer endpoint distribution into the one of the simple random walk with a multiplicative error, that converges to one. 

A result for the case of $d=1$ was presented in papers by Alberts, Khanin and Quastel. In \cite{ALKQ14} they showed that transition probabilities of the discrete polymer measure admit a random limit when space-time is scaled diffusively. Because every $\beta>0$ lies in the strong disorder regime, they also relied on an intermediate disorder scaling: $\beta_N\sim \widehat{\beta}N^{-1/4}$.
In \cite{ALKQ13} they constructed the corresponding continuum polymer measure using the random field of transition probabilities that arose in \cite{ALKQ14}. 
As opposed to $d\geq 3$, the limiting distribution of polymer paths turned out to be singular w.r.t. the Wiener measure when scaled diffusively, while maintaining the same basic properties as Brownian motion.

An analogous result was also shown for the continuum disordered pinning model in \cite{CSZ_PTRF}. The advances, both for the pinning model and the $(1+1)$-dimensional directed polymer, then motivated to provide a general skeleton for the study of weak disorder scaling limits of discrete systems, see \cite{CSZ_EMS}.\\


A natural extension of the invariance principle in Theorem~\ref{theo_main} is to strengthen the result to $\dP$-a.s. convergence, similar to the results for the $L^2$-phase in $d\geq 3$ \cite{AZ96,SZ96}.
We want to point out that such results, holding with probability one, usually exploit the fact that the sequence of partition functions $(Z_{\beta, N})_N$ forms a martingale.
In the two-dimensional case, this property is lost due to the dependency $\beta= \beta_N$, which is why we do not expect our methods to yield an almost sure invariance principle.

The statement of the local limit theorem, Proposition~\ref{prop_polymer_llt}, is reminiscent of the construction in \cite{ALKQ13}, where the limiting field of partition functions was used to construct the continuum directed polymer in $(1+1)$-dimension. 
Because in $d=1$ the random (macroscopic) field induced by the limiting partition functions is continuous in its time and space points, the constructed polymer measure is the correct limiting object.
In $d=2$ this is however not the case anymore, which leads to substantial different behaviour of the polymer paths on a microscopic and macroscopic level. This (rough) structure of the partition function requires substantial work in order to establish the limiting polymer's behaviour.

\begin{remark}
Theorem~\ref{theo_fdd_quenched} and \ref{theo_main} should hold for a larger class of symmetric random walks which satisfy a local limit theorem in the sense of \cite[Hypothesis 2.4]{CSZ17} and their replica overlap fulfils $R_N \to \infty$ as a slowly varying function. 
\end{remark}

\begin{remark}\label{rem_cont}
Instead of studying discrete polymers, we could have similarly worked with polymers in the continuous space-time domain $[0,1]\times \R^2$ where the simple random walk is replaced by a Brownian motion and the disorder is given by a space-time white noise $\xi$. The corresponding polymer measure is then defined by the following Gibbsian tilt
$$
\pP^{\xi}_{\beta_{\varepsilon},\varepsilon}(d X)
\propto
\exp\left(
\beta_{\varepsilon}\int_0^1 \xi^{\varepsilon}(s,X_s)ds
- \frac{1}{2}\beta_{\varepsilon}^2 \|j\|_2^2 \varepsilon^{-2} 
\right)
\pP(d X),
$$
where $\xi^{\varepsilon}=\xi(t,\cdot)\ast j_{\varepsilon}$
with $j_{\varepsilon}=\varepsilon^{-2}j(\cdot/\varepsilon)$ for some $j\in C_c^{\infty}(\R^2)$.
As we turn off the mollification $\varepsilon\to 0$, we will need to tune the disorder strength like $ \beta_{\varepsilon}= \widehat{\beta}\sqrt{2\pi /\log \varepsilon^{-1}}$, similar to the discrete polymer. The corresponding version of Theorem~\ref{theo_main} in the continuum then reads as follows: for every $\widehat{\beta}\in(0,1)$
\begin{align}
\pP^{\xi}_{\beta_{\varepsilon},\varepsilon}
\weakconv \pP, \qquad \text{as }\varepsilon\to 0
\quad \text{in }\dP\text{-probability}.
\end{align}
Since all properties of the partition function of the discrete polymer also hold in the continuum, see \cite{CSZ17,CSZ_KPZ}, 
the proof of the above fact should follow along the same lines.
However, we refrain from giving a full proof to keep the paper at a reasonable length.
\end{remark}

\subsection{Background and outline of the proof}\label{sec_background}

Before presenting the outline of the proofs, we want to motivate the choice of the disorder scaling $\beta_N$. 
A second moment calculation of the partition function, which was already performed in \cite{IS88}, however for a polymer model of slightly different form, provides the following heuristic: 
\begin{align*}
    \dE[Z_{\beta,N}(0,0,N,\star)^2]
    &=\pE_N^{\otimes 2}  \Big[\prod_{n=1}^N e^{(\lambda(2\beta)-2\lambda(\beta))\mathds{1}_{S_n=S_{n}'}}\Big]
    =
    \pE_N^{\otimes 2}  \Big[\prod_{n=1}^N (1+\sigma^2\mathds{1}_{S_n=S_{n}'})\Big]\\
    &= \sum_{k=0}^N \sigma^{2k} \sum_{1\leq n_1<\cdots<n_k\leq N}
    \pE_N^{\otimes 2}  \Big[\prod_{i=1}^k \mathds{1}_{S_{n_i}=S_{n_i}'}\Big],
\end{align*}
where $S$ and $S'$ are two independent random walks of length $N$ and $\sigma$ is given by 
\begin{align}\label{eq_def_sigma}
    \sigma= \sigma(\beta):= \sqrt{e^{\lambda(2\beta)-2\lambda(\beta)}-1}.
\end{align}
\sloppy Note that $\lambda(2\beta)-2\lambda (\beta)\sim~\beta^2$  for small $\beta>0$ and therefore ${\lim_{\beta\to 0} \beta/ \sigma(\beta)=1}$ \cite[Equation (2.15)]{CSZ_KPZ}.
We upper bound the sum by ignoring the ordering of $n_i$'s  which yields 
\begin{align*}
    \dE[Z_{\beta,N}(0,0,N,\star)^2]
    \leq
    \sum_{k=0}^N \sigma^{2k} \Big(\sum_{n=1}^N
    \sum_{z\in\Z^2} \pP_N(S_n=z)^2\Big)^k.
\end{align*}
Recalling the replica overlap from \eqref{eq_R_N} and considering the fact that $\sigma(\beta)\sim \beta$ for small $\beta>0$, this suggests that the correct rescaling is given by $\beta= \beta_N :=\widehat{\beta}/\sqrt{R_N}$, whenever $\widehat{\beta}\in(0,1)$. Throughout the paper we will write $\sigma_N:=\sigma(\beta_N)$.

Indeed, it was proven in \cite{CSZ17} that under $\beta_N$ the partition function $Z_{\beta_N,N}(0,0,N,\star)$ converges to a non-trivial (random) limit whenever $\widehat{\beta}\in (0,1)$, see also \eqref{eq_csz_weakconv_of_p2l} below. 
Moreover, they noticed the existence of a transition on the finer scale with $\widehat{\beta}_c=1$ denoting the critical point where the $L^2$-norm of the partition function blows up in the limit.
Whenever $\widehat{\beta}\in (0,1)$ and $\beta_N$ is scaled as above, we speak of the \emph{intermediate weak disorder regime}.
\\

Theorem~\ref{theo_main}, in the present paper, states that rescaled polymer paths in the intermediate weak disorder limit behave like the ones of Brownian motion.
 However, this is not a straightforward consequence of the positivity of the limiting partition function,
but requires precise estimates quantifying the correlation structure of the limiting field.
The first step towards the main result is to show convergence of the finite-dimensional distributions of the (rescaled) directed polymer measure to the ones of Brownian motion, cf. Theorem~\ref{theo_fdd_quenched}. 
We begin by observing that for $m_1,\ldots,m_k\in \mathbb{N}$ and $z_1,\ldots,z_k\in\mathbb{Z}^2$ 
\begin{align*}
    &\pP_{\beta_N,N}^{\w}(S_{m_1}=z_1,\ldots, S_{m_k}=z_k)\\&\qquad=
    \frac{1}{Z_{\beta_N,N}(0,0,N,\star)}
    \prod_{j=1}^{k+1}
\pE_N\left[
e^{\sum_{n=m_{j-1}+1}^{m_j} \left( \beta_N w_{n,S_n}-\lambda(\beta_N) \right)}
\mathds{1}_{S_{m_j}=z_j}\Big| S_{m_{j-1}}=z_{j-1}\right],
\end{align*}
where $m_0=z_0=0$, $m_{k+1}=N$ and $z_{k+1}=\star$. For our purposes it will be more convenient to rewrite the above expression in terms of the expectations conditioned on both the start and end point, i.e. 
\begin{align}\label{eq_disc_fdd_as_p2p_product}
    &\pP_{\beta_N,N}^{\w}(S_{m_1}=z_1,\ldots, S_{m_k}=z_k) \nonumber
    \\&\qquad=
    \frac{1}{Z_{\beta_N,N}(0,0,N,\star)}
    \prod_{j=1}^{k+1}
    \ptp_{\beta_N,N}(m_{j-1},z_{j-1}\mid m_j,z_j)\,
q_{m_j-m_{j-1}}(z_{j}-z_{j-1}),
\end{align}
where we introduced the point-to-point partition functions
\begin{align*}
    \ptp_{\beta_N,N}(m_{j-1},z_{j-1}\mid m_j,z_j)
    :=
    \pE_N\left[
e^{\sum_{n=m_{j-1}+1}^{m_j} \left( \beta_N \w_{n,S_n}-\lambda(\beta_N) \right)}
\Big| S_{m_{j-1}}=z_{j-1},\, S_{m_j}=z_j\right]
\end{align*}
and the shorthand $q_n(z)$ denoting the transition probability $\pP_N(S_n=z)$ of the simple-random walk.

\begin{remark}\label{rem_p2p_endpoint}
Note that the point-to-point partition function $\ptp_{\beta_N,N}(0,0 \mid N,z)$ also takes the disorder at the endpoint into consideration. However, as it will turn out, it is more natural to compare the product of point-to-plane partition functions to the point-to-point partition function 
\begin{align}\label{eq_def_p2p_wo_endpoint}
    Z_{\beta_N,N}(0,0 \mid N,z)
    :=
     \pE_N\left[
e^{\sum_{n=1}^{N-1} \left( \beta_N \w_{n,S_n}-\lambda(\beta_N) \right)}
\Big| S_{0}=0,\, S_{N}=z\right],
\end{align}
not taking the endpoint-disorder into account. 
The distinction of the point-to-point partition functions' notation may be very subtle, but so is the difference between them.
In fact, the difference between \eqref{eq_def_p2p_wo_endpoint} and $\ptp_{\beta_N,N}(0,0 \mid N,z)$ vanishes in $L^2(\dP)$: 
\begin{align*}
    \| \ptp_{\beta_N,N}(0,0 \mid N,z)-Z_{\beta_N,N}(0,0 \mid N,z)\|_2^2
    &=
    \|Z_{\beta_N,N}(0,0 \mid N,z)\|_2^2 \dE[(e^{\beta_N \w_{N,z}-\lambda(\beta_N)}-1)^2]\\
    &=\|Z_{\beta_N,N}(0,0 \mid N,z)\|_2^2 (e^{\lambda(2\beta_N)-2\lambda(\beta_N)}-1),
\end{align*}
where we used the independence property of the disorder in the first equality. Because the first term on the r.h.s. is uniformly bounded in $N$ and $\lambda(2\beta_N)-2\lambda(\beta_N)\sim \beta_N^2$, the $L^2$-difference vanishes.
We refer to both $Z_{\beta_N,N}(0,0 \mid N,z)$ and $\ptp_{\beta_N,N}(0,0 \mid N,z)$ as point-to-point partition function since the meaning will be clear from the context. 
\end{remark}

With the slight abuse of notation, we will write for $0\leq s<t\leq 1$
$$
\ptp_{\beta_N,N}( sN ,y\mid tN,z) 
$$
instead of $\ptp_{\beta_N,N}( \lfloor sN \rfloor ,y\mid \lfloor tN \rfloor,z)$. Similarly, for the point-to-plane partition function. Furthermore, for future reference, we introduce the plane-to-point partition function, which is defined as 
\begin{align}\label{eq_def_l2p}
    Z_{\beta_N,N}( \widetilde{m}, \star, m,z)
    :=\pE_N\left[
e^{\sum_{n=\widetilde{m}}^{m-1} \left( \beta w_{n,S_n}-\lambda(\beta) \right)}
\Big| S_{m}=z\right].
\end{align}
One can think of it as the partition function of a polymer starting in $(m,z)$ and evolving backwards in time.
For convenience, we will refer to both point-to-line and line-to-point as point-to-plane whenever the context is clear. \\

Having representation \eqref{eq_disc_fdd_as_p2p_product} at hand, we see the necessity to understand the limiting behaviour of point-to-point partition functions, before analysing the finite-dimensional distributions of the polymer measure.
In order to manage such point-to-point partition functions, we prove a local limit theorem (Proposition~\ref{prop_fact_part_func}) which states that they can be approximated by the product of two point-to-plane partition functions, i.e.
\begin{align}\label{eq_llt_p2p_intro}
    \ptp_{\beta_N,N}(0,0\,|\,N,z)
    =
    Z_{\beta_N,N}(0,0,\tfrac{N}{2},\star)
    Z_{\beta_N,N}(\tfrac{N}{2},\star,N,z)
    +\varepsilon_{N}\,,
\end{align}
with $\varepsilon_N$ vanishing in $L^2(\dP)$. 
Factorisations of this nature were proven in $d\geq 3$ \cite{Si95, Va06,CNN20,LZ20}.
During completion of this paper, Nakajima and Nakashima \cite{NN21} proved independently a result similar to Proposition~\ref{prop_fact_part_func} in the continuous space-time setting, see also Remark~\ref{rem_NN}. They use  the local limit theorem to extend the class of SPDEs and initial conditions that admit Edwards-Wilkinson fluctuations.

We want to put particular emphasis on \cite{Si95} because of the similarities in their work and our proof of  Proposition~\ref{prop_fact_part_func}.
Their proof of the local limit theorem of the form \eqref{eq_llt_p2p_intro} uses the fact that polynomial chaos components,  \eqref{eq_polyn_chaos_p2p}, can be factorised using a single random walk transition probability which is of order $q_N(z)$  \cite[Theorem~2]{Si95}. Our proof of Lemma~\ref{cor_A_is_enough} resembles this approach.
Additionally, they explain how a central limit theorem for the end-point distribution can be obtained from the above factorisation \cite[Theorem~4]{Si95}.
It is interesting to note that Sinai proved the local limit theorem under the condition 
\begin{align}\label{eq_condition_sinai}
   \sigma(\beta)^2 R_\infty<1\,,
\end{align}
where $R_{\infty}:=\lim_{N\to\infty}R_N$, which is only finite in dimension $d\geq 3$. 
Formally generalising condition \eqref{eq_condition_sinai} w.r.t. the weak disorder limit in $d=2$, it reads $\sigma(\beta_N)^2 R_N<1$ which is equivalent to our assumption $\widehat{\beta}<1$ due to $\sigma(\beta_N)\sim \beta_N$.


As a consequence of  \eqref{eq_llt_p2p_intro}, the limiting distribution of point-to-point partition functions can be deduced from the corresponding point-to-plane partition functions approximating it. 
In \cite[Theorem 2.12]{CSZ17}, Caravenna, Sun and Zygouras proved that finite families of partition functions of directed polymers converge jointly to a multivariate log-normal distribution. More precisely, consider the collection of  space-time points $((n_i,z_i))_{1\leq i \leq k}= ((n_i(N),z_i(N)))_{1\leq i \leq k}$ such that for every $1\leq~i,j~\leq~k$
\begin{align*}
    &R_{N-n_i}/R_N\to 1, \quad\text{ as }N\to\infty,\\
    \text{ and }\quad 
    \lim_{N\to\infty}&
    R_{|n_i-n_j|\vee |z_i-z_j|^2}/R_N = \zeta_{i,j}\in [0,1] \quad\text{ exists}
    ,
\end{align*}
then 
\begin{align}\label{eq_csz_weakconv_of_p2l}
    \left(Z_{\beta_N,N}(n_i,z_i, N, \star)\right)_{1\leq i\leq k}\to (:e^{Y_i}:)_{1\leq i\leq k},
\end{align}
where $(Y_i)_{1\leq i\leq k}$ is a multivariate Gaussian with 
\begin{align}\label{eq_csz_weakconv_of_p2l_2}
    \dE[Y_i]=0 \quad \text{and}\quad 
    \dE[Y_i Y_j]=\log \frac{1-\widehat{\beta}^2\zeta_{i,j}}{1-\widehat{\beta}^2}
    \qquad \forall 1\leq i,j\leq k.
\end{align}
Particularly, the result holds for space-time points $((n_i,z_i))_{1\leq i \leq k}$ having positive macroscopic distance, in which case the corresponding  tuple $(Y_i)_{1\leq i\leq k}$ consists of independent Gaussians because $\zeta_{i,j}=\mathds{1}_{i\neq j}$. 

\begin{remark}\label{rem_self_averaging}
This fact is precisely the reason for the different behaviour of the limiting polymer law in $d=2$ compared to $d=1$.
In the one-dimensional case the partition functions have non-trivial dependency in the large $N$ limit for macroscopically separated space-time points, see \cite{ALKQ14}, leading to a path measure singular w.r.t. the Wiener measure.
In the two-dimensional setting \eqref{eq_csz_weakconv_of_p2l_2} implies that partition functions started from macroscopically separated points will have independent limits, leading to an self-averaging effect for the polymer measure. 
However, for points having vanishing macroscopic distance, the limiting field will have non-trivial dependency. 
\end{remark}


After having dealt with the approximation of point-to-point partition functions, we can move on to the convergence of quenched polymer marginals.
The greatest difficulty when dealing with the finite-dimensional marginals of the form \eqref{eq_disc_fdd_as_p2p_product} is the fact that none of the point-to-point partition functions is independent of the denominator $Z_{\beta_N,N}(0,0,N,\star)$.
We outline the approach taken in this paper; for the sake of simplicity, we only explain the following for the end-point distribution.

First, we prove that the limiting annealed polymer marginal, i.e. $\lim_{N\to\infty}\dE[\pP^\w_{\beta_N,N}(\tfrac{1}{\sqrt{N}}S_N \in \cdot )]$, agrees with the ones of Brownian motion, cf. Lemma~\ref{lem_annealed_law}. 
In fact, we show the much stronger result that the quenched marginal can be approximated in $L^1(\dP)$ by a simplified representation, without the partition function $Z_{\beta_N,N}(0,0,N,\star)$ in the denominator:
\begin{align}\label{eq_L1_approx_epmarginal}
    \lim_{N\to\infty}
    \Big\|
    \pP^\w_{\beta_N,N}(\tfrac{1}{\sqrt{N}}S_N \in B)
    - 
    \sum_{z\in \sqrt{N}B} Z_{\beta_N,N}(\tfrac{N}{2},\star, N,z) q_N(z)
    \Big\|_1=0\,,
\end{align}
using the factorisation in \eqref{eq_llt_p2p_intro}.
Here and throughout the paper $\sqrt{N}B$ denotes the set $\{ z\in\mathbb{Z}^2 \, :\, \tfrac{z}{\sqrt{N}}\in B\}$.
The expectation of the latter representation is immediate, which yields the annealed limit $\pP(\tfrac{1}{\sqrt{2}}W_1 \in B)$.

In order to conclude convergence of the quenched marginal, the natural next step would be to prove that $\sum_{z\in \sqrt{N}B} Z_{\beta_N,N}(\tfrac{N}{2},\star, N,z) q_N(z)$ converges to its mean in $L^1(\dP)$. Instead, we show the stronger convergence in $L^2(\dP)$, i.e. 
\begin{align*}
    \lim_{N\to\infty}
    \Big\|
    \sum_{z\in \sqrt{N}B} Z_{\beta_N,N}(\tfrac{N}{2},\star, N,z) q_N(z)
    - 
    \pP(\tfrac{1}{\sqrt{2}}W_1 \in B)
    \Big\|_2=0\,,
\end{align*}
since this reduces to a second moment calculation.
Estimating the second moment of $\sum_{z\in \sqrt{N}B} Z_{\beta_N,N}(\tfrac{N}{2},\star, N,z) q_N(z)$, requires careful evaluation of the partition functions' covariance structure, cf. Lemma~\ref{lem_self_avging}, leading to a law-of-large number like behaviour.
Together with \eqref{eq_L1_approx_epmarginal} this yields $L^1(\dP)$-convergence of the quenched end-point distribution.

When calculating second moments, we essentially  introduce a second independent copy of the polymer before averaging over the environment, see also \cite{CY06}.
A difference in our approach is that we exploit the fact that in $d=2$ the subcritical regime coincides with the $L^2$-regime. 
In other words, the positivity of the limiting partition function happens exactly in the regime where the $L^2(\dP)$-norm remains uniformly bounded. 
On the other hand, in $d\geq 3$ this is not the case as the $L^2$-regime is a strict subset of the weak disorder regime. Thus, Comets-Yoshida constructed, taking only advantage of the positivity of the limiting partition function, a (random) inhomogeneous Markov chain characterising the limiting discrete polymer measure with infinite time horizon.\\

After proving convergence of finite-dimensional distributions, cf. Theorem~\ref{theo_fdd_quenched}, we show that for any fixed function $F\in C_b(C[0,1])$ 
\begin{align}\label{eq_fclt}
    \pi_N^{\ast}\pE^\w_{\beta_N,N}[F(X)]
\to \pE[F(\tfrac{1}{\sqrt{2}}W )], \qquad \text{as }N\to\infty,
\quad \text{in }\dP\text{-probability},
\end{align}
by blending in ideas from the classical Donsker's invariance principle:
using tightness of the annealed polymer measure, it suffices to restrict the polymer paths to a compact set $K\subset C[0,1]$ when testing against a function $F\in C_b(C[0,1])$. 
The Stone-Weierstrass theorem then states that  $F$ can be approximated uniformly by cylinder functions on $K$, i.e. functions that only depend on finitely many marginals of the polymer path. 
Together with Theorem~\ref{theo_fdd_quenched} this yields the functional central limit theorem.

Lastly, we prove equivalence of functional central limit theorem and invariance principle by using a countable weak convergence determining family of functions, cf. Proposition~\ref{prop_equiv_fclt_ivp}.
This yields weak convergence of the polymer measures as stated in Theorem~\ref{theo_main}.
Note that the same argument allows rewriting the functional central limit theorem for $d\geq 3$ in \cite{CY06} in terms of an invariance principle, cf. Corollary~\ref{cor_ivp_dgeq3}.

\begin{remark}\label{rem_flct_ivp}
We want to stress that \eqref{eq_fclt} is not yet a classical invariance principle stating convergence of the polymer measure, but only a central limit theorem stating convergence when paths are tested against individual test functions. To emphasise this point, we note that a `true' invariance principle (in $\dP$-probability) reads as follows: for every sequence $(N_j)_{j\in\N}$ in $\N$ there exists a subsequence $(N_{j_m})_{m\in\N}$ and a set $\overline{\Omega}\subset \Omega$ of full measure such that 
$$
\pi_{N_{j_m}}^{\ast}\pE^\w_{\beta_{N_{j_m}},{N_{j_m}}}[F(X)]\to \pE[F(\tfrac{1}{\sqrt{2}}W)]\quad \forall F\in C_b(C[0,1]),
$$
for every $\w\in \overline{\Omega}$. In \eqref{eq_fclt} on the other hand, we fix a function $F\in C_b(C[0,1])$ for which there exists a subsequence and a set $\overline{\Omega}$ of full mass, on which the convergence holds. 
The dependency of $(N_{j_m})_{m\in\N}$ and $\overline{\Omega}$ on $F$ does not allow us to exchange the order of quantifiers without further reasoning.
\end{remark}

\subsection{Structure of the article}
The remainder of the paper is structured as follows. In Section~\ref{sec_p2p} we prove that point-to-point partition functions can be approximated by the point-to-plane partition functions which we introduced above. 
In Section~\ref{sec_annealed} we use this fact to prove convergence of the annealed finite-dimensional distributions of the polymer measure to the ones of a Brownian motion with diffusion matrix $\tfrac{1}{\sqrt{2}}I_2$, where $I_2$ denotes the identity matrix in $\R^{2\times 2}$. 
Together with a tightness argument this yields an annealed invariance principle, cf. Proposition~\ref{prop_annealed_clt}.
Section~\ref{sec_quenched} is divided into three parts. First, we prove  Theorem~\ref{theo_fdd_quenched}, where we use the self-averaging behaviour described above. Next, we present the proof of the invariance principle, cf. Theorem~\ref{theo_main}, where we exploit the tightness of the annealed polymer measure.
Lastly, we show the local limit theorem for the polymer marginals on microscopic scales, Proposition~\ref{prop_polymer_llt}.

\subsection{Notation}
Throughout the paper, $q_n(z)$ denotes the shorthand for the transition probabilities ${\pP_N(S_n=z)}$ of the simple random walk. Its continuous counterpart, the density of a centred Gaussian variable on $\R^2$ with variance $t$ is written as $p_t(x)$; we also write $\lambda(\cdot)=\lambda^{(d)}(\cdot)$ for the Lebesgue measure on $\R^d$.
Moreover, $\|\cdot\|_p$ denotes the $L^p(\dP)$-norm, i.e. 
$$
\|\cdot\|_p^p:= \dE[|\cdot|^p],
$$
for every $p>0$. 
Lastly, we will write $a_N\sim b_N$ for sequences if $\lim_{N\to\infty} \,a_N/ b_N=1$.

\subsection*{Acknowledgements}
We thank Nikos Zygouras for suggesting the problem and much valuable advice, as well as
Francesco Caravenna, Francis Comets, Dimitris Lygkonis and Rongfeng Sun for useful comments. We especially thank Stefan Junk for pointing out an error in an earlier version of the manuscript and the anonymous referees for valuable comments, that helped improving the main result to an invariance principle and simplify the presentation.
The author acknowledges financial support from EPSRC through grant EP/R513374/1.

%% file: sections/p2p.tex
\section{Approximation of the point-to-point partition function}\label{sec_p2p}

This section's main result is given by the following Proposition which states that point-to-point partition functions can be locally uniformly approximated by the product of a point-to-plane and a plane-to-point partition function.

\begin{proposition}\label{prop_fact_part_func}
Let $\widehat{\beta}\in (0,1)$,  
$x\in\R^2$ and $r>0$ be arbitrary, 
then 
for $0<s^+< t^-< 1$ we have
\begin{align*}
\sup_{\substack{z\in \sqrt{N}B(x,r)\\ \text{s.t. } q_N(z)>0}}
\left\|
\ptp_{\beta_N,N}(0,0\mid N,z) 
- Z_{\beta_N,N}(0,0,s^+N,\star) Z_{\beta_N,N}(t^-N,\star,N,z)
\right\|_{2}\to 0.
\end{align*}
The statement remains true when replacing $\ptp_{\beta_N,N}(0,0\mid N,z)$ with $Z_{\beta_N,N}(0,0\mid N,z)$, which we introduced in~\eqref{eq_def_p2p_wo_endpoint}.
\end{proposition}

\begin{remark}
Proposition~\ref{prop_fact_part_func} also holds for $s^+=1$ and $t^-=0$ w.r.t. $L^{1+\delta}(\dP)$-convergence for some $\delta>0$ small enough. This can be shown following the same steps in the proof of Proposition~\ref{prop_fact_part_func}. 
After completion of this paper, it was proved that $\sup_{N\in\N}\dE[\ptp_{\beta_N,N}(0,0, N,\star)^p]<\infty$ for arbitrary $p>0$ \cite{CZ21,LZ21}. The stronger moment estimates, allow to lift the mode of convergence from $L^{1+\delta}(\dP)$ to $L^{2}(\dP)$.
\end{remark}

Throughout this section, the point $x\in \R^2$ plays the role of the macroscopic endpoint of the polymer path. In particular, $\sqrt{N}B(x,r)$ includes all microscopic points which are close to $x$ on a macroscopic scale in the large $N$ limit which is important in Sections~\ref{sec_annealed} and \ref{sec_quenched}.
Note that Proposition~\ref{prop_fact_part_func} remains true when replacing the initial time $s=0$ and final time $t=1$ with arbitrary values $0\leq s< t\leq 1$, i.e. when considering partition functions $Z_{\beta_N,N}(sN,0 \mid tN,z)$ or $\ptp_{\beta_N,N}(sN,0 \mid tN,z)$.

\begin{remark}\label{rem_NN}
\sloppy A similar result was obtained recently and independently by Nakajima and Nakashima \cite[Theorem 2.8]{NN21}. They proved that the point-to-point partition function of a directed random polymer in the continuum can be approximated in $L^2(\dP)$ by the product of point-to-plane partition functions with mesoscopic time-horizon if the distance between start and terminal space-point is not too large.
Similar to the present paper, they show that contributions to the point-to-point partition function only come from the environment close to start and endpoint, before they replace the Brownian Bridge measure by two Brownian motions running independently forward and backward in time.
For the partition function of a polymer of length $N$, an equivalent result to \cite[Theorem 2.8]{NN21} for the discrete case would read as follows
\begin{align}\label{eq_nn_discrete}
    \sup_{\substack{z\in \Z^2:\, |z|\leq \sqrt{N\log N}\\ \text{s.t. }q_N(z)>0}}
    \left\|
\ptp_{\beta_N,N}(0,0\mid N,z) 
- Z_{\beta_N,N}(0,0,l_N,\star) Z_{\beta_N,N}(N-l_N,\star,N,z)
\right\|_{2}^2 \to 0\,,
\end{align}
for
$l_N=N^{1-(\log N)^{\gamma -1}}$ with $\gamma\in (0,1)$. 
In contrast to Proposition~\ref{prop_fact_part_func}, 
here the radius of uniformity is  $\sqrt{N\log N}$. This is due to keeping track of vanishing rates in their proof, which allows to strengthen the result. 
Likewise, exact evaluation of the quantities in Lemma~\ref{lem_rw_tk_properties} using the local limit theorem should allow to increase the radius of uniformity in  Proposition~\ref{prop_fact_part_func} to the same order. However, in regards of our main result this is not necessary.

\end{remark}

Before continuing, we remind the reader that the partition function $Z_{\beta_N,N}(0,0,N,\star)$ can be written in terms of a discrete chaos expansion \cite{CSZ17,CSZ_KPZ}:
\begin{align}\label{eq_p2l_chaos_expansion}
    Z_{\beta_N,N}(0,0,N,\star)
     &=
    \pE_N\Big[
    e^{\sum_{n=1}^N
    \sum_{z\in \Z^2}( \beta_N\w_{n,z}-\lambda(\beta_N))\mathds{1}_{S_n=z}}
    \Big]  \\
    &=
    \pE_N\Big[
    \prod_{n=1}^N
    \prod_{z\in \Z^2}
    (
    1+\sigma_N\eta_{n,z}\mathds{1}_{S_n=z})
    \Big]
    =
    1+
    \sum_{k=1}^N  Z_{\beta_N,N}^{(k)}(0,0, N,\star),\nonumber
\end{align} 
where $Z_{\beta_N,N}^{(k)}(0,0, N,\star)$ is defined below and $$\eta_{n,z}=\eta_{n,z}^{(N)}:=\frac{1}{\sigma_N}(e^{\beta_N \w_{n,z}-\lambda(\beta_N)}-1)$$ being centred i.i.d. random variables with unit variance.  
In the last equality of \eqref{eq_p2l_chaos_expansion} we expanded the products which gives rise to the $k$-th homogeneous chaos denoted by
\begin{align}\label{eq_p2l_kth_chaos}
    Z_{\beta_N,N}^{(k)}(0,0, N,\star):= \sigma_N^k
    \sumtwo{1\leq n_1<\cdots <n_k\leq N}{z_1,\ldots, z_k\in\Z^2}
    \left(\prod_{i=1}^k  q_{n_i-n_{i-1}}(z_i-z_{i-1})
    \eta_{n_i,z_i}\right),
\end{align}
with $(n_0,z_0)$ denoting the origin $(0,0)\in \N\times \Z^2$.
Note that the terms in the series expansion above are orthogonal in the sense that $\dE[Z_{\beta_N,N}^{(k)}(0,0, N,\star) Z_{\beta_N,N}^{(j)}(0,0, N,\star)]=0$, whenever $k\neq j$, due to the different number of disorder-terms considered. Throughout the paper, we will use this fact in second moment computations without further explanation.

An analogous expansion holds for the plane-to-point partition function $Z_{\beta_N,N}(0,\star,N,z)$ with 
\begin{align}\label{eq_l2p_kth_chaos}
    Z^{(k)}_{\beta_N,N}(0,\star,N,z):=
    \sigma_N^k
    \sumtwo{0\leq n_1<\cdots <n_k\leq N-1}{z_1,\ldots, z_k\in\Z^2}
    \left(\prod_{i=1}^k  q_{n_{i+1}-n_{i}}(z_{i+1}-z_{i})
    \eta_{n_i,z_i}\right),
\end{align}
where we assumed $(n_{k+1},z_{k+1})=(N,z)$ and used the symmetry of the transition probabilities of the simple random walk.
Similarly, for the point-to-point partition function we write
\begin{align*}
    Z_{\beta_N,N}(0,0\mid N,z)
    =
    1+
    \sum_{k=1}^N  Z_{\beta_N,N}^{(k)}(0,0\mid N,z)
\end{align*}
and 
\begin{align}\label{eq_polyn_chaos_p2p}
    Z_{\beta_N,N}^{(k)}(0,0\mid N,z)
    :=
    \sigma_N^k
    \sumtwo{1\leq n_1<\cdots <n_k\leq N-1}{z_1,\ldots, z_k\in\Z^2}
    \left(\prod_{i=1}^k  q_{n_i-n_{i-1}}(z_i-z_{i-1})
    \eta_{n_i,z_i}\right) \frac{q_{N-n_k}(z-z_k)}{q_N(z)}.
\end{align}
We point out that the point-to-plane partition function can be recovered by taking the average over all possible endpoints and include the endpoint-disorder: 
\begin{align*}
    Z_{\beta_N,N}(0,0, N,\star)
    =
    \sum_{z\in \Z^2}\ptp_{\beta_N,N}(0,0\mid N,z)q_N(z)
    =
    \sum_{z\in \Z^2}Z_{\beta_N,N}(0,0\mid N,z)e^{\beta_N \w_{N,z}-\lambda(\beta_N)}q_N(z).
\end{align*}

The proof of Proposition~\ref{prop_fact_part_func} relies on the fact that the main contribution of the point-to-point partition function comes from two mesoscopic sized subsets of the space-time domain around the start and terminal point. 
These space-time areas are the same as the ones giving main contribution to the partition functions $Z_{\beta_N,N}(0,0,N,\star)$, see \cite{CSZ_KPZ}. 
For a microscopic reference point $(n,z)\in \mathbb{N}\times \Z^2$, we define such sets both forward and backward in time: 
$$
A^{\pm}_N(n,z):=
\{
(m,y)\, : \, 
|y-z|\leq {N}^{1/2-a_N/4}
\; \text{and}\; 0 \leq \pm(m - n) \leq N^{1-a_N}
\},
$$
where $a_N:=(\log N )^{\gamma -1}$ for some $\gamma\in (0,1)$.

As already mentioned above, only samples inside of $A^{+}_N(0,0)
\cup A^{-}_N(N,z)$ will contribute to the $L^2$-limit of $Z_{\beta_N,N}^{(k)}(0,0\mid N,z)$. On this account, we introduce the following decomposition
$$
Z_{\beta_N,N}^{(k)}(0,0\mid N,z)
=
Z_{\beta_N,N}^{(k),A}(0,0\mid N,z)
+
\widehat{Z}_{\beta_N,N}^{(k)}(0,0\mid N,z),
$$
where $Z_{\beta_N,N}^{(k),A}(0,0\mid N,z)$ denotes the sum
\begin{align}\label{eq_Z_k_in_A}
\sigma_N^k
    \sumthree{1\leq n_1<\cdots <n_k\leq N-1}{z_1,\ldots, z_k\in\Z^2}{(n_i,z_i)\in A^{+}_N(0,0)
\cup A^{-}_N(N,z) }
    \left(\prod_{i=1}^k \eta_{n_i,z_i}\, q_{n_i-n_{i-1}}(z_i-z_{i-1})\right) \frac{q_{N-n_k}(z-z_k)}{q_N(z)}    
\end{align}
and $\widehat{Z}_{\beta_N,N}^{(k)}(0,0\mid N,z)$ the corresponding remainder. 
This restriction can be thought of as `turning off' the disorder outside the two boxes; for visualisation the reader may refer to Figure~\ref{figr} below. 
Similarly, we separate
\begin{align}\label{eq_decomp_p2l}
Z_{\beta_N,N}(0,0\mid N,z)
=
Z_{\beta_N,N}^{A}(0,0\mid N,z)
+
\widehat{Z}_{\beta_N,N}(0,0\mid N,z)
\end{align}
with $Z_{\beta_N,N}^{A}(0,0\mid N,z)= \sum_{k=0}^N Z_{\beta_N,N}^{(k),A}(0,0\mid N,z)$. 
Sometimes we will also use this notation for point-to-plane partition functions. Then, $Z_{\beta_N,N}^{(j),A}(0,0, N,\star)$ and $Z_{\beta_N,N}^{(j),A}(0,\star, N,z)$ will stand for the multi-linear polynomials in \eqref{eq_p2l_kth_chaos} and \eqref{eq_l2p_kth_chaos} with disorder restricted to $A^+_N(0,0)$ or  $A^-_N(N,z)$, respectively.
Moreover, we want to mention that $Z_{\beta_N,N}^{(j),A}(0,0, N,\star)=0$ for all $j>N^{1-a_N}$, hence, 
\begin{align}\label{eq_Z_A_only_up to_N_a_N}
Z_{\beta_N,N}^{A}(0,0, N,\star)
=
    \sum_{j=0}^{\lfloor N^{1-a_N} \rfloor}Z_{\beta_N,N}^{(j),A}(0,0, N,\star)\,,
\end{align}
similarly for the line-to-point partition function.
Lastly, note that the size of $A^{\pm}_N$ only depends on the level of approximation $N$ and not the point-to-point partition function's final time.
\\

We will frequently make use of higher-order and negative moment estimates of the partition function. We summarise equations (3.12), (3.14) and (3.15) from \cite{CSZ_KPZ} in the following lemma.

\begin{lemma}[Caravenna-Sun-Zygouras]\label{lem_hypcont_p2l}
Let $\widehat{\beta}\in (0,1)$, then
\begin{enumerate}[label=(\roman*)]
    \item  there exists a $\delta=\delta(\widehat{\beta})>0$ and a constant $C'_{\widehat{\beta}}<\infty$ such that for every $p\in [2,2+\delta]$
\begin{align*}
\sup_{N\in \mathbb{N}}\dE&\left[
Z_{\beta_N,N}(0,0,N,\star)^{p}
\right]\leq C'_{\widehat{\beta}}\, ,
\quad
\sup_{N\in \mathbb{N}}\dE\left[
Z_{\beta_N,N}^{A}(0,0,N,\star)^{p}
\right]\leq C'_{\widehat{\beta}}\\
&\quad\text{ and } \quad
\sup_{N\in \mathbb{N}}\dE\left[
\widehat{Z}_{\beta_N,N}(0,0,N,\star)^{p}
\right]\leq C'_{\widehat{\beta}}(a_N)^{p/2}.
\end{align*}

    \item for every $p> 0$ there exists a constant $\widehat{C}_{\widehat{\beta},p}>0$ such that 
    \begin{align*}
        \sup_{N\in \mathbb{N}}\dE&\left[
Z_{\beta_N,N}(0,0,N,\star)^{-p}
\right]\leq \widehat{C}_{\widehat{\beta},p}\quad \text{and}
\quad
\sup_{N\in \mathbb{N}}\dE\left[
Z_{\beta_N,N}^{A}(0,0,N,\star)^{-p}
\right]\leq \widehat{C}_{\widehat{\beta},p}.
    \end{align*}
\end{enumerate}
\end{lemma}
Because $Z_{\beta_N,N}(0,0,N,\star)$ and $Z_{\beta_N,N}(0,\star,N,z)$ have the same distribution, all statements in Lemma~\ref{lem_hypcont_p2l} hold for the plane-to-point partition function, too.

\begin{remark}\label{rem_weaker_neg_moments}
That Lemma~\ref{lem_hypcont_p2l}(ii) holds for arbitrary negative moments, is a consequence of
the concentration inequality \eqref{eq_conc_cond} and a left-tail estimate, see \cite[Proposition 3.1]{CSZ_KPZ}.
We will use this control, in the proof of Lemma~\ref{lem_annealed_law} and Proposition~\ref{prop_polymer_llt}, to separate products of partition functions using H\"older's inequality, cf. \eqref{eq_push_neg2pos_moments}. From the estimate \eqref{eq_push_neg2pos_moments} we see that control of either large negative or positive moments of the partition function is sufficient. Hence,
with an improved control on positive moments, it should be possible to push control from negative moments to positive ones.
It was recently shown that all positive moments of the partition function are uniformly bounded in $N$
\cite{CZ21,LZ21} (without the need of assumption \eqref{eq_conc_cond}), and it is expected that all statements in Lemma~\ref{lem_hypcont_p2l}(i) continue to hold for arbitrary $p>0$. 
Thus, assumption \eqref{eq_conc_cond} could then be replaced by  $\sup_{N\in \mathbb{N}}\dE\left[
Z_{\beta_N,N}(0,0,N,\star)^{-p}
\right]<\infty$ for some $p>2$. 
\end{remark}

\sloppy Although we will follow a similar approach in evaluating the limits of $Z_{\beta_N,N}^{A}(0,0\mid N,z)$ 
and $\widehat{Z}_{\beta_N,N}(0,0\mid N,z)$, we present their respective proofs separately for the sake of a more approachable presentation. We start with the chaos decomposition restricted to the macroscopically vanishing set $A^{+}_N(0,0)
\cup A^{-}_N(N,z)$.

\subsection{A single jump factorises the term}

For the term $Z_{\beta_N,N}^{(k),A}(0,0\mid N,z)$ we only consider sample points $(n_i,z_i)_{i=1}^{k}$ that lie inside  $A^{+}_N(0,0)
\cup A^{-}_N(N,z)$ by restricting the sum in \eqref{eq_Z_k_in_A}. Particularly, this implies for every polymer path with $k$ intermediate samples the existence of a single `large' jump across the valley separating $(0,0)$ and $(N,z)$.
This jump divides the polymer samples into two chains, one close to the starting point $(0,0)$, the other one close to the terminal point $(N,z)$. These two chains can be analysed independently. 
\begin{figure}[H]
    \centering
    \begin{tikzpicture}[scale = 0.8]
\draw (0,0) rectangle (2.5,2);
\draw (0,-0.5) -- (0,4);

\draw (12,-0.5) -- (12,4);
\draw (9.5,1.5) rectangle (12,3.5);

\filldraw[black] (0,1) circle (2pt) node[anchor=east] {\small $(0,0)$};
\filldraw[black] (12,2.5) circle (2pt) node[anchor=west] {\small $(N,z)$};
\filldraw[black] (2.3,1.5) circle (2pt) node[anchor=north west] {\small $\; (n_j,z_{j})$};
\filldraw[black] (9.8,3.3) circle (2pt) node[anchor=south east]  {\small $(n_{j+1},z_{j+1})\;\;$};

\draw   (0,1)  to[out=-20,in=-70] (0.7,1.5);
\filldraw[black] (0.7,1.5) circle (2pt);
\draw    (0.7,1.5) to[out=-20,in=-70] (1,0.7);
\filldraw[black] (1,0.7) circle (2pt);
\draw   (1,0.7) to[out=-20,in=-120] (1.8,0.7);
\draw[dotted]   (1.8,0.7) to[out=60,in=-120] (2,1);
\draw   (2,1) to[out=60,in=170] (2.3,1.5);

\draw   (9.8,3.3) to[out=-80,in=-150] (10.5,2.8);
\filldraw[black] (10.5,2.8) circle (2pt);
\draw   (10.5,2.8) to[out=-90,in=120] (10.7,2.4);
\draw[dotted]   (10.7,2.4) to[out=-69,in=-180] (11,2.3);
\draw   (11,2.3) to[out=0,in=120] (11.5,2);
\filldraw[black] (11.5,2) circle (2pt);
\draw   (11.5,2) to[out=0,in=-120] (12,2.5);

\draw (1.5,2) node[anchor=south] {\small $A_N^+(0,0)$};
\draw (10.5,1.5) node[anchor=north] {\small $A_N^-(N,z)$};

\draw[dashed] (2.3,1.5) -- (9.8,3.3);

\end{tikzpicture}
    \caption{Restricting samples to the macroscopic vanishing boxes $A_N^{\pm}$ implies the existence of a large jump over time.}
\end{figure}
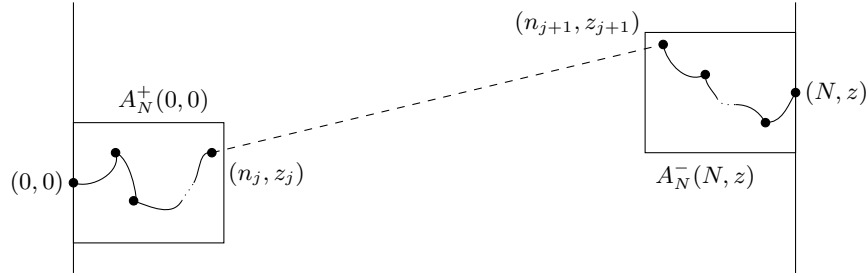
\label{figr}

Following the above explanation, it is reasonable to rewrite \eqref{eq_Z_k_in_A} by summing over all possible positions of this `large' jump:
\begin{align}\label{eq_rep_Z_k_with_j}
&Z_{\beta_N,N}^{(k),A}(0,0\mid N,z)\\
&\qquad=
\sigma_N^k
\sum_{j=0}^k
\sum_{
\substack{(n_i,z_i)\in A^{+}_N(0,0) \forall i\leq j \\ (n_i,z_i)\in A^{-}_N(N,y) \forall i > j\\ \text{s.t. }1\leq n_1<\cdots<n_k\leq N-1}}
\left(\prod_{i=1}^k 
q_{n_i-n_{i-1}}(z_i-z_{i-1}) \eta_{n_i,z_i}
\right)  \frac{q_{N-n_k}(z-z_k)}{q_N(z)},\nonumber
\end{align}
where $j$ is the largest index before the jump, see also Figure~\ref{figr}.
For the sake of clarity, we will abbreviate the conditions in the sum by $(n_i,z_i)_{i=1}^{k}\in \textgoth{A}(j)$ and write $\sum_{(n_i,z_i)_{i=1}^{k}\in \textgoth{A}(j)}$ coherently.
Most notably, when omitting the ratio ${q_{n_{j+1}-n_j}(z_{j+1}-z_j)}/{q_N(z)}$ in \eqref{eq_rep_Z_k_with_j} the r.h.s. simplifies to 
\begin{align*}
    \sum_{j=0}^k
    Z_{\beta_N,N}^{(j),A}(0,0, N,\star)
Z_{\beta_N,N}^{(k-j),A}(0,\star, N,z),
\end{align*}
and the two partition functions inside the sum are stochastically independent. Here, we used the notation defined in \eqref{eq_p2l_kth_chaos} and \eqref{eq_l2p_kth_chaos} and enhanced thereafter.
This motivates the statement of the following lemma. 
\begin{lemma}\label{cor_A_is_enough}
For all $x\in \mathbb{R}^2$ and $r>0$ we have
\begin{align*}
\sup_{\substack{z\in \sqrt{N}B(x,r)\\ \text{s.t. } q_N(z)>0}}
\left\|
Z_{\beta_N,N}^{A}(0,0\mid N,z)
-
Z_{\beta_N,N}^{A}(0,0,N,\star) Z_{\beta_N,N}^{A}(0,\star,N,z)
\right\|_{2}\to 0,
\end{align*}
as $N$ tends to infinity.
\end{lemma}

\begin{proof}
Let $x\in \mathbb{R}^2$ and $r>0$. We begin by noting that 
\begin{align}\label{eq_cor_A_enough_first_approx}
Z_{\beta_N,N}^{A}(0,0,N,\star) Z_{\beta_N,N}^{A}(0,\star,N,z)
=
    \sum_{k=0}^N 
\sum_{j=0}^k  Z_{\beta_N,N}^{(j),A}(0,0,N,\star) Z_{\beta_N,N}^{(k-j),A}(0,\star, N,z)\,.
\end{align}
This follows from the fact that $Z_{\beta_N,N}^{(j),A}(0,0,N,\star)=0$ whenever $j>N^{1-a_N}$ (similarly for $Z_{\beta_N,N}^{(k-j),A}(0,\star, N,z)$):
\begin{align*}
    \sum_{k=0}^N 
\sum_{j=0}^k  Z_{\beta_N,N}^{(j),A}(0,0,N,\star) Z_{\beta_N,N}^{(k-j),A}(0,\star, N,z)
=
\sum_{j=0}^{N^{1-a_N}} 
Z_{\beta_N,N}^{(j),A}(0,0,N,\star)
\sum_{k=j}^{j+N^{1-a_N}}   Z_{\beta_N,N}^{(k-j),A}(0,\star, N,z)\,,
\end{align*}
where we changed the order of the sums and added the restrictions $j,k-j\leq N^{1-a_N}$. After an index shift in the inner sum, we see that \eqref{eq_cor_A_enough_first_approx} holds true due to identity \eqref{eq_Z_A_only_up to_N_a_N}.

Thus, expanding the r.h.s. of  \eqref{eq_cor_A_enough_first_approx}, we have
\begin{align}\label{eq_prod_p2l_A_expanded}
   &  Z_{\beta_N,N}^{A}(0,0,N,\star) Z_{\beta_N,N}^{A}(0,\star,N,z) \nonumber \\
&\qquad=
\sum_{k=0}^N 
  \sum_{j=0}^k 
\sigma_N^{2j}
 \sum_{\substack{(n_i,z_i)_{i=1}^{j}\in A^+_N(0,0)\\ \text{s.t. }1\leq n_1<\cdots<n_j}}
 \bigg(\prod_{i=1}^j
q_{n_i-n_{i-1}}(z_i-z_{i-1}) \eta_{n_i,z_i}
\bigg)\\
&\qquad\qquad \times\sigma_N^{2(k-j)}\sum_{\substack{(\widetilde{n}_i,\widetilde{z}_i)_{l=1}^{k-j}\in A^-_N(N,z)\\ \text{s.t. }\widetilde{n}_1<\cdots< \widetilde{n}_{k-j}\leq N-1}}
\bigg(\prod_{i=1}^{k-j}
q_{\widetilde{n}_{i+1}-\widetilde{n}_{i}}(\widetilde{z}_{i+1}-\widetilde{z}_{i}) \eta_{\widetilde{n}_i,\widetilde{z}_i}
\bigg)
\,.\nonumber
\end{align}
\if{false}{
\red{
Let $x\in \mathbb{R}^2$, $r>0$ and assume that the following result holds 
\begin{align}\label{cc}
\sup_{\substack{z\in \sqrt{N}B(x,r)\\ \text{s.t. } q_N(z)>0}}
\left\|
Z_{\beta_N,N}^{A}(0,0\mid N,z)
-
\sum_{k=0}^N 
\sum_{j=0}^k  
Z_{\beta_N,N}^{(j),A}(0,0, N,\star)
Z_{\beta_N,N}^{(k-j),A}(0,\star, N,z)
\right\|_2
\to 0,
\end{align}
as $N\to\infty$. 
The lemma then follows by adding and subtracting the intermediate term introduced in \eqref{eq_split_poly_chain} and applying the triangle inequality, after proving that
\begin{align*}
    \sup_{\substack{z\in \sqrt{N}B(x,r)\\ \text{s.t. } q_N(z)>0}}
    \left\|
\sum_{k=0}^N 
\sum_{j=0}^k  Z_{\beta_N,N}^{(j),A}(0,0,N,\star) Z_{\beta_N,N}^{(k-j),A}(0,\star, N,z) - Z_{\beta_N,N}^{A}(0,0,N,\star) Z_{\beta_N,N}^{A}(0,\star,N,z)
\right\|_2
\end{align*}
vanishes in the large $N$ limit.
The statement in the previous equation follows using the following chain of (in)equalities 
\begin{align}\label{eq_approx_Z_A_using_imcomp_sum}
    &\left\|
\sum_{k=0}^N 
\sum_{j=0}^k  Z_{\beta_N,N}^{(j),A}(0,0,N,\star) Z_{\beta_N,N}^{(k-j),A}(0,\star, N,z) - Z_{\beta_N,N}^{A}(0,0,N,\star) Z_{\beta_N,N}^{A}(0,\star,N,z)
\right\|_2^2 \nonumber\\
&\qquad\qquad =
\sum_{j=0}^N  \dE[Z_{\beta_N,N}^{(j),A}(0,0,N,\star)^2]
\sum_{k=N-j+1}^N \dE[Z_{\beta_N,N}^{(k),A}(0,\star, N,z)^2]\nonumber
\\
&\qquad\qquad \leq
\sum_{j=0}^N (\sigma_N^2 R_N)^j
\sum_{k=N-j+1}^N  \dE[(Z_{\beta_N,N}^{(k),A}(0,\star, N,z) )^2],
\end{align}
where we used the fact that all $Z_{\beta_N,N}^{(j),A}(0,0,N,\star)$ and $Z_{\beta_N,N}^{(k),A}(0,\star, N,z)$
are pairwise uncorrelated in the first equality, before employing $$\dE[Z_{\beta_N,N}^{(j),A}(0,0,N,\star)^2]\leq (\sigma_N^2 R_N)^j,$$ for every $j\in\mathbb{N}$.
Note in particular that $\sum_{j=0}^{\infty} (\sigma_N^2 R_N)^j<\infty$ since $\sigma_N\sim \beta_N$ as $N$ grows large.
Because the second sum is uniformly bounded in $N$ and converges to zero for every choice of $j\in\mathbb{N}$, we may apply the dominated convergence theorem which yields that \eqref{eq_approx_Z_A_using_imcomp_sum} vanishes uniformly in $z$. \\}
}\fi
Now, using representations \eqref{eq_rep_Z_k_with_j} and \eqref{eq_prod_p2l_A_expanded},
we can estimate the second moment in the statement of the lemma, due to orthogonality, in terms of 
\begin{align}\label{eq_pre_crude_est}
&\sum_{k=0}^N
\dE\Big[\Big(
 Z_{\beta_N,N}^{(k),A}(0,0\mid  N,z)
-
\sum_{j=0}^k  
Z_{\beta_N,N}^{(j),A}(0,0, N,\star)
Z_{\beta_N,N}^{(k-j),A}(0,\star, N,z)
\Big)^2 \Big]\\
& \leq
\sum_{k=0}^N \sum_{j=0}^k
\dE\Big[
\big(
Z_{\beta_N,N}^{(j),A}(0,0, N,\star)
Z_{\beta_N,N}^{(k-j),A}(0,\star, N,z)
\big)^2
\sup_{\substack{(n_j,z_j)\in A^+(0,0) \\(n_{j+1},z_{j+1})\in A^-(N,z)}}
\Big(
\frac{q_{n_{j+1}-n_j}(z_{j+1}-z_j)}{q_N(z)}
-1
\Big)^2
\Big]
.\nonumber
\end{align}
Therefore, we only need to show that the ratio ${q_{n_{j+1}-n_j}(z_{j+1}-z_j)}/{q_N(z)}$ is negligible, which is intuitively clear since $(n_j,z_j)$ and $(n_{j+1},z_{j+1})$ are on a macroscopic level close to $(0,0)$ and $(N,z)$, respectively. Using the estimate
\begin{align}\label{eq_crude_estimate_ratio}
    \sup_{\substack{(n_j,z_j)\in A^+(0,0) \\(n_{j+1},z_{j+1})\in A^-(N,z)}}
    \left|
    \frac{q_{n_{j+1}-n_j}(z_{j+1}-z_j)}{q_N(z)}-1\right| 
    &\leq 
    \sup_{\substack{|N-n|< 2 N^{1-a_N} \\ |z-y|< 2{N}^{1/2-a_N/4}}}
\left|
\frac{q_n(y)}{q_N(z)}-1
\right|,
\end{align}
we can upper bound \eqref{eq_pre_crude_est} further (note that the r.h.s. of \eqref{eq_crude_estimate_ratio} does not depend on $k$ or $j$) and finally have 
\begin{align*}
    & \sum_{k=0}^N
\dE\Big[\Big(
 Z_{\beta_N,N}^{(k),A}(0,0\mid  N,z)
-
\sum_{j=0}^k  
Z_{\beta_N,N}^{(j),A}(0,0, N,\star)
Z_{\beta_N,N}^{(k-j),A}(0,\star, N,z)
\Big)^2 \Big]\\
&\leq
\sup_{\substack{|N-n|< 2 N^{1-a_N} \\ |z-y|< 2{N}^{1/2-a_N/4}}}
\left|
\frac{q_n(y)}{q_N(z)}-1
\right|^2
 \sum_{k=0}^N 
 \sum_{j=0}^k  
\dE\left[\left(
Z_{\beta_N,N}^{(j),A}(0,0, N,\star)
Z_{\beta_N,N}^{(k-j),A}(0,\star, N,z)
\right)^2 \right].
\end{align*}
Again, the sum of the right hand side is uniformly bounded in $N$ and $z$, cf. Lemma~\ref{lem_hypcont_p2l}. 
The last step consists of showing that the ratio of random-walk-transition-kernels is indeed uniformly close to one, i.e. 
\begin{align}\label{eq_another_trans_kern_est}
    \sup_{\substack{z\in \sqrt{N}B(x,r)\\ \text{s.t. } q_N(z)>0}}
    \sup_{\substack{|N-n|< 2 N^{1-a_N} \\ |z-y|< 2{N}^{1/2-a_N/4}}}
\left|
\frac{q_n(y)}{q_N(z)}-1
\right|\to 0, \quad \text{as }N\to\infty.
\end{align}
This follows directly from Lemma~\ref{lem_rw_tk_properties}(i) and finishes the proof. 
\end{proof}

\subsection{Multiple exceptional jumps are negligible}

Lemma~\ref{cor_A_is_enough} states that instead of looking at the point-to-point partition function restricted to have a single large jump from $A_N^+(0,0)$ to $A_N^-(N,z)$, it suffices to look at the product of two point-to-plane partition functions, one looking forward the other one looking backward in time. 
In order to prove the stronger result in Proposition~\ref{prop_fact_part_func} it remains to show that samples with points outside of $A_N^+(0,0) \cup A_N^-(N,z)$ do not contribute to the $L^2$-limit of the partition function. 
We say that samples of this kind have an exceptional jump, if 
there exists an index $1\leq \widehat{j}\leq k$ such that $(n_{\widehat{j}},z_{\widehat{j}})\notin A_N^+(0,0) \cup A_N^-(N,z)$. We will use 
$(n_i,z_i)_{i=1}^{k}\notin \textgoth{A}$ as a shorthand.

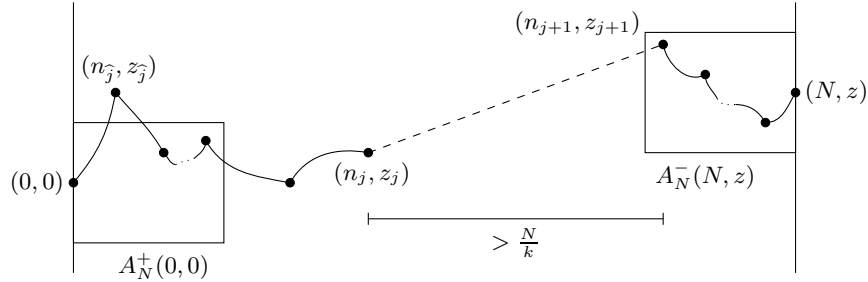
\begin{figure}[H]
    \centering
    \begin{tikzpicture}[scale = 0.8]
\draw (0,0) rectangle (2.5,2);
\draw (0,-0.5) -- (0,4);

\draw (4.9,0.3) -- (4.9,0.5);
\draw (9.8,0.3) -- (9.8,0.5);
\draw (4.9,0.4) -- (9.8,0.4);
\draw (7.35,0.4) node[anchor=north] {\small $>\frac{N}{k}$};

\draw (12,-0.5) -- (12,4);
\draw (9.5,1.5) rectangle (12,3.5);

\filldraw[black] (0,1) circle (2pt) node[anchor=east] {\small $(0,0)$};
\filldraw[black] (12,2.5) circle (2pt) node[anchor=west] {\small $(N,z)$};
\filldraw[black] (2.2,1.7) circle (2pt) ;
\filldraw[black] (4.9,1.5) circle (2pt) node[anchor=north] {\small $\; (n_j,z_{j})$};
\filldraw[black] (9.8,3.3) circle (2pt) node[anchor=south east]  {\small $(n_{j+1},z_{j+1})\;\;$};

\draw   (0,1)  to[out=50,in=-100] (0.7,2.5);
\filldraw[black] (0.7,2.5) circle (2pt) node[anchor=south ] {\small $\; (n_{\widehat{j}},z_{\widehat{j}})$};
\draw    (0.7,2.5) to[out=-50,in=120] (1.5,1.5);
\filldraw[black] (1.5,1.5) circle (2pt);
\draw   (1.5,1.5) to[out=-45,in=170] (1.7,1.3);
\draw[dotted]   (1.7,1.3) to[out=30,in=-160] (2,1.4);
\draw   (2,1.4) to[out=40,in=-90] (2.2,1.7);
\draw   (2.2,1.7) to[out=-60,in=170] (3.6,1);
\filldraw[black] (3.6,1) circle (2pt);
\draw   (3.6,1) to[out=60,in=170] (4.9,1.5);

\draw   (9.8,3.3) to[out=-80,in=-150] (10.5,2.8);
\filldraw[black] (10.5,2.8) circle (2pt);
\draw   (10.5,2.8) to[out=-90,in=120] (10.7,2.4);
\draw[dotted]   (10.7,2.4) to[out=-69,in=-180] (11,2.3);
\draw   (11,2.3) to[out=0,in=120] (11.5,2);
\filldraw[black] (11.5,2) circle (2pt);
\draw   (11.5,2) to[out=0,in=-120] (12,2.5);

\draw (1.5,0) node[anchor=north] {\small $A_N^+(0,0)$};
\draw (10.5,1.5) node[anchor=north] {\small $A_N^-(N,z)$};

\draw[dashed] (4.9,1.5) -- (9.8,3.3);

\end{tikzpicture}
    \caption{
    When considering $k$ samples, the existence of a time-jump that is at least of length $\tfrac{N}{k}$, still allows us to split samples into two groups (while paying a multiplicative constant $k$).
    At least one of the two groups will contain samples outside the boxes $A_N^{\pm}$.
    }
\end{figure}
\label{figr2}

\begin{lemma}\label{lem_non_box_limit_vanishes}
For every $x\in\R^2$ and $r>0$ we have
\begin{align*}
\lim_{N\to\infty}
\sup_{z\in \sqrt{N}B(x,r)}
\dE[
(\widehat{Z}_{\beta_N,N}(0,0\mid N,z))^2
]
= 0.
\end{align*}
\end{lemma}

\begin{proof}
Similar to the case where all samples are sufficiently close to the start and end point, we can again partition samples into two groups. 
First note that for every collection of $k$ samples there is at least one index $0\leq j\leq k$ such that $n_{j+1}-n_j\geq \frac{N}{k}$, cf. Figure~\ref{figr2} (here we use again the notation $n_{k+1}=N$). Consequently, we may write 
\begin{align*}
&\dE[
(\widehat{Z}_{\beta_N,N}^{(k)}(0,0\mid N,z))^2
]
\leq 
\sigma_N^{2k}
\sum_{j=0}^k
\sum_{(n_i,z_i)_{i=1}^{k}\notin \textgoth{A}}
 \mathds{1}\left\{
n_{j+1}-n_j\geq \frac{N}{k}
\right\}
\\
&\qquad \times \left(\prod_{i=1}^j
q_{n_i-n_{i-1}}^2(z_i-z_{i-1}) 
\right) 
\frac{q_{n_{j+1}-n_j}^2(z_{j+1}-z_j)}{q_N^2(z)}
\left(\prod_{i=j+2}^{k+1}
q_{n_i-n_{i-1}}^2(z_{i}-z_{i-1})
\right).
\end{align*}
Once more, we want to separate indices into two groups which requires the removal of 
a ratio of random walk transition kernels. 
We apply Lemma~\ref{lem_rw_tk_properties}(ii) which yields 
\begin{align}\label{eq_trans_k_estim}
    \frac{q_{n_{j+1}-n_j}^2(z_{j+1}-z_j)}{q_N^2(z)}
    \leq 
    \sup_{\substack{z\in \sqrt{N}B(x,r)\\ \text{s.t. } q_N(z)>0}}
        \sup_{\substack{n\geq N/k\\ z_{j+1}-z_j\in\Z^2}}
        \frac{q^2_n(z_{j+1}-z_j)}{q^2_N(z)}\leq C k^2,
\end{align}
with $C$ being a non-negative constant, independent of $N$ and $k$. 
\\

In the remainder of the proof we do not explicitly state results to hold locally uniformly in $z$. However, the reader should note that the statements remain true when adding the supremum over $\{z\in \sqrt{N}B(x,r) :\, q_N(z)>0 \}$ in front of all expressions below.\\

We just proved that
\begin{align}\label{eq_z_hat_estim_bfr_drop}
\dE[
(\widehat{Z}_{\beta_N,N}^{(k)}(0,0\mid N,z))^2
]
\leq C
k^2\,
\sigma_N^{2k}
\sum_{j=0}^k
\sumtwo{(n_i,z_i)_{i=1}^{k}\notin \textgoth{A}}{n_{j+1}-n_{j}\geq N/k}&
\left(\prod_{i=1}^j
q_{n_i-n_{i-1}}^2(z_i-z_{i-1}) 
\right) \\&\times
\left(\prod_{i=j+2}^{k+1}
q_{n_i-n_{i-1}}^2(z_{i}-z_{i-1})
\right),\nonumber
\end{align}
for some positive constant $C$. Since $(n_i,z_i)_{i=1}^{k}\notin \textgoth{A}$, we know there exists a $\widehat{j}\in\{1,\ldots, k\}$ such that $(n_{\widehat{j}},z_{\widehat{j}})\notin  A^+_N(0,0) \cup A^-_N(N,z)$. Note that $\widehat{j}$ does not necessarily agree with $j$ and 
we consider the two cases $\widehat{j}\leq j$ and $\widehat{j}> j$ separately. First, assume that  $\widehat{j}\leq j$, then we can estimate the product which does not contain $\widehat{j}$ as follows
\begin{align*}
    \sumtwo{0\leq n_{j+1}<\cdots <n_{k}\leq N-1}{z_{j+1},\ldots, z_{k}\in\Z^2}
    \left(\prod_{i=j+2}^{k+1}  q^2_{n_i-n_{i-1}}(z_i-z_{i-1})\right) 
    \leq 
    \Big(
    \sumtwo{1\leq n\leq N}{\widetilde{z}\in\Z^2}
     q^2_{n}(\widetilde{z})\Big)^{k-j}
    =
    R_N^{k-j},
\end{align*}
where we dropped the constraint of the $n_i$'s being ordered and also used the symmetry of the simple random walk transition probabilities in the first step. Thus, 
\begin{align}\label{eq_estimate_widehatjgeqj}
C
k^2\,
\sigma_N^{2k}
\sum_{j=0}^k
\sumtwo{(n_i,z_i)_{i=1}^{k}\notin \textgoth{A}}{n_{j+1}-n_{j}\geq N/k}& \mathds{1}_{\widehat{j}\leq j}
\left(\prod_{i=1}^j
q_{n_i-n_{i-1}}^2(z_i-z_{i-1}) 
\right) 
\left(\prod_{i=j+2}^{k+1}
q_{n_i-n_{i-1}}^2(z_{i}-z_{i-1})
\right) \nonumber \\
&\leq 
C
k^2\,
\sum_{j=0}^k
(\sigma_N^2 R_N)^{k-j}
\sigma_N^{2j}
\sumthree{1\leq n_1<\cdots<n_j\leq N}{(z_i)_{i=1}^{j}\in (\Z^2)^j}{\exists\, \widehat{j} \text{ s.t. } (n_{\widehat{j}},z_{\widehat{j}})\notin A^+_N(0,0)}
\left(\prod_{i=1}^j
q_{n_i-n_{i-1}}^2(z_i-z_{i-1}) 
\right) 
\\
&=
C
k^2\,
\sum_{j=0}^k
(\sigma_N^2 R_N)^{k-j}
\dE[(
\widehat{Z}_{\beta_N,N}^{(j)}(0,0,N,\star))^2].\nonumber
\end{align}
We can perform the same estimates in the case of $\widehat{j}> j$ (in fact, we overestimate with $\widehat{j}\geq j$) with the roles of $A_N^+(0,0)$ and $A_N^-(N,z)$ reversed. Together with \eqref{eq_z_hat_estim_bfr_drop} and \eqref{eq_estimate_widehatjgeqj} this yields
\begin{align}\label{eq_hatZ_k_upperbound}
    \dE[
(\widehat{Z}_{\beta_N,N}^{(k)}(0,0\mid N,z))^2
]
&\leq 
Ck^2
\sum_{j=0}^k
\Big(
(\sigma_N^2 R_N)^{k-j}
\dE[(
\widehat{Z}_{\beta_N,N}^{(j)}(0,0,N,\star))^2]
+
(\sigma_N^2 R_N)^{j}
\dE[(
\widehat{Z}_{\beta_N,N}^{(k-j)}(0,\star,N,z))^2]
\Big)\nonumber\\
&=2Ck^2
\sum_{j=0}^k
(\sigma_N^2 R_N)^{k-j}
\dE[(
\widehat{Z}_{\beta_N,N}^{(j)}(0,0,N,\star))^2]\,,
\end{align}
where we used $\dE[(
\widehat{Z}_{\beta_N,N}^{(j)}(0,0,N,\star))^2]=\dE[(
\widehat{Z}_{\beta_N,N}^{(j)}(0,\star,N,z))^2]$ in the last step.
Overall, using \eqref{eq_hatZ_k_upperbound} and orthogonality of the polynomial chaos components, we have 
\begin{align}
        \dE[
(\widehat{Z}_{\beta_N,N}(0,0\mid N,z))^2
]&\leq  2C
\sum_{k=1}^N k^2\,
\sum_{j=1}^k
(\sigma_N^2 R_N)^{k-j}
\dE[(
\widehat{Z}_{\beta_N,N}^{(j)}(0,0,N,\star))^2].
\end{align}
Switching the order of the two sums and performing an index shift, we have 
\begin{align*}
\sum_{k=1}^N k^2\,
\sum_{j=1}^k&
(\sigma_N^2 R_N)^{k-j}
\dE[(
\widehat{Z}_{\beta_N,N}^{(j)}(0,0,N,\star))^2]
= 
\sum_{j=1}^N 
\dE[(
\widehat{Z}_{\beta_N,N}^{(j)}(0,0,N,\star))^2]
\sum_{k=0}^{N-j}
(k+j)^2\,
(\sigma_N^2 R_N)^{k}
\end{align*}
which overall yields 
\begin{align}\label{eq_p2p_widehatZ_vanishes_upperbound}
     \dE[
(\widehat{Z}_{\beta_N,N}(0,0\mid N,z))^2
]&\leq  C 
\sum_{j=1}^N j^2
\dE[(
\widehat{Z}_{\beta_N,N}^{(j)}(0,0,N,\star))^2].
\end{align}
Here we used the fact that $$
\sum_{k=0}^{N-j}
(k+j)^2\, (\sigma_N^2 R_N)^{k}
\leq j^2\Big(1 +4 \sum_{k=1}^{\infty}
k^2\, (\sigma_N^2 R_N)^{k}\Big)
\to j^2\Big(1+4 \sum_{k=1}^{\infty}
k^2\, \widehat{\beta}^{2k}\Big), \qquad \text{as }N\to\infty,$$
for all $j\geq 1$. Because the series converges, recall that $\lim_{N\to \infty} \sigma_N^2 R_N =\widehat{\beta}^2<1$, the series in on the l.h.s. is uniformly bounded in $N$ and can be absorbed in the constant $C$.
In Lemma~\ref{lem_exp_decay_hatZ_k}, in the Appendix, we deduce that
$$
\dE[(
\widehat{Z}_{\beta_N,N}^{(j)}(0,0,N,\star))^2]
\leq 
C\, j^2\, (\sigma_N^2 R_N)^{\frac{j}{2}} \, a_N
\sim 
C\, j^2\, \widehat{\beta}^{j} \, a_N\,,
$$
which is a direct implication of the estimates in \cite[Section 3.4]{CSZ_KPZ}. Hence, 
\begin{align*}
    \sum_{j=1}^N j^2
\dE[(
\widehat{Z}_{\beta_N,N}^{(j)}(0,0,N,\star))^2]
\leq 
C \, a_N\, \sum_{j=1}^N j^4 \widehat{\beta}^j\,,
\end{align*}
which yields that \eqref{eq_p2p_widehatZ_vanishes_upperbound} vanishes in the large $N$ limit, since $a_N\to 0$. This completes the proof.
\end{proof}
It is worth pointing out that the estimates used in \eqref{eq_another_trans_kern_est} and \eqref{eq_trans_k_estim} hold because $z\in\sqrt{N}B(x,r)$ is macroscopically bounded which prevents blow ups of $(N\,q_N(z))^{-1}$. This is the reason why Proposition~\ref{prop_fact_part_func} only holds locally uniformly.\\

We are now ready to prove this section's main result, namely the point-to-point partition function can in fact be factorised into a product of two point-to-plane partition functions.

\begin{proof}[Proof of Proposition~\ref{prop_fact_part_func}]
First, we ignore the $\sup_{z\in \sqrt{N}B(x,r),\, q_N(z)>0}$ and only consider a single, fixed $z\in \Z^2$ such that $q_N(z)>0$.
As already mentioned in Remark~\ref{rem_p2p_endpoint}, $\ptp_{\beta_N,N}(0,0\mid N,z)$ can be approximated arbitrary well in $L^2(\dP)$ by $Z_{\beta_N,N}(0,0\mid N,z)$ in the large $N$ limit, which is why we can restrict ourselves to the latter in this proof.

Recall that $0<s^+< t^-<1$.
Using the triangle inequality yields the estimate
\begin{align*}
&\left\|
{Z_{\beta_N,N}(0,0\mid N,z)} - Z_{\beta_N,N}(0,0,s^+N,\star) Z_{\beta_N,N}(t^-N,\star,N,z)
\right\|_{2}\\
&\qquad\qquad =
\left\|
({\widehat{Z}_{\beta_N,N}(0,0\mid N,z)}+ {Z_{\beta_N,N}^{A}(0,0\mid N,z)})- Z_{\beta_N,N}(0,0,s^+N,\star) Z_{\beta_N,N}(t^-N,\star,N,z)
\right\|_{2}\\
&\qquad\qquad\leq
\left\|
{\widehat{Z}_{\beta_N,N}(0,0\mid N,z)}
\right\|_{2}
+
\left\|
{Z_{\beta_N,N}^{A}(0,0\mid N,z)} - Z_{\beta_N,N}^{A}(0,0,s^+N,\star) Z_{\beta_N,N}^{A}(t^-N,\star,N,z)
\right\|_{2}\\
&\qquad\qquad\qquad+
\left\|
Z_{\beta_N,N}^{A}(0,0,s^+N,\star) Z_{\beta_N,N}^{A}(t^-N,\star,N,z) - Z_{\beta_N,N}(0,0,s^+N,\star) Z_{\beta_N,N}(t^-N,\star,N,z)
\right\|_{2}.
\end{align*}
The first and second term vanish due to Lemma~\ref{lem_non_box_limit_vanishes} and Lemma~\ref{cor_A_is_enough}, respectively. 
Intuitively, the third term is negligible due to the results in \cite{CSZ17} which state that only samples inside the boxes $A^{\pm}_N$ contribute to the $L^2$-limit of the point-to-plane partition functions. We begin by estimating 
\begin{align*}
&\left\|
Z_{\beta_N,N}^{A}(0,0,s^+N,\star) Z_{\beta_N,N}^{A}(t^-N,\star,N,z) - Z_{\beta_N,N}(0,0,s^+N,\star) Z_{\beta_N,N}(t^-N,\star,N,z)
\right\|_{2}\\
&\qquad\qquad
\leq \left\|\left(
Z_{\beta_N,N}^{A}(0,0,s^+N,\star)  - Z_{\beta_N,N}(0,0,s^+N,\star) \right) Z_{\beta_N,N}^{A}(t^-N,\star,N,z)\right\|_{2}\\
&\qquad\qquad\qquad+
\left\|Z_{\beta_N,N}(0,0,s^+N,\star)
\left(
Z_{\beta_N,N}^{A}(t^-N,\star,N,z) -  Z_{\beta_N,N}(t^-N,\star,N,z)
\right)
\right\|_{2}\\
&\qquad\qquad =
\left\|
\widehat{Z}_{\beta_N,N}(0,0,s^+N,\star)\right\|_{2}
\left\|Z_{\beta_N,N}^{A}(t^-N,\star,N,z)\right\|_{2}\\
&\qquad \qquad \qquad+
\left\|Z_{\beta_N,N}(0,0,s^+N,\star)\right\|_{2}
\left\|\widehat{Z}_{\beta_N,N}(t^-N,\star,N,z)
\right\|_{2},
\end{align*}
where we used independence of the disorder on the disjoint time intervals $[0,s^+]$ and $[t^-,1]$ in the last step.
The second moments of $\widehat{Z}_{\beta_N,N}$ converge to zero, see Lemma~\ref{lem_hypcont_p2l}, whereas the 2-norms of ${Z}^A_{\beta_N,N}$ and ${Z}_{\beta_N,N}$ are uniformly bounded in $N$.
Note that all statements in the proof hold uniformly on balls $\sqrt{N}B(x,r)$, therefore the local uniformity 
in the statement of the Proposition follows. This completes the proof.
\end{proof}

Before we move on to the next section, we state and prove the analogous result of Lemma~\ref{lem_hypcont_p2l}(i) for the point-to-point partition function:

\begin{corollary}[Hypercontractivity]\label{cor_hypercontract}
Let $x\in \mathbb{R}^2$ and $r>0$ be arbitrary. 
For every $\widehat{\beta}<1$ there exists a $\delta=\delta(\widehat{\beta})>0$ and $C>0$ such that 
\begin{enumerate}[label=(\roman*)]
    \item $
\sup_{N\in \mathbb{N}}\sup_{z\in \sqrt{N}B(x,r)} \dE[
Z_{\beta_N,N}(0,0\mid N,z)^{2+\delta}
]\leq C,
$
\item $
\sup_{N\in \mathbb{N}}\sup_{z\in \sqrt{N}B(x,r)} \dE[
\ptp_{\beta_N,N}(0,0\mid N,z)^{2+\delta}
]\leq C.
$
\end{enumerate}
\end{corollary}

\begin{proof}
Let $x\in \mathbb{R}^2$ and $r>0$. 
First, we prove statement (i) for the point-to-point partition function $Z_{\beta_N,N}(0,0\mid N,z)$, ignoring the endpoint-disorder. 
From the proof of  Proposition~\ref{prop_fact_part_func} and Lemma~\ref{lem_hypcont_p2l}, we know that $\sup_{z\in \sqrt{N}B(x,r)}\dE[
Z_{\beta_N,N}(0,0\mid N,z)^{2}
]$ is uniformly bounded in $N$.
Moreover, the point-to-point partition function still has the form of a multi-linear polynomial.
In order to lift the boundedness to the $(2+\delta)$-moment for $\delta>0$ sufficiently small, we apply the hypercontractivity property, see for example Appendix B in \cite{CSZ_KPZ}. It is easy to check that the necessary conditions are still satisfied. Thus, we may estimate the $p$-th moment, $p\geq 2$, by 
\begin{align}\label{eq_hypercontract_p2p}
     \dE[
Z_{\beta_N,N}(0,0\mid N,z)^{p}
]
\leq \Big(\sum_{k=0}^N c_{p}^{2k}
\dE[
Z_{\beta_N,N}^{(k)}(0,0\mid N,z)^{2}
]
\Big)^{\frac{p}{2}}.
\end{align}
Here, $c_{p}$ is a constant which only depends on $p$ and the law of the noise $\eta^{(N)}$.
In \cite[Theorem B.1]{CSZ_KPZ} it was additionally shown that $\lim_{p\downarrow 2}c_p=1$. 
In order to show that the sum is uniformly bounded in $N$, we split it into the two (by now well known) parts:
\begin{align}\label{eq_hyp_contrac_split}
\sum_{k=0}^N c_{p}^{2k}
\dE[
Z_{\beta_N,N}^{(k)}(0,0\mid N,z)^{2}
]
= 
\sum_{k=0}^N c_{p}^{2k}
\dE[
Z_{\beta_N,N}^{A,(k)}(0,0\mid N,z)^{2}
]
+
\sum_{k=0}^N c_{p}^{2k}
\dE[
\widehat{Z}_{\beta_N,N}^{(k)}(0,0\mid N,z)^{2}
],
\end{align}
due to orthogonality.
The second term on the r.h.s. can be upper bounded using \eqref{eq_hatZ_k_upperbound} and Lemma~\ref{lem_exp_decay_hatZ_k}
\begin{align*}
    \sum_{k=0}^N c_{p}^{2k}
\dE[
\widehat{Z}_{\beta_N,N}^{(k)}(0,0\mid N,z)^{2}
]
&\leq C
\sum_{k=0}^N c_{p}^{2k}
k^2
\sum_{j=0}^k
(\sigma_N^2 R_N)^{k-j}
\dE[(
\widehat{Z}_{\beta_N,N}^{(j)}(0,0,N,\star))^2]\\
&\leq 
a_N\, C
\sum_{k=0}^N c_{p}^{2k}
k^2
\sum_{j=0}^k
(\sigma_N^2 R_N)^{k-j}
 (\sigma_N^2 R_N)^{j/2} j^2
 \leq 
 a_N\, C
\sum_{k=0}^N 
k^5
(c_{p}^{4}\sigma_N^2 R_N)^{k/2}
\,,
\end{align*}
and vanishes because $c_p^2 \widehat{\beta}<1$ for $p\geq 2$ small enough, where we used 
the fact that $c_p^2 \sigma_N \sqrt{R_N} \sim c_p^2\widehat{\beta}$. The first term on the r.h.s. of \eqref{eq_hyp_contrac_split}, on the other hand, can be estimated by
\begin{align*}
    &\sum_{k=0}^N c_{p}^{2k}
\dE[
Z_{\beta_N,N}^{A,(k)}(0,0\mid N,z)^{2}
]\\
&\qquad\leq 
\bigg(
\sup_{\substack{|N-n|<2N^{1-a_N}\\ |z-y|<2N^{1/2-a_N/4}}} \frac{q_n(y)}{q_N(z)}\bigg)
\sum_{k=0}^N c_{p}^{2k}
\sum_{j=0}^k
\dE[
Z_{\beta_N,N}^{A,(j)}(0,0, N,\star)^{2}
]
\dE[
Z_{\beta_N,N}^{A,(k-j)}(0,\star, N,z)^{2}
]\\
&\qquad\leq 
\bigg( \sup_{\substack{|N-n|<2N^{1-a_N}\\ |z-y|<2N^{1/2-a_N/4}}} \frac{q_n(y)}{q_N(z)} \bigg)
\sum_{k=0}^N (k+1) c_{p}^{2k} (\sigma_N^{2}R_N)^k,
\end{align*}
where we performed a similar estimate as in the proof of Lemma~\ref{cor_A_is_enough} in the first inequality, recall display \eqref{eq_rep_Z_k_with_j}. In the second step we simply used the fact that $\dE[
Z_{\beta_N,N}^{A,(j)}(0,0, N,\star)^{2}
]\leq (\sigma_N^{2}R_N)^j$  (and similarly for the line-to-point partition function).
Again for $p\geq 2$ small enough, $c_p$ can be absorbed into $\widehat{\beta}$, whereas Lemma~\ref{lem_rw_tk_properties}(i) gives that the supremum converges to one. Overall, this yields that \eqref{eq_hypercontract_p2p} is uniformly bounded in $N$.\\

The estimate for $\ptp_{\beta_N,N}(0,0\mid N,z)$ follows now from (i). Recall that 
\begin{align*}
    \ptp_{\beta_N,N}(0,0\mid N,z)=
    e^{\beta_N\w_{N,z}-\lambda(\beta_N)} Z_{\beta_N,N}(0,0\mid N,z),
\end{align*}
where we included again the disorder at the endpoint separately. We then have
\begin{align*}
    \dE[
\ptp_{\beta_N,N}(0,0\mid N,z)^{p}
]
=
\dE[
(e^{\beta_N\w_{N,z}-\lambda(\beta_N)})^{p}
]
 \dE[
Z_{\beta_N,N}(0,0\mid N,z)^{p}
],
\end{align*}
using the independence property of the disorder $\w$. Choosing $p<2+\delta$, where $\delta=\delta(\widehat{\beta})$ from part (i), we know the second term is uniformly bounded. It remains to estimate 
$
    \dE[
(e^{\beta_N\w_{N,z}-\lambda(\beta_N)})^{p}
]
=
e^{\lambda(p\beta_N)-p\lambda(\beta_N)}.
$
By Taylor expansion we can write
$
\lambda(p\beta_N)-p\lambda(\beta_N)
\sim \frac{1}{2}p(p-1)\beta_N^2,
$
which yields boundedness of the exponential. 
This concludes the proof of (ii).
\end{proof}

%% file: sections/annealed.tex
\section{The annealed polymer measure}\label{sec_annealed}

After proving the factorisation of point-to-point partition functions, we can finally start analysing the limiting polymer measure. 
As dealing with the annealed polymer measures first will substantially simplify the required steps in the quenched case,
we define the disorder-averaged measure $\mu_{\beta_N,N}$ on $(C[0,1],\mathcal{F})$,  by 
\begin{align}
    \mu_{\beta_N,N}(B):=
    \dE[\pi_N^{*}\pP_{\beta_N,N}^{\w}(B)]
    \qquad \forall B\in\mathcal{F}.
\end{align}
This section's main result is an invariance principle for the paths of the annealed polymer measure. 

\begin{proposition}[Annealed invariance principle]\label{prop_annealed_clt}
For $\widehat{\beta}\in (0,1)$, we have 
\begin{align*}
    \mu_{\beta_N,N}\weakconv \pP(\tfrac{1}{\sqrt{2}}W\in \cdot\,), \quad \text{as }N\to\infty,
\end{align*}
where we recall that $\pP$ denotes the Wiener measure on $C[0,1]$.
\end{proposition}

We begin by showing that the limiting finite-dimensional distributions of $\mu_{\beta_N,N}$ agree with the ones of a Brownian motion with diffusion matrix $\tfrac{1}{\sqrt{2}}I_2$.
In the section's second part we prove the required tightness in $C[0,1]$. Together with the identification of finite-dimensional distributions, this yields Proposition~\ref{prop_annealed_clt}.\\

Instead of determining the finite-dimensional distributions of the interpolated paths, it suffices to work with the corresponding starting point of the interpolation. To see this, recall that under $\pi_N^{*}\pP_{\beta_N,N}^{\w}$ we have 
\begin{align}\label{eq_def_interpolated_paths}
    X_{t}= \frac{1}{\sqrt{N}} S_{\lfloor tN \rfloor} +\frac{1}{\sqrt{N}} (tN-\lfloor tN \rfloor)(S_{\lfloor tN \rfloor+1}-  S_{\lfloor tN \rfloor} ).
\end{align}
Because the simple random walk has a finite range transtition kernel, the second term vanishes for $N$ large and we are left with $\tfrac{1}{\sqrt{N}} S_{\lfloor tN \rfloor}$.
Thus, the weak limits of $\tfrac{1}{\sqrt{N}}S_{\lfloor tN \rfloor}$ and $X_t$ (under $\pi_N^{*}\pP_{\beta_N,N}^{\w}$) must coincide if they exist. We will assume $tN\in \N$ for the sake of notation. \\

Let $0< t_1<t_2<\ldots<t_k\leq 1$. Throughout the remaining sections it will be convenient to fix a partitioning of the time intervals $(t_{i-1},t_i]$ and $(t_k,1]$ using  $t_i^{\pm}$'s such that
\begin{align}\label{eq_partitioning_of_time}
    t_{i-1}^{+}<t_i^-<t_i<t_i^+.
\end{align}
Now, recalling the representation in \eqref{eq_disc_fdd_as_p2p_product}, the idea when proving convergence of the marginal distributions is to replace the first point-to-point partition function using Proposition~\ref{prop_fact_part_func}:
\begin{align}\label{eq_annealed_ratio_approx_preproof}
    \frac{\ptp_{\beta_N,N}(0,0\mid t_1N,z_1)}{Z_{\beta_N,N}(0,0,N,\star)}
    \simeq 
    \frac{Z_{\beta_N,N}(0,0,t_1N,\star)Z_{\beta_N,N}(0,\star,t_1N,z_1)}{Z_{\beta_N,N}(0,0,N,\star)},
\end{align}
where `$\simeq$' should be understood as approximation in $L^1(\dP)$ for large $N$. Because both the denominator and the first term in the numerator only depend on the disorder in a small neighbourhood around the starting point, they converge to the same limit and should cancel as $N$ diverges.
We are then left with the remaining point-to-point partition functions which can be analysed separately due to independence of the disorder on disjoint time intervals.

Having the above approach in mind, we introduce a shorthand notation for the following constellation of terms.
Consider $0< t_1<t_2<\ldots<t_k< 1$, we define
\begin{align}\label{eq_def_Q}
&Q_{\beta_N,N}^\w \big((t_i,B_i)_{i=1}^k\big)\nonumber\\
    &\qquad:=
    \sum_{\substack{z_i\in \sqrt{N}B_i\\1\leq i\leq k}}
    Z_{\beta_N,N}(t_1^-N,\star, t_1N,z_1)
    \Big(
    \prod_{j=1}^{k-1}
    Z_{\beta_N,N}(t_jN,z_j, t_j^+N,\star)
    Z_{\beta_N,N}(t_{j+1}^-N,\star, t_{j+1}N,z_{j+1})
    \Big)\\
    &\qquad \qquad\qquad \times
    Z_{\beta_N,N}(t_kN,z_k, t_k^+N,\star)
    \prod_{j=1}^{k}q_{(t_j-t_{j-1})N}(z_{j}-z_{j-1}).\nonumber
\end{align}
Whenever $t_k=1$, the partition function  $Z_{\beta_N,N}(t_kN,z_k, t_k^+N,\star)$ is dropped, since it depends on the disorder outside of $\{0,\ldots, N\}\times \Z^2$.

\begin{lemma}\label{lem_annealed_law}
For $\widehat{\beta}\in (0,1)$ and $0< t_1<t_2<\ldots<t_k\leq 1$ we have 
\begin{align}\label{eq_L1_approx_quenched}
    \lim_{N\to\infty}\Big\|
    \pP_{\beta_N,N}^{\w}\big(\tfrac{1}{\sqrt{N}}S_{t_1N}\in  B_1,\ldots,\tfrac{1}{\sqrt{N}}S_{t_kN}\in  B_k\big)
    - 
    Q_{\beta_N,N}^\w \big((t_i,B_i)_{i=1}^k\big)
    \Big\|_1 = 0,
\end{align}
for every choice of bounded measurable sets $(B_i)_{i=1}^k\in (\R^2)^{k}$ satisfying 
$\lambda(\partial B_i)=0$.

In particular, this implies
\begin{align}\label{eq_annealed_fdd_goal}
\lim_{N\to\infty}
\dE
\big[
\pP_{\beta_N,N}^{\w}\big(\tfrac{1}{\sqrt{N}}S_{t_1N}\in  B_1,\ldots,\tfrac{1}{\sqrt{N}}S_{t_kN}\in  B_k\big)
\big]= \pP(\tfrac{1}{\sqrt{2}}W_{t_1}\in  B_1,\ldots,\tfrac{1}{\sqrt{2}}W_{t_k}\in  B_k).
\end{align}
\end{lemma}

\begin{proof}
We only elaborate the steps for the case $k=2$ and $t_2=1$; the general statement follows along the same lines modulo more involved notation. Let $(B_i)_{1\leq i\leq 2}$ be bounded sets with boundary of Lebesgue-measure zero.
Throughout the proof, we assume $B(x_i,r_i)$'s to be balls large enough such that they cover the bounded sets $B_i$.
We write $t_0=0$ and $z_0=0$ for the starting point of the random polymer.\\

First, we will replace the point-to-point partition function  $\ptp_{\beta_N,N}(0,0\mid t_1N,z_1)$ in \eqref{eq_disc_fdd_as_p2p_product} by its point-to-plane counterparts, using Proposition~\ref{prop_fact_part_func}. Afterwards, we exchange the arising term ${Z_{\beta_N,N}(0,0, t_0^+N,\star)}$  with ${Z_{\beta_N,N}(0,0, N,\star)}$ to cancel the point-to-plane partition function in the denominator. 
In other words, we want to show that 
\begin{align}\label{eq_too_tired}
\sup_{\substack{z_i\in \sqrt{N}B(x_i,r_i)\\ \text{s.t. } q_{t_1 N}(z_1)>0}}
    &\bigg\|
\frac{\ptp_{\beta_N,N}(t_{1}N,z_{1}\mid N,z_2)}{Z_{\beta_N,N}(0,0,N,\star)}
\left(
\ptp_{\beta_N,N}(0,0\mid t_1N,z_1)
 -  {Z_{\beta_N,N}(0,0, N,\star)}{Z_{\beta_N,N}(t_1^-N,\star, t_1N,z_1)}\right)
\bigg\|_1
\end{align}
vanishes in the large $N$ limits.
Afterwards, the remaining point-to-point partition function can be replaced using again Proposition~\ref{prop_fact_part_func}, which yields \eqref{eq_L1_approx_quenched}.

Instead of showing the above convergence in \eqref{eq_too_tired} directly,  we divide the statement into two, more manageable, terms by introducing the intermediate term ${Z_{\beta_N,N}(0,0, t_0^+N,\star)}{Z_{\beta_N,N}(t_1^-N,\star, t_1N,z_1)}$.
\begin{itemize}
    \item We start with the replacement of the point-to-point partition function following Proposition~\ref{prop_fact_part_func}. First, we apply the Cauchy-Schwartz inequality which yields 
    \begin{align}\label{eq_anneal_first_est}
        \qquad \qquad &\bigg\|
        \frac{\ptp_{\beta_N,N}(t_{1}N,z_{1}\mid N,z_2)}{Z_{\beta_N,N}(0,0,N,\star)}
        \left(
        \ptp_{\beta_N,N}(0,0\mid t_1N,z_1)
         -  {Z_{\beta_N,N}(0,0, t_0^+N,\star)}{Z_{\beta_N,N}(t_1^-N,\star, t_1N,z_1)}\right)
        \bigg\|_1 \nonumber \\
        &\leq 
        \|
        Z_{\beta_N,N}(0,0,N,\star)^{-1}
        \|_{2}
        \|
        \ptp_{\beta_N,N}(t_{1}N,z_{1}\mid N,z_2)
        \|_{2}\\
        &\qquad\qquad\times
        \|
        \ptp_{\beta_N,N}(0,0\mid t_1N,z_1)
         -  {Z_{\beta_N,N}(0,0, t_0^+N,\star)}{Z_{\beta_N,N}(t_1^-N,\star, t_1N,z_1)}
        \|_{2}.\nonumber
    \end{align}
   Here we used independence of the disorder on the disjoint time intervals $(0,t_1N]$ and $(t_1N,N]$ in addition to the fact
     that partition functions $\ptp_{\beta_N,N}(sN,0\mid tN,z)$ only depend on the disorder $\w$ in $(sN,tN]\times \Z^2$.
     The last term in \eqref{eq_anneal_first_est} vanishes uniformly in $z\in \sqrt{N}B(x_1,r_1)$ due to Proposition~\ref{prop_fact_part_func}, whereas the first and second term are uniformly bounded, cf. Lemma~\ref{lem_hypcont_p2l} and Corollary~\ref{cor_hypercontract}, respectively.
    
    \item Lastly, we can estimate the remaining norm using again the fact that we have arbitrary good control of negative moments of $Z_{\beta_N,N}(0,0,N,\star)$. 
    Using 
    \begin{align*}
        {Z_{\beta_N,N}(0,0, t_0^+N,\star)}
         -  {Z_{\beta_N,N}(0,0, N,\star)}
         &= {\widehat{Z}_{\beta_N,N}(0,0, t_0^+N,\star)}
         -  {\widehat{Z}_{\beta_N,N}(0,0, N,\star)},
    \end{align*}
    the Cauchy-Schwarz inequality gives 
    \begin{align*}
        &\left\|
        \frac{
        {Z_{\beta_N,N}(0,0, t_0^+N,\star)}
         -  {Z_{\beta_N,N}(0,0, N,\star)}}{Z_{\beta_N,N}(0,0,N,\star)}
        {Z_{\beta_N,N}(t_1^-N,\star, t_1N,z_1)}\ptp_{\beta_N,N}(t_{1}N,z_{1}\mid N,z_2)
        \right\|_1\\
        & \leq 
        \left\|\frac{
        {\widehat{Z}_{\beta_N,N}(0,0, t_0^+N,\star)}
         -  {\widehat{Z}_{\beta_N,N}(0,0, N,\star)}}{Z_{\beta_N,N}(0,0,N,\star)}\right\|_2
        \left\|{Z_{\beta_N,N}(t_1^-N,\star, t_1N,z_1)}\ptp_{\beta_N,N}(t_{1}N,z_{1}\mid N,z_2)\right\|_2.
    \end{align*}
    The second term on the r.h.s. is uniformly bounded by Lemma~\ref{lem_hypcont_p2l}(ii). The first term, on the other hand, can be estimated using H\"older's inequality and the trivial estimate $\|a-b\|\leq \|a\|+\|b\|$:
    \begin{align}\label{eq_push_neg2pos_moments}
        &\left\|\frac{
        \widehat{Z}_{\beta_N,N}(0,0, t_0^+N,\star)
         -  {\widehat{Z}_{\beta_N,N}(0,0, N,\star)}}{Z_{\beta_N,N}(0,0,N,\star)}\right\|_2\\
        &\leq \|
        Z_{\beta_N,N}(0,0,N,\star)^{-1}
        \|_{2+\delta^{-1}}
        (
        \|
        \widehat{Z}_{\beta_N,N}(0,0, t_0^+N,\star)
        \|_{2+\delta}\,
        +\,
        \|
        \widehat{Z}_{\beta_N,N}(0,0, N,\star)
        \|_{2+\delta}
        ),\nonumber
    \end{align}
    where $\delta>0$ is choosen sufficently small. Because $\|
        \widehat{Z}_{\beta_N,N}(0,0, N,\star)
        \|_{2+\delta}$ and $\|
        \widehat{Z}_{\beta_N,N}(0,0, t_0^+N,\star)
        \|_{2+\delta}$ converge to zero using Lemma~\ref{lem_hypcont_p2l}(i), while the first term is uniformly bounded, this proves convergence to zero of the expression in \eqref{eq_too_tired}. 
\end{itemize}
Having proven that \eqref{eq_too_tired} vanishes, it is only left to replace $\ptp_{\beta_N,N}(t_{1}N,z_{1}\mid N,z_2)$ by $Z_{\beta_N,N}(t_{1}N,z_{1}, t_{1}^+ N,\star)Z_{\beta_N,N}(t_{2}^-N,\star, N,z_2)$
which holds again by means of Proposition~\ref{prop_fact_part_func}. This concludes the proof of \eqref{eq_L1_approx_quenched}.

Overall, \eqref{eq_L1_approx_quenched} implies 
\begin{align*}
\lim_{N\to\infty}
    \dE
[
\pP_{\beta_N,N}^{\w}(S_{t_1N}\in \sqrt{N}B_1, S_{N}\in \sqrt{N}B_2)
]
&=
\lim_{N\to\infty}
    \dE
[
Q_{\beta_N,N}^{\w}((t_i,B_i)_{i=1}^2)
]\\
&=\lim_{N\to\infty}
\sum_{\substack{z_i\in \sqrt{N}B_i\\1\leq i\leq 2}}
\prod_{j=1}^{2}q_{(t_j-t_{j-1})N}(z_{j}-z_{j-1}),
\end{align*}
where we used the fact that $\dE[
Z_{\beta_N,N}(t_1^-N,\star,t_1 N,z_1)
Z_{\beta_N,N}(t_{1}N,z_{1}, t_{1}^+ N,\star)Z_{\beta_N,N}(t_{2}^-N,\star, N,z_2)
]=1$. 
It is well known that the simple random walk in $d=2$ converges in law to Brownian motion with diffusion matrix $\tfrac{1}{\sqrt{2}}I_2$, see for example \cite[Theorem 7.6.1]{LL10}. 
\end{proof}

\if{false}{
\red{Note that Lemma~\ref{lem_annealed_law} is not enough to deduce weak convergence of finite-dimensional distributions, because the statement is restricted to bounded sets. However, assuming the finite-dimensional distribution's limit exists, Lemma~\ref{lem_annealed_law} allows us to uniquely identify it. 
Tightness of the annealed polymer measures, on the other hand, yields existence of weak accumulation points on $\mathcal{M}_1(C[0,1])$. This implies in particular accumulation points of finite-dimensional distributions, which will be enough to prove Proposition~\ref{prop_annealed_clt}.}
}\fi

Note that Lemma~\ref{lem_annealed_law} is enough to deduce weak convergence of finite-dimensional distributions. In order to lift the convergence result to the annealed polymer measures, we require tightness.

\begin{lemma}\label{lem_annealed_tightness}
The family of annealed polymer measures $(\mu_{\beta_N,N})_N$ is tight in $\mathcal{M}_1(C[0,1])$.
\end{lemma}

\begin{proof}
We prove the tightness of $(\mu_{\beta_N,N})_N$ by following \cite[Theorem 16.5]{Ka02} which states that it suffices to show that for every $\varepsilon>0$
\begin{align}\label{eq_mod_cont_tight_annealed}
\lim_{\delta\to 0}\limsup_{N\to\infty}
\mu_{\beta_N,N}( m_{\delta}(X)  \geq \varepsilon)
=
\lim_{\delta\to 0}\limsup_{N\to\infty}
\dE[
\pi_{N}^{*}\pP_{\beta_{N},{N}}^{\w}(
m_{\delta}(X)\geq \varepsilon 
)]
=0,
\end{align}
where $m_{\delta}(\varphi)$ denotes the modulus of continuity on the Wiener space, i.e. 
\begin{align*}
m_{\delta}(\varphi):= \sup_{\substack{0\leq s,t \leq 1 \\ |t-s|<\delta}} \left| 
\varphi_t-\varphi_s
\right|
\qquad\forall\; \varphi \in C[0,1].
\end{align*}
Next, we use the same trick as stated in 
\cite[Lemma 4.2]{CY06} for random polymers in $d\geq 3$. See also 
\cite[Lemma 4.2]{ALKQ13} for a similar application in the case of the continuum random polymer in $d=1$. For every $\pP_N$-integrable $Y$, we have
\begin{align*}
\dE[
Z_{\beta_N,N}(0,0,N,\star) \pE^{\w}_{\beta_N,N}[
Y]
]
=
\pE_{N}[Y].
\end{align*}
Thus, for arbitrary $\varepsilon>0$ we can write
\begin{align}\label{eq_mod_cont_srw}
\dE[
Z_{\beta_N,N}(0,0,N,\star) \,
\pi_N^{*}\pP_{\beta_N,N}^{\w}(
m_{\delta}(X)\geq \varepsilon 
)]
=
\pi_N^{*}\pP_{N}(
m_{\delta}(X)\geq \varepsilon
).
\end{align}
However, the sequence of rescaled, interpolated simple random walks is known to be tight in the Wiener space. Hence, the right hand side of equation \eqref{eq_mod_cont_srw} vanishes when taking first the limit superior $N\to\infty $ and then the limit $\delta\to 0$.
In particular, this implies for any $\lambda>0$ that 
\begin{align*}
\lambda
\dE[
\mathds{1}\{
Z_{\beta_N,N}(0,0,N,\star) \geq \lambda\} \,
\pi_N^{*}\pP_{\beta_N,N}^{\w} (
m_{\delta}(X)\geq \varepsilon 
)
]\rightarrow 0, \qquad \text{as }N\to\infty,\; \delta\to 0. 
\end{align*}
On the other hand, 
\begin{align*}
\dE[
\mathds{1}\{
Z_{\beta_N,N}(0,0,N,\star) < \lambda\}\,
\pi_N^{*}\pP_{\beta_N,N}^{\w}(
m_{\delta}(X)\geq \varepsilon 
)]
&\leq \lambda \dE[
Z_{\beta_N,N}(0,0,N,\star)^{-1}
],
\end{align*}
where we first dropped the inner probability  before applying Markov's inequality. Once more we make use of the fact that $\dE[
Z_{\beta_N,N}(0,0,N,\star)^{-1}
]$  is uniformly bounded in $N$ which overall yields
\begin{align*}
    \mu_{\beta_N,N}
(
m_{\delta}(X)\geq \varepsilon 
)
&=
\dE[
\pi_N^{*}\pP_{\beta_N,N}^{\w}(
m_{\delta}(X)\geq \varepsilon 
)
]\\
&\leq 
\dE[
\mathds{1}\{
Z_{\beta_N,N}(0,0,N,\star) \geq \lambda\}
\pi_N^{*}\pP_{\beta_N,N}^{\w}(
m_{\delta}(X)\geq \varepsilon 
)
]
+
\lambda \widehat{C},
\end{align*}
for some $\widehat{C}= \widehat{C}_{\widehat{\beta}}>0$. 
Taking first the limits $N\to\infty$, $\delta\to 0$ and lastly the limit $\lambda\to 0$ finally yields \eqref{eq_mod_cont_tight_annealed}.
Thus, tightness of $(\mu_{\beta_N,N})_N$ follows.
\end{proof}

Combining now both Lemma~\ref{lem_annealed_law} and Lemma~\ref{lem_annealed_tightness}, we can show that annealed polymer paths converge in distribution to the ones of Brownian motion.

\begin{proof}[Proof of Proposition~\ref{prop_annealed_clt}]

We first note that Lemma~\ref{lem_annealed_law} implies convergence of the finite-dimensional marginals of the annealed polymer measure. 
Let $0< t_1<\ldots<t_k\leq 1$ be arbitrary.
Tightness of the marginals follows from the fact that, for every $\varepsilon>0$, we can choose bounded  continuity sets $B_i\subset \R^2$ such that 
$$
\pP(\tfrac{1}{\sqrt{2}}W_{t_1}\in B_1, \ldots , \tfrac{1}{\sqrt{2}}W_{t_k}\in B_k)> 1-\varepsilon\,.
$$
Thus, using Lemma~\ref{lem_annealed_law}, the corresponding polymer marginals satisfy 
\begin{align*}
   \inf_{N>\overline{N}} \dE
\big[
\pP_{\beta_N,N}^{\w}\big(\tfrac{1}{\sqrt{N}}S_{\lfloor t_1N \rfloor}\in  B_1,\ldots,\tfrac{1}{\sqrt{N}}S_{\lfloor t_kN\rfloor }\in  B_k\big)
\big]>1-\varepsilon\,,
\end{align*}
for some $\overline{N}=\overline{N}_\varepsilon\in \N$. This can be lifted to all $N\in\N$ by extending the sets $B_i$. 
Therefore, the laws of
$\tfrac{1}{\sqrt{N}}(S_{\lfloor t_1N \rfloor},\ldots,S_{\lfloor t_kN\rfloor })$ are tight and converge weakly to the law of 
$\tfrac{1}{\sqrt{2}}(W_{t_1},\ldots, W_{t_k})$, since bounded continuity sets suffice to identify the limiting measure uniquely. 
The same holds for $(X_{t_1},\ldots, X_{t_k})$, see \eqref{eq_def_interpolated_paths} and the discussion below, i.e.
\begin{align*}
    \mu_{\beta_{N},N}(X_{t_1}\in \cdot, \ldots , X_{t_k}\in \cdot) \weakconv \pP(\tfrac{1}{\sqrt{2}}W_{t_1}\in \cdot, \ldots , \tfrac{1}{\sqrt{2}}W_{t_k}\in \cdot).
\end{align*}
Recall now from Lemma~\ref{lem_annealed_tightness} that the annealed polymer measures $(\mu_{\beta_N,N})_N$ are tight. Hence, convergence of the annealed polymer measures is immediate, because weak accumulation points of $(\mu_{\beta_N,N})_N$ are determined by their finite-dimensional distributions, see for example \cite[Theorem 7.1]{Bi99}.
\end{proof}

%% file: sections/quenched.tex
\section{The quenched polymer measure}\label{sec_quenched}

After proving the invariance principle for the annealed polymer measures in the previous section, we can now proceed to prove our main results for the quenched polymer measures. 
In the first part of the section, we show convergence of finite-dimensional distributions of the rescaled, interpolated polymer paths, cf. Theorem~\ref{theo_fdd_quenched}. 
Afterwards, we prove Proposition~\ref{prop_functional_clt} which is a functional central limit theorem for the polymer paths when tested against individual test functions in $C_b(C[0,1])$.
We then explain how the functional central limit theorem can be lifted to an invariance principle, cf. Proposition~\ref{prop_equiv_fclt_ivp}, using a countable convergence determining class of functions.
Lastly, we prove the local limit theorem, Proposition~\ref{prop_polymer_llt}, for the polymer marginals on microscopic scales.

\subsection{Finite-dimensional distributions}

\subsubsection{Self-averaging of the random polymer}

We begin by showing convergence of the quenched polymer marginals in $L^1(\dP)$, evaluated on fixed bounded continuity sets $B$, that can be factorised:

\begin{proposition}\label{prop_fdd_quenched}
For $\widehat{\beta}\in (0,1)$ and $0< t_1<\ldots<t_k\leq 1$ we have 
\begin{align}\label{eq_prop_fdd_quenched_statement}
\big\|
\pP_{\beta_N,N}^{\w}\big(\tfrac{1}{\sqrt{N}}S_{t_1N}\in B_1,&\ldots, \tfrac{1}{\sqrt{N}}S_{t_kN}\in B_k\big)
- 
\pP(\tfrac{1}{\sqrt{2}}W_{t_1}\in  B_1,\ldots,\tfrac{1}{\sqrt{2}}W_{t_k}\in  B_k)\big\|_1\to 0,
\end{align}
for any choice of bounded measurable sets $(B_i)_{i=1}^k\in (\R^2)^{k}$ satisfying 
$\lambda(\partial B_i)=0$. 
\end{proposition}


We recall from Lemma~\ref{lem_annealed_law} that marginals of $\pi_N^{*}\pP_{\beta_N,N}^{\w}$ can be approximated in terms of $Q_{\beta_N,N}^{\w}$.
Thus, it is enough to show convergence of $Q_{\beta_N,N}^{\w}((t_i, B_i)_{i=1}^k)$ to the corresponding transition probabilities of Brownian motion, in order to lift Lemma~\ref{lem_annealed_law} to the quenched marginals. 

In the remainder of this section, we restrict ourselves to the case of $k=2$ with $t_2=1$.
The general statement of Proposition~\ref{prop_fdd_quenched} follows along the same lines.

\begin{lemma}\label{lem_upper_bnd_limit_fdd}
Let $t_1\in (0,1)$ and $B_1,B_2\subset \R^2$ be arbitrary bounded continuity sets, then 
\begin{align*}
    \lim_{N\to\infty}&
    \|Q_{\beta_N,N}^\w((t_i,B_i)_{i=1}^2)\|_2^2
= 
 \pP(\tfrac{1}{\sqrt{2}}W_{t_1}\in  B_1,\tfrac{1}{\sqrt{2}}W_{1}\in  B_2)^2.
\end{align*}
\end{lemma}

Recall the definition of $Q_{\beta_N,N}^\w$ in \eqref{eq_def_Q}. Independence of the occurring partition functions, yields that $\|Q_{\beta_N,N}^\w((t_i,B_i)_{i=1}^2)\|_2^2$ agrees with
\begin{align}\label{eq_uppbound_sec_mome_tran_prob}
   & \sum_{\substack{y_i\in \sqrt{N}B_i\\1\leq i\leq 2 }}
\sum_{\substack{z_i\in \sqrt{N}B_i\\1\leq i\leq 2 }}
\dE[
Z_{\beta_N,N}(t_1^-N,\star, t_1N,y_1)
{Z_{\beta_N,N}(t_1^-N,\star, t_1N,z_1)}
]\nonumber\\ 
&\qquad\qquad\times
  \dE[
Z_{\beta_N,N}(t_{1}N,y_{1}, t_1^+N,\star)
{Z_{\beta_N,N}(t_{1}N,z_{1}, t_1^+N,\star)}
]\\
&\qquad\qquad\times
\dE[
Z_{\beta_N,N}(t^-_2N,\star, N,y_2)
{Z_{\beta_N,N}(t^-_2N,\star, N,z_2)}
]\nonumber\\ 
&\qquad\qquad\times\prod_{j=1}^{2} q_{(t_j-t_{j-1})N}(y_j-y_{j-1})q_{(t_j-t_{j-1})N}(z_j-z_{j-1})\,. \nonumber
\end{align}
Hence, in order to get a sharp bound, we need precise estimates of mixed moments of the form 
\begin{align*}
    \dE[
        Z_{\beta_N,N}(sN,y,tN,\star)
        Z_{\beta_N,N}(sN,z,tN,\star)
        ]\,,
\end{align*}
locally uniformly in space, as $N$ tends to infinity.

{Recall from \eqref{eq_csz_weakconv_of_p2l} that partition functions starting at macroscopically separated points are asymptotically independent.} Therefore, one expects a law-of-large-number-like behaviour 
when ignoring the contribution of starting points $y_i$ and $z_i$ lying `too close' to each other. This idea encourages us to divide the points $(z_i)_{1\leq i\leq 2}$ into two groups, while keeping $(y_i)_{1\leq i\leq 2}$  fixed. 
We either have $|z_i-y_i|>2N^{1/2-a_N/4}$ for every $1\leq i \leq 2$
 or there exists an index $1\leq j\leq 2$ such that $|z_j-y_j|\leq 2N^{1/2-a_N/4}$.\\
 
Before stating the full proof of Lemma~\ref{lem_upper_bnd_limit_fdd}, we prove the following covariance estimate of two point-to-line partition functions starting from separate points of at least distance $2N^{\frac{1}{2}-\frac{a_N}{4}}$.

\begin{lemma}\label{lem_self_avging}
Let $s,t\in [0,1]$ such that $s<t$, $x\in\R^2$ and $r>0$, then 
\begin{align*}
     \lim_{N\to\infty}
\sup_{\substack{y,z\in \sqrt{N}B(x,r)\\ |y-z|> 2N^{\frac{1}{2}-\frac{a_N}{4}}}}
        \dE[
        Z_{\beta_N,N}(sN,y,tN,\star)
        Z_{\beta_N,N}(sN,z,tN,\star)
        ]
    = 1.
\end{align*}
The same statement also holds for line-to-point partition functions.
\end{lemma}
\begin{proof}
Let $y,z\in \sqrt{N}B(x,r)$ such that $|y-z|> 2N^{\frac{1}{2}-\frac{a_N}{4}}$. First, we expand the expectation in the statement 
\begin{align*}
    &\dE[
        Z_{\beta_N,N}(sN,y,tN,\star)
        Z_{\beta_N,N}(sN,z,tN,\star)
        ]\\
        &\qquad=
        \dE[
        Z_{\beta_N,N}^A(sN,y,tN,\star)
        Z_{\beta_N,N}^A(sN,z,tN,\star)
        ]
        +
        \dE[
        \widehat{Z}_{\beta_N,N}(sN,y,tN,\star)
        Z_{\beta_N,N}^A(sN,z,tN,\star)
        ]\\
        &\quad \qquad +
        \dE[
        Z_{\beta_N,N}^A(sN,y,tN,\star)
        \widehat{Z}_{\beta_N,N}(sN,z,tN,\star)
        ]
        +
        \dE[
        \widehat{Z}_{\beta_N,N}(sN,y,tN,\star)
        \widehat{Z}_{\beta_N,N}(sN,z,tN,\star)
        ].
\end{align*}
Using the Cauchy-Schwarz inequality and Lemma~\ref{lem_hypcont_p2l}, the last three terms on the r.h.s. vanish uniformly over $y,z\in \sqrt{N}B(x,r)$ in the large $N$ limit.
For the remaining term we use the fact that $|y-z|> 2N^{\frac{1}{2}-\frac{a_N}{4}}$ which implies $A^+_N(sN,y)\cap A^+_N(sN,z)=\emptyset $, whence
\begin{align*}
    \dE[
        Z_{\beta_N,N}^A(sN,y,tN,\star)
        Z_{\beta_N,N}^A(sN,z,tN,\star)
        ]
        =
        \dE[
        Z_{\beta_N,N}^A(sN,y,tN,\star)
        ]
        \dE[
        Z_{\beta_N,N}^A(sN,z,tN,\star)
        ]=1.
\end{align*}
This concludes the proof.
\end{proof}

\begin{proof}[Proof of Lemma \ref{lem_upper_bnd_limit_fdd}]
As in the proof of Lemma~\ref{lem_annealed_law}, we  assume that the sets $B_i$ are covered by open balls $B(x_i,r_i)$. 
Recall that the second moment of $Q_{\beta_N,N}^{\w}((t_i,B_i)_{i=1}^2)$ is given by the expression \eqref{eq_uppbound_sec_mome_tran_prob}.

First we show that points $y_i$ and $z_i$ lying `too close' to each other are negligible. Thereafter we prove the partition functions' self-averaging effect for the remaining points. Overall, this yields a sharp enough estimate for \eqref{eq_uppbound_sec_mome_tran_prob}.
Note that in the case where $B_i$'s have empty interior, self-averaging does not take place necessarily. However, in this case we can follow the same steps as in the first bullet below to deduce that \eqref{eq_prop_fdd_quenched_statement} still holds since both the random walk and polymer marginals evaluated on such sets converge to zero.

\begin{itemize}
    \item 
     We begin by analysing points lying `too close', i.e. at least one of the $z_i$'s lies in a $2N^{1/2-a_N/4}$ neighbourhood around $y_i$.
Without loss of generality we assume that this is the case for $z_1$.
Note that all expectation arising in \eqref{eq_uppbound_sec_mome_tran_prob} are uniformly bounded in $N$ and $y_i,z_i\in \sqrt{N}B(x_i,r_i)$, say by a constant $C>0$. This follows by an application of the Cauchy-Schwarz inequality and Lemma~\ref{lem_hypcont_p2l}. Thus, when considering the r.h.s. of \eqref{eq_uppbound_sec_mome_tran_prob} with the second sum restricted to $z_1\in \sqrt{N}B_1$ satisfying $|y_1-z_1|\leq 2N^{1/2-a_N/4}$, we can upper bound the expression by
\begin{align*}
C^3&
\sum_{\substack{y_i\in \sqrt{N}B_i\\1\leq i\leq 2}}
\sum_{\substack{z_1\in \sqrt{N}B_1\\ |y_1-z_1|\leq 2N^{1/2-a_N/4}\\ z_2\in \sqrt{N}B_2}}
\prod_{j=1}^{2} q_{(t_j-t_{j-1})N}(y_j-y_{j-1})q_{(t_j-t_{j-1})N}(z_j-z_{j-1})\\
 &\leq C^3\sum_{\substack{y_1\in \sqrt{N}B_1}} q_{t_1N}(y_1)
\sum_{\substack{z_1\in \sqrt{N}B_1\\ |y_1-z_1|\leq 2N^{1/2-a_N/4}}}
 q_{t_1N}(z_1),
\end{align*}
where we dropped the probability kernels which don't depend on $z_1$ or $y_1$ in the last step.
    However, $N^{1/2-a_N/4}$ vanishes in the macroscopic limit and therefore the right hand side converges to zero as $N\to\infty$, because
    \begin{align*}
         \sup_{\substack{y_1\in \sqrt{N}B_1}} 
         \sum_{\substack{z_1\in \sqrt{N}B_1\\ |y_1-z_1|\leq 2N^{1/2-a_N/4}}}
         q_{t_1N}(z_1)
         &\simeq \sup_{\substack{y_1\in \sqrt{N}B_1}} 
         \frac{1}{N}\sum_{\substack{z_1\in \sqrt{N}B_1\\ |y_1-z_1|\leq 2N^{1/2-a_N/4}}}
         2\,p_{\frac{t_1}{2}}(\tfrac{z_1}{\sqrt{N}})
         \to 0  ,
    \end{align*}
    where we used the local limit theorem for the simple random walk.
    Consequently, whenever the space points lie `too close' to each other, in the sense that there exists an index $j$ such that $|z_j-y_j|<2N^{1/2-a_N/4}$, they do not contribute to the limiting marginal distribution of the polymer path.\\

    \item It only remains to estimate the second part of the decomposition of \eqref{eq_uppbound_sec_mome_tran_prob}, where we restrict the sum over $z_i$'s such that every $z_i$ has at least distance $2N^{1/2-a_N/4}$ from $y_i$, i.e. 
    \begin{align}\label{eq_pre_est_quenchend_far_away}
        &\sum_{\substack{y_i\in \sqrt{N}B_i\\1\leq i\leq 2 }}
        \sum_{\substack{z_i\in \sqrt{N}B_i\\ |y_i-z_i|>2N^{1/2-a_N/4} \\1\leq i\leq 2 }}
        \dE[
Z_{\beta_N,N}(t_1^-N,\star, t_1N,y_1)
{Z_{\beta_N,N}(t_1^-N,\star, t_1N,z_1)}
]\nonumber\\ 
&\qquad\qquad\times
  \dE[
Z_{\beta_N,N}(t_{1}N,y_{1}, t_1^+N,\star)
{Z_{\beta_N,N}(t_{1}N,z_{1}, t_1^+N,\star)}
]\\
&\qquad\qquad\times
\dE[
Z_{\beta_N,N}(t^-_2N,\star, N,y_2)
{Z_{\beta_N,N}(t^-_2N,\star, N,z_2)}
]\nonumber\\ 
        &\qquad\qquad\times\prod_{j=1}^{2} q_{(t_j-t_{j-1})N}(y_j-y_{j-1})q_{(t_j-t_{j-1})N}(z_j-z_{j-1}).\nonumber
    \end{align}
    Using Lemma~\ref{lem_self_avging}, each occurring expectation in the expression above converges to $1$, when taking the large $N$ limit.
    
    \end{itemize}
Coming back to the term we wanted to estimate originally, 
we can now write (after taking  $\lim_{N\to\infty}$ on both sides)
\begin{align*}
    \lim_{N\to\infty}&
    \|Q_{\beta_N,N}^{\w}((t_i,B_i)_{i=1}^2)\|_2^2
= 
 \pP(\tfrac{1}{\sqrt{2}}W_{t_1}\in  B_1,\tfrac{1}{\sqrt{2}}W_{1}\in  B_2)^2,
\end{align*}
where we started with  \eqref{eq_uppbound_sec_mome_tran_prob}, neglected the space-points lying `too close' to each other and proved convergence of the remaining sum \eqref{eq_pre_est_quenchend_far_away} using Lemma~\ref{lem_self_avging}.
\end{proof}

Having just proven convergence of the  second moment of $Q_{\beta_N,N}^{\w}((t_i,B_i)_{i=1}^2)$ to its squared mean, it is an immediate consequence that $Q_{\beta_N,N}^{\w}((t_i,B_i)_{i=1}^2)$ converges to its mean in $L^2(\dP)$:
\begin{corollary}\label{cor_Q_conv_in_L2}
Let $t_1\in(0,1)$ and $B_1,B_2\subset \R^2$ be arbitrary bounded continuity sets, then 
\begin{align*}
    \lim_{N\to\infty}\big\|
    Q_{\beta_N,N}^{\w}((t_i,B_i)_{i=1}^2)
    -
    \pP(\tfrac{1}{\sqrt{2}}W_{t_1}\in  B_1,\tfrac{1}{\sqrt{2}}W_{1}\in  B_2)
    \big\|_2=0.
\end{align*}
\end{corollary}


Finally, we can summarise the results of above lemmas in the proof of Proposition~\ref{prop_fdd_quenched}:

\begin{proof}[Proof of Proposition~\ref{prop_fdd_quenched}]

Application of the triangle inequality yields 
\begin{align*}
    \|\pP_{\beta_N,N}^{\w}&(S_{t_1N}\in \sqrt{N}B_1, S_{N}\in \sqrt{N}B_2)
    -
    \pP(\tfrac{1}{\sqrt{2}}W_{t_1}\in  B_1,\tfrac{1}{\sqrt{2}}W_{1}\in  B_2)
    \|_1\\
    &\leq 
    \|\pP_{\beta_N,N}^{\w}(S_{t_1N}\in \sqrt{N}B_1, S_{N}\in \sqrt{N}B_2)
    -
    Q_{\beta_N,N}^{\w}((t_i,B_i)_{i=1}^2)
    \|_1\\
    &\qquad\qquad+
    \|Q_{\beta_N,N}^{\w}((t_i,B_i)_{i=1}^2)
    -
    \pP(\tfrac{1}{\sqrt{2}}W_{t_1}\in  B_1,\tfrac{1}{\sqrt{2}}W_{1}\in  B_2)
    \|_1.
\end{align*}
The first term on the r.h.s. vanishes by Lemma~\ref{lem_annealed_law}, whereas the second term vanishes due to Corollary~\ref{cor_Q_conv_in_L2}.
\end{proof}

Using standard estimates, one can show that the transition probabilities of the interpolated polymer path are well approximated by the corner points of the discrete path. Thus, Proposition~\ref{prop_fdd_quenched} also holds for the rescaled and interpolated polymer paths under $\pi^{\ast}_N
\pP_{\beta_N,N}^{\w}$: 


\begin{corollary}\label{cor_conv_fdd_quenched}
For $\widehat{\beta}\in (0,1)$ and $0< t_1<\ldots<t_k\leq 1$ we have 
\begin{align*}
\lim_{N\to\infty}
\big\|\pi^{\ast}_N
\pP_{\beta_N,N}^{\w}\big(X_{t_1}\in B_1,&\ldots, X_{t_k}\in B_k\big)
- 
\pP(\tfrac{1}{\sqrt{2}}W_{t_1}\in  B_1,\ldots,\tfrac{1}{\sqrt{2}}W_{t_k}\in  B_k)\|_1= 0,
\end{align*}
for every choice of bounded measurable sets $(B_i)_{i=1}^k\in (\R^2)^{k}$ satisfying 
$\lambda(\partial B_i)=0$. 
\end{corollary}

\if{false}{
\begin{proof}
Let $0\leq t_1<\ldots<t_k\leq 1$ and $(B_i)_{i=1}^k\in (\R^2)^{k}$ with 
$\lambda(\partial B_i)=0$.
In order to prove the corollary, it suffices to show that 
\begin{align*}
    \big\|
\pP_{\beta_N,N}^{\w}\big(X_{t_1}^{(N)}\in B_1,\ldots, X_{t_k}^{(N)}\in B_k\big) - \pP_{\beta_N,N}^{\w}\big(\tfrac{1}{\sqrt{N}}S_{\lfloor t_1N\rfloor }\in B_1,\ldots, \tfrac{1}{\sqrt{N}}S_{\lfloor t_kN \rfloor}\in B_k\big) \big\|_2\to 0,
\end{align*}
due to Proposition~\ref{prop_fdd_quenched}. 
First note that 
\begin{align}\label{eq_approx_interpol_sigma_add}
    &\big\|
\pP_{\beta_N,N}^{\w}\big(X_{t_1}^{(N)}\in B_1,\ldots, X_{t_k}^{(N)}\in B_k\big) - \pP_{\beta_N,N}^{\w}\big(\tfrac{1}{\sqrt{N}}S_{\lfloor t_1N\rfloor }\in B_1,\ldots, \tfrac{1}{\sqrt{N}}S_{\lfloor t_kN \rfloor}\in B_k\big) \big\|_2\nonumber \\
&\qquad = \big\|
\pP_{\beta_N,N}^{\w}\big(X_{t_1}^{(N)}\in B_1,\ldots, X_{t_k}^{(N)}\in B_k, \, \exists \,i : \,\tfrac{1}{\sqrt{N}}S_{\lfloor t_iN\rfloor }\notin B_i \big) \nonumber \\
&\qquad \qquad- \pP_{\beta_N,N}^{\w}\big(\tfrac{1}{\sqrt{N}}S_{\lfloor t_1N\rfloor }\in B_1,\ldots, \tfrac{1}{\sqrt{N}}S_{\lfloor t_kN \rfloor}\in B_k, \, \exists \,i : \,X_{t_i}^{(N)}\notin B_i\big) \big\|_2\\
&\qquad \leq 
\sum_{i=1}^k 
\big\|
\pP_{\beta_N,N}^{\w}\big(X_{t_i}^{(N)}\in B_i, \tfrac{1}{\sqrt{N}}S_{\lfloor t_iN\rfloor }\notin B_i \big)\big\|_2 
+
\big\|\pP_{\beta_N,N}^{\w}\big(\tfrac{1}{\sqrt{N}}S_{\lfloor t_iN\rfloor }\in B_i, \,X_{t_i}^{(N)}\notin B_i\big) \big\|_2.\nonumber
\end{align}
Because $|X_{t_i}^{(N)}-\tfrac{1}{\sqrt{N}}S_{\lfloor t_iN\rfloor}|\leq \tfrac{1}{\sqrt{N}}$ under $\pP_{\beta_N,N}^{\w}$, we know that both $X_{t_i}^{(N)}$ and $\tfrac{1}{\sqrt{N}}S_{\lfloor t_iN\rfloor}$ must be close to the boundary of $B_i$ whenever only one of them lies in $B_i$. More precisely, for an arbitrary $\varepsilon>0$ and all $N$ large enough we have
\begin{align}\label{eq_sets_inclusion_approx_interpo}
    \big\{X_{t_i}^{(N)}\in B_i, \tfrac{1}{\sqrt{N}}S_{\lfloor t_iN\rfloor }\notin B_i\big\}\cup \big\{\tfrac{1}{\sqrt{N}}S_{\lfloor t_iN\rfloor }\in B_i, \,X_{t_i}^{(N)}\notin B_i\big\}\subset 
    \big\{X_{t_i}^{(N)}, \tfrac{1}{\sqrt{N}}S_{\lfloor t_iN\rfloor } \in (\partial B)^{(\varepsilon)}\big\},
\end{align}
where $(\partial B_i)^{(\varepsilon)}$ denotes a set
$$
(\partial B_i)^{(\varepsilon)}:=
\bigcup_{j=1}^{m_{\varepsilon}}B\big(x_i^{(\varepsilon)},\varepsilon\big)\supset \partial B_i,
$$
for some $x_1^{(\varepsilon)},\ldots,x_{m_{\varepsilon}}^{(\varepsilon)}\in \partial B_i$.
Such a family of points exist due to compactness of $\partial B_i$.
Now,  
\begin{align*}
    \big\|\pP_{\beta_N,N}^{\w}\big(\tfrac{1}{\sqrt{N}}S_{\lfloor t_iN\rfloor }\in  (\partial B_i)^{(\varepsilon)}\big)\big\|_2
    &\leq \|\pP_{\beta_N,N}^{\w}\big(\tfrac{1}{\sqrt{N}}S_{\lfloor t_iN\rfloor }\in  (\partial B_i)^{(\varepsilon)}\big)- \pP\big(\tfrac{1}{\sqrt{2}}W_{t_i}\in (\partial B_i)^{(\varepsilon)} \big) \|_2\\
    &\qquad + \pP\big(\tfrac{1}{\sqrt{2}}W_{t_i}\in (\partial B_i)^{(\varepsilon)} \big).
\end{align*}
Because $(\partial B_i)^{(\varepsilon)}$ is an open set with $\lambda(\partial(\partial B_i)^{(\varepsilon)} )=0$,
we can apply Proposition~\ref{prop_fdd_quenched} and the first term on the r.h.s. vanishes as $N\to\infty$, which only leaves us to estimate the second term. Because $\varepsilon>0$ was chosen arbitrary, we can write
\begin{align}\label{eq_interp_close_boundary}
\lim_{\varepsilon\to 0}
\lim_{N\to\infty}
\big\|\pP_{\beta_N,N}^{\w}\big(\tfrac{1}{\sqrt{N}}S_{\lfloor t_iN\rfloor }\in  (\partial B_i)^{(\varepsilon)}\big)\big\|_2
\leq
\lim_{\varepsilon\to 0}
    \pP\big(\tfrac{1}{\sqrt{2}}W_{t_i}\in (\partial B_i)^{(\varepsilon)} \big) = \pP\big(\tfrac{1}{\sqrt{2}}W_{t_i}\in \partial B_i \big)=0,
\end{align}
since $B_i$ is a continuity set of the normal distribution by assumption.
Combining now \eqref{eq_approx_interpol_sigma_add}, \eqref{eq_sets_inclusion_approx_interpo} and \eqref{eq_interp_close_boundary},  we showed that the transition probabilities of the interpolated paths are well approximated (in $L^2(\dP)$) by the corner points of the discrete polymer path. This finishes the proof.
\end{proof}
}\fi

\subsubsection{Tightness and uniqueness of the limit}

In Corollary~\ref{cor_conv_fdd_quenched} we showed convergence of polymer marginals evaluated on bounded, factorised, continuity sets in $L^1(\dP)$. This can be lifted to unbounded measurable sets $U\subset(\R^2)^k$ satisfying the same properties. 
However, we want to stress that this does not imply weak convergence of the polymer marginals yet, since exceptional points of the disorder can depend on the choice of sets $U$.
Nevertheless, we are able to show weak convergence of quenched finite-dimensional distributions in probability, since probability measures on $(\R^2)^k$ are uniquely identified by evaluation on a countable family of sets.

We begin by recalling a standard result: a sequence of random variables on a metric space converges in probability if and only if every subsequence has a further subsequence which converges almost surely, see for example \cite[Lemma 4.2]{Ka02}.
Thus, in order to show the convergence in probability of the marginal distributions, it suffices to prove tightness along sufficiently many subsequences
and identify the limit points using a $\pi$-system. This step is motivated by the recent article \cite{Ju21}, where Junk showed convergence of the polymer-endpoint distribution in bond disorder for $d\geq 3$ using a very similar approach. 

For the $\pi$-system, on which we will identify the limiting finite-dimensional distributions, we choose half-open cylinders on $(\R^2)^k$:
\begin{align*}
\mathcal{E}^k
:=
\{
[a_1,b_1)\times \cdots \times [a_{2k},b_{2k})\subset (\R^2)^k\, :\, a_i,b_i \in \mathbb{Q} \text{ and } a_i< b_i \text{ for } 1\leq i \leq 2k
\},
\end{align*}
which generates the Borel-sigma-algebra on $(\R^2)^k$. Note that $\mathcal{E}^k$ has countably many elements and let $\{E_i\}_{i=1}^\infty$ be an arbitrary enumeration of them. We omitted the dependency of $E_i$'s on $k$ for the sake of a lighter notation.

We start by showing that there exist sufficiently many subsequences along which  finite-dimensional distributions evaluated on the $\mathcal{E}^k$ converge almost surely.

\begin{lemma}\label{lem_pi_system_fdd}
Let $0\leq t_1 <\ldots <t_k \leq 1$, then for every sequence $(N_j)_{j\in\N}$ in $\N$ there exists a subsequence $(N_{j_m})_{m\in\N}$ and $\Omega_{1} =\Omega_{1}((N_j)_{j\in\N},(t_i)_{1\leq i \leq k})\subset \Omega$ with $\dP(\Omega_{1})=1$ such that for every $\w \in \Omega_{1}$
\begin{align*}
\lim_{m\to\infty}
\pi_{N_{j_m}}^{*}\pP_{\beta_{N_{j_m}},{N_{j_m}}}^{\w}((X_{t_1}, \ldots, X_{t_k})\in E)
=
\pP(
\tfrac{1}{\sqrt{2}}(W_{t_1}, \ldots, W_{t_k})\in E)
\quad \forall E\in \mathcal{E}^k.
\end{align*}
\end{lemma}

\begin{proof}
Let $(N_j)_{j\in\N}$ be an arbitrary sequence in $\N$.
We prove the lemma only for a single marginal $t\in[0,1]$, the multi-marginal case follows along the same lines.
Corollary~\ref{cor_conv_fdd_quenched} implies that for every $\varepsilon>0$ and $E_i\in\mathcal{E}$
\begin{align*}
\lim_{j\to\infty}\dP(|
    \pi_{N_{j}}^{*}\pP_{\beta_{N_{j}},{N_{j}}}^{\w}(X_t\in E_i)
    -
    \pP(\tfrac{1}{\sqrt{2}}W_t\in E_i)|>\varepsilon)=0.
\end{align*}
In particular, for every $i,m\in\N$ there exists a $M_{i,m}\in\N$ such that 
\begin{align*}
    \dP(|
    \pi_{N_{j}}^{*}\pP_{\beta_{N_{j}},{N_{j}}}^{\w}(X_t\in E_i)
    -
    \pP(\tfrac{1}{\sqrt{2}}W_t\in E_i)|>\tfrac{1}{m})\leq m^{-1}2^{-m} \quad \forall j\geq M_{i,m}. 
\end{align*}
We define a subsequence of $(N_j)_j$ using $j_m:=j_{m-1}\vee \max_{1\leq i\leq m }M_{i,m}$, then for every $m\in\N$
\begin{align*}
    \dP(\exists i\leq m \text{ with }|
    \pi_{N_{j_m}}^{*}\pP_{\beta_{N_{j_m}},{N_{j_m}}}^{\w}(X_t\in E_i)
    -
    \pP(\tfrac{1}{\sqrt{2}}W_t\in E_i)|>\tfrac{1}{m})\leq 2^{-m},
\end{align*}
which is summable in $m$. The Borel-Cantelli lemma then yields  
\begin{align*}
    \dP\big( \sup_{i\leq m}|
    \pi_{N_{j_m}}^{*}\pP_{\beta_{N_{j_m}},{N_{j_m}}}^{\w}(X_t\in E_i)
    -
    \pP(\tfrac{1}{\sqrt{2}}W_t\in E_i)| >\tfrac{1}{m} \ \text{ infinitely often}\big)=0,
\end{align*}
which implies that
\begin{align*}
\Omega_1:=
     \big\{\w\in\Omega \,:\,\lim_{m\to\infty} \sup_{i\leq m}|
    \pi_{N_{j_m}}^{*}\pP_{\beta_{N_{j_m}},{N_{j_m}}}^{\w}(X_t\in E_i)
    -
    \pP(\tfrac{1}{\sqrt{2}}W_t\in E_i)|=0\big\}
\end{align*}
has full mass,
i.e. $\dP(\Omega_1)=1$.
This concludes the proof.
\end{proof}

\if{false}{
The previous lemma allows us to identify possible limits along subsequences. In order to show existence of such limits, we need to show tightness.
\red{
\begin{lemma}\label{lem_fdd_tightness}
Let $0\leq t_1 <\ldots <t_k \leq 1$, then for every sequence $(N_j)_{j\in \N}$ in $\N$ there exists a subsequence $(N_{j_l})_{l\in \N}$ and $\Omega_{2} =\Omega_{2}((N_j)_{j\in\N},(t_i)_{1\leq i \leq k})\subset \Omega$ with $\dP(\Omega_{2})=1$ such that 
\begin{align}
\Big(\pi_{N_{j_l}}^{*}\pP_{\beta_{N_{j_l}},{N_{j_l}}}^{\w}\big((X_{t_1}, \ldots, X_{t_k})\in \cdot)\big)\Big)_{l\in\N}
\end{align}
is tight for every $\w \in \Omega_{2}$.
\end{lemma}
\begin{proof}
Again we restrict the presentation to the case of a single marginal at $t\in[0,1]$. The general case follows along the same lines. \\
Let $(N_j)_{j\in\N}$ be a sequence in $\N$.
We begin by showing that the existence of subsequences $(N_{j_l})_{l\in\N}$, along which the family of marginal distributions is tight, is an implication of  
\begin{align}\label{eq_fdd_tight_alt_statement}
    \lim_{l\to\infty} \sup_{r\in\N} \big|
    \pi_{N_{j_l}}^{*}\pP_{\beta_{N_{j_l}},{N_{j_l}}}^{\w}(X_{t}\in K_r^c)
    -
    \pP(
\tfrac{1}{\sqrt{2}}W_{t}\in K_r^c)
    \big|=0 ,
\end{align}
where $K_r$ denotes the compact ball $\overline{B}(0,r)\subset\R^2$ of radius $r$.
First, note that for every $L\in\N$
\begin{align*}
    &\lim_{r\to\infty} \sup_{l\in \N} \big|
    \pi_{N_{j_l}}^{*}\pP_{\beta_{N_{j_l}},{N_{j_l}}}^{\w}(X_{t}\in K_r^c)
    -
    \pP(
\tfrac{1}{\sqrt{2}}W_{t}\in K_r^c)
    \big|\\
    &\qquad =
    \lim_{r\to\infty} \sup_{l\geq L} \big|
    \pi_{N_{j_l}}^{*}\pP_{\beta_{N_{j_l}},{N_{j_l}}}^{\w}(X_{t}\in K_r^c)
    -
    \pP(
\tfrac{1}{\sqrt{2}}W_{t}\in K_r^c)
    \big|\\
    &\qquad \leq
    \sup_{r\in\N} \sup_{l\geq L} \big|
    \pi_{N_{j_l}}^{*}\pP_{\beta_{N_{j_l}},{N_{j_l}}}^{\w}(X_{t}\in K_r^c)
    -
    \pP(
\tfrac{1}{\sqrt{2}}W_{t}\in K_r^c)
    \big|,
\end{align*}
where we used in the first equality that $\lim_{r\to\infty}|
    \pi_{N_{j_l}}^{*}\pP_{\beta_{N_{j_l}},{N_{j_l}}}^{\w}(X_{t}\in K_r^c)
    -
    \pP(
\tfrac{1}{\sqrt{2}}W_{t}\in K_r^c)
    |= 0$ for any fixed $l$. We can now exchange the order of suprema and take the limit $L\to\infty$ which yields 
\begin{align}\label{eq_conseq_exch_lim_and_sup}
    &\lim_{r\to\infty} \sup_{l\in \N} \big|
    \pi_{N_{j_l}}^{*}\pP_{\beta_{N_{j_l}},{N_{j_l}}}^{\w}(X_{t}\in K_r^c)
    -
    \pP(
\tfrac{1}{\sqrt{2}}W_{t}\in K_r^c)
    \big|\\
    &\qquad 
    \leq 
    \limsup_{L\to\infty} \sup_{r\in \N} \big|
    \pi_{N_{j_l}}^{*}\pP_{\beta_{N_{j_l}},{N_{j_l}}}^{\w}(X_{t}\in K_r^c)
    -
    \pP(
\tfrac{1}{\sqrt{2}}W_{t}\in K_r^c)
    \big|=0,\nonumber
\end{align}
because the limit on the right must agree with \eqref{eq_fdd_tight_alt_statement}. Now, 
\begin{align*}
    &\sup_{l\in\N}\pi_{N_{j_l}}^{*}\pP_{\beta_{N_{j_l}},{N_{j_l}}}^{\w}(X_{t}\in K_r^c)
     \leq 
    \sup_{l\in\N}\big|
    \pi_{N_{j_l}}^{*}\pP_{\beta_{N_{j_l}},{N_{j_l}}}^{\w}(X_{t}\in K_r^c)
    -
    \pP(
\tfrac{1}{\sqrt{2}}W_{t}\in K_r^c)
    \big|
    +
    \pP(
\tfrac{1}{\sqrt{2}}W_{t}\in K_r^c).
\end{align*}
Assuming that \eqref{eq_fdd_tight_alt_statement} holds, both terms on the right vanish when taking the limit $r\to\infty$ as a consequence of \eqref{eq_conseq_exch_lim_and_sup}. Thus, it suffices to prove \eqref{eq_fdd_tight_alt_statement} along subsequences to conclude the statement of the lemma. \\
The existence of subsequences $(N_{j_l})_{l\in\N}$ along which \eqref{eq_fdd_tight_alt_statement} holds is equivalent to showing convergence in probability, i.e. for every $\varepsilon>0$
\begin{align}\label{eq_alt_stmt_in_prob}
\lim_{N\to\infty}
    \dP\big(
    \sup_{r\in\N} \big|
    \pi_{N}^{*}\pP_{\beta_{N},{N}}^{\w}(X_{t}\in K_r^c)
    -
    \pP(
\tfrac{1}{\sqrt{2}}W_{t}\in K_r^c)
    \big|
    >\varepsilon
    \big)=0.
\end{align}
For any $R\in\N$ we have 
\begin{align}\label{eq_conv_prob_alt_statem}
&\lim_{N\to\infty}
    \dP\big(
    \sup_{r\in\N} \big|
    \pi_{N}^{*}\pP_{\beta_{N},{N}}^{\w}(X_{t}\in K_r^c)
    -
    \pP(
\tfrac{1}{\sqrt{2}}W_{t}\in K_r^c)
    \big|
    >\varepsilon
    \big)\\
    &\qquad\leq 
    \lim_{N\to\infty}
    \dP\big(
    \sup_{r\geq R} \big|
    \pi_{N}^{*}\pP_{\beta_{N},{N}}^{\w}(X_{t}\in K_r^c)
    -
    \pP(
\tfrac{1}{\sqrt{2}}W_{t}\in K_r^c)
    \big|
    >\varepsilon
    \big),\nonumber
\end{align}
where we used Corollary~\ref{cor_conv_fdd_quenched}, after adding and subtracting $1$ inside the absolute value to replace $K^c_r$ by $K_r$, to show that 
\begin{align*}
& \lim_{N\to\infty}\sum_{r=1}^{R-1}\dP(|\pi_{N}^{*}\pP_{\beta_{N},{N}}^{\w}(X_{t}\in K_r)-\pP(\tfrac{1}{\sqrt{2}}W_{t}\in K_r)|>\varepsilon)\\
    &\qquad =\lim_{N\to\infty}\sum_{r=1}^{R-1}\dP(|\pi_{N}^{*}\pP_{\beta_{N},{N}}^{\w}(X_{t}\in K_r^c)-\pP(\tfrac{1}{\sqrt{2}}W_{t}\in K_r^c)|>\varepsilon)=0.
\end{align*} 
We can now upper bound \eqref{eq_conv_prob_alt_statem} further, using
\begin{align}\label{eq_est_alt_statement_sec2fin}
    &
    \lim_{N\to\infty}
    \dP\big(
    \sup_{r\geq R} \big|
    \pi_{N}^{*}\pP_{\beta_{N},{N}}^{\w}(X_{t}\in K_r^c)
    -
    \pP(
\tfrac{1}{\sqrt{2}}W_{t}\in K_r^c)
    \big|
    >\varepsilon
    \big)\nonumber\\
    &\qquad\leq
    \lim_{N\to\infty}
    \dP\big(
    \pi_{N}^{*}\pP_{\beta_{N},{N}}^{\w}(X_{t}\in K_R^c)
    +
    \pP(
\tfrac{1}{\sqrt{2}}W_{t}\in K_R^c)
    >\varepsilon
    \big)\\
    &\qquad\leq
    \lim_{N\to\infty}
    \dP\big(
    \pi_{N}^{*}\pP_{\beta_{N},{N}}^{\w}(X_{t}\in K_R^c)
    >\tfrac{\varepsilon}{2}
    \big),\nonumber
\end{align}
for any $R$ large enough such that $\pP(\tfrac{1}{\sqrt{2}}W_t\in K_R^c)<\varepsilon/2$. Tightness of the annealed polymer measure from Lemma~\ref{lem_annealed_tightness} then yields
\begin{align*}
\lim_{N\to\infty}
    \dP\big(
    \pi_{N}^{*}\pP_{\beta_{N},{N}}^{\w}(X_{t}\in K_R^c)
    >\tfrac{\varepsilon}{2}
    \big)
    \leq \frac{2}{\varepsilon} \sup_{N\in\N}\dE[\pi_{N}^{*}\pP_{\beta_{N},{N}}^{\w}(X_{t}\in K_R^c)] \to 0, \quad \text{ as }R\to\infty.
\end{align*}
Together with \eqref{eq_conv_prob_alt_statem} and \eqref{eq_est_alt_statement_sec2fin}, we showed convergence in probability \eqref{eq_alt_stmt_in_prob} which implies the existence of a subsequence and a set $\Omega_2$ of full measure such that \eqref{eq_fdd_tight_alt_statement} holds for every $\w\in\Omega_2$. This finishes the proof.
\end{proof}
}
}\fi

\begin{proof}[Proof of Theorem~\ref{theo_fdd_quenched}]
We fix $0\leq t_1 <\ldots <t_k \leq 1$ and let $(N_{j})_{j\in\N}$ be a sequence in $\N$. In Lemma~\ref{lem_pi_system_fdd} we proved the existence of a subsequence $(N_{j_m})_{m\in\N}$ and disorders $\Omega_1$, with $\dP(\Omega_1)=1$, such that 
\begin{align}\label{eq_conv_pi_sys_statement}
\lim_{m\to\infty}
\pi_{N_{j_m}}^{*}\pP_{\beta_{N_{j_m}},{N_{j_m}}}^{\w}((X_{t_1}, \ldots, X_{t_k})\in E)
=
\pP(
\tfrac{1}{\sqrt{2}}(W_{t_1}, \ldots, W_{t_k})\in E)
\quad \forall E\in \mathcal{E}^k,
\end{align}
for every $\w\in\Omega_1$.
Tightness of the sequence $(\pi_{N_{j_m}}^{*}\pP_{\beta_{N_{j_m}},{N_{j_m}}}^{\w}((X_{t_1}, \ldots, X_{t_k})\in \cdot))_{m\in \N}$ is an immediate consequence of \eqref{eq_conv_pi_sys_statement}, cf. proof of Proposition~\ref{prop_annealed_clt}, and the limiting probability measure is uniquely determined by the $\pi$-system $\mathcal{E}^k$.

Overall, we showed that for every sequence $(N_{j})_{j\in\N}$ in $\N$ there exists a subsequence $(N_{j_{m}})_{m\in\N}$ along which the finite-dimensional distributions converge almost surely to the ones of Brownian motion with diffusion matrix $\tfrac{1}{\sqrt{2}}I_2$. This is equivalent to weak convergence of the polymer marginals in $\dP$-probability.
\end{proof}

\subsection{An invariance principle for polymer paths}

We are finally ready to prove Theorem~\ref{theo_main}. 
Using tightness of the annealed polymer paths, cf. Lemma~\ref{lem_annealed_tightness}, and convergence of the finite-dimensional distributions, cf. Theorem~\ref{theo_fdd_quenched}, we can prove the desired result.
The steps resemble very much the ones when proving the invariance principle for the simple random walk. However, due the double randomness of paths and the environment, cf. Remark~\ref{rem_flct_ivp}, we require an additional argument to conclude the full invariance principle.  

We begin by proving a functional central limit theorem:

\begin{proposition}\label{prop_functional_clt}
Let $\widehat{\beta}\in (0,1)$ and $\beta_N$ as in \eqref{eq_R_N}.
Then for every $F\in C_b(C[0,1])$
$$
\pi_N^{\ast}\pE^\w_{\beta_N,N}[F(X)]
\to \pE[F(\tfrac{1}{\sqrt{2}}W )], \qquad \text{as }N\to\infty,
\quad \text{in }\dP\text{-probability},
$$ 
where $\pE$ is the expectation w.r.t. the Wiener measure on $C[0,1]$.
\end{proposition}

\begin{remark}
Note that the convergence in the functional central limit theorem above also holds in $L^1(\dP)$, since the random variables $(\pi_N^{\ast}\pE^\w_{\beta_N,N}[F(X)])_N$ are uniformly bounded by $\|F\|_\infty$. In particular, this implies convergence of expectations which is equivalent to the annealed invariance principle, cf. Proposition~\ref{prop_annealed_clt}.
\end{remark}

\begin{proof}
Let $F\in C_b(C[0,1])$. The statement of the theorem is equivalent to 
\begin{align}\label{eq_conv_in_prob_path_tested}
    \lim_{N\to\infty}
    \dP(
    |
    \pi_N^{\ast}\pE^\w_{\beta_N,N}[F(X)]
- \pE[F(\tfrac{1}{\sqrt{2}}W )]
    |>\varepsilon)=0\quad \forall \varepsilon>0.
\end{align}
In the following, we fix $\varepsilon>0$ and choose $\delta>0$ arbitrary. Let then $K=K(\delta,\varepsilon,F)\subset C[0,1]$ be a compact set such that
\begin{align}\label{eq_choice_K_tight}
\pP(\tfrac{1}{\sqrt{2}}W\in K^c)\leq \frac{\varepsilon}{8 \|F\|_{\infty}}
\quad \text{ and }\quad 
\sup_{N\in\N}
\mu_{\beta_N,N}(K^c)=
\sup_{N\in\N}
\dE[
\pi_N^{*}\pP_{\beta_N,N}(K^c)]
\leq 
\frac{\delta\, \varepsilon}{8 \|F\|_{\infty}}.
\end{align}
A set $K$ with these properties exists due to tightness of the Wiener measure and tightness of the annealed polymer measure, see Lemma~\ref{lem_annealed_tightness}.

Throughout the remainder of the proof,  we will denote by $\Pi_{t_1,\ldots,t_k}: C[0,1]\mapsto (\R^2)^k$ the projection of a path onto previously chosen coordinates $0\leq t_1<\ldots<t_k\leq 1$, i.e. 
$$
\Pi_{t_1,\ldots,t_k}(X)= (X_{t_1},\ldots, X_{t_k})\in (\R^2)^k.
$$
Restricting the domain of $F$ to the compact set $K$, we can approximate $F$ uniformly by cylinder functions on $C(K,\R)$, i.e. functions that only depend on finitely many coordinates of the path, using the Stone-Weierstrass theorem \cite[Theorem 13.4]{La02}. More precisely, there exist $0\leq t_1<\ldots<t_k\leq 1$ and a continuous $f:\Pi_{t_1,\ldots,t_k}(K)  \mapsto \R$ such that $\|f\circ \Pi_{t_1,\ldots,t_k}\|_{K,\infty}\leq \|F\|_\infty$ and
\begin{align}\label{eq_est_approx_cyl}
    \|F-f\circ \Pi_{t_1,\ldots,t_k}\|_{K,\infty}
    :=
    \sup_{\varphi \in K} |F(\varphi)-f(\varphi_{t_1},\ldots, \varphi_{t_k}) |<\frac{\varepsilon}{4}.
\end{align}
Using the Tietze extension theorem \cite[Theorem 20.4]{Ru86}, we can extend $f$ continuously from $\Pi_{t_1,\ldots,t_k}(K) $ to $(\R^2)^k$. At the same time, this yields a continuous extension of $f\circ \Pi_{t_1,\ldots,t_k}$ on $C[0,1]$. Without loss of generality, the extension can be chosen in such a way that $\|f\circ\Pi_{t_1,\ldots,t_k}\|_{\infty}\leq \|f\circ\Pi_{t_1,\ldots,t_k}\|_{K,\infty} \leq  \|F\|_{\infty}$.
Estimating now $(F-f\circ\Pi_{t_1,\ldots,t_k})$ on $K$ and $K^c$ respectively, after applying the triangle inequality, yields 
\begin{align*}
    &|\pi_N^{\ast}\pE^\w_{\beta_N,N}[(\mathds{1}_K+\mathds{1}_{K^c})(F-f\circ\Pi_{t_1,\ldots,t_k} )(X)]
- \pE[(\mathds{1}_K+\mathds{1}_{K^c})(F-f\circ\Pi_{t_1,\ldots,t_k} )(\tfrac{1}{\sqrt{2}}W )]|\\
&\qquad\leq 
2\|F-f\circ\Pi_{t_1,\ldots,t_k}\|_{K,\infty}+ 2\|F\|_{\infty} 
|\pi_N^{\ast}\pP^\w_{\beta_N,N}(K^c)
+ \pP(\tfrac{1}{\sqrt{2}}W\in K^c) |.
\end{align*}
Thus, due to \eqref{eq_est_approx_cyl} and the choice of $K$ in \eqref{eq_choice_K_tight}, we have 
\begin{align}\label{eq_est_diff_f_and_F}
    &|\pi_N^{\ast}\pE^\w_{\beta_N,N}[(F-f\circ\Pi_{t_1,\ldots,t_k} )(X)]
- \pE[(F-f\circ\Pi_{t_1,\ldots,t_k} )(\tfrac{1}{\sqrt{2}}W )]|
\leq 
\frac{3\varepsilon}{4}+ 2 \|F\|_{\infty} \pi_N^{\ast}\pP^\w_{\beta_N,N}(K^c).
\end{align}
All together, we can upper bound the term of interest
\begin{align*}
    |\pi_N^{\ast}\pE^\w_{\beta_N,N}[F(X)]
- \pE[F(\tfrac{1}{\sqrt{2}}W )]|
&\leq 
|\pi_N^{\ast}\pE^\w_{\beta_N,N}[f(X_{t_1},\ldots, X_{t_k})]
- \pE[f(\tfrac{1}{\sqrt{2}}W_{t_1},\ldots, \tfrac{1}{\sqrt{2}}W_{t_k} )]|\\
&\qquad+
|\pi_N^{\ast}\pE^\w_{\beta_N,N}[(F-f\circ\Pi_{t_1,\ldots,t_k} )(X)]
- \pE[(F-f\circ\Pi_{t_1,\ldots,t_k} )(\tfrac{1}{\sqrt{2}}W )]|,
\end{align*}
which implies in particular that 
\begin{align*}
    &\dP(|\pi_N^{\ast}\pE^\w_{\beta_N,N}[F(X)]
- \pE[F(\tfrac{1}{\sqrt{2}}W )]|>2\varepsilon)\\
&\qquad\leq 
\dP(|
\pi_N^{\ast}\pE^\w_{\beta_N,N}[f(X_{t_1},\ldots, X_{t_k})]
- \pE[f(\tfrac{1}{\sqrt{2}}W_{t_1},\ldots, \tfrac{1}{\sqrt{2}}W_{t_k} )]
|>\varepsilon)\\
&\qquad\qquad  +\dP(|
\pi_N^{\ast}\pE^\w_{\beta_N,N}[(F-f\circ\Pi_{t_1,\ldots,t_k} )(X)]
- \pE[(F-f\circ\Pi_{t_1,\ldots,t_k} )(\tfrac{1}{\sqrt{2}}W )]
|>\varepsilon).
\end{align*}
The first term on the right vanishes as $N\to\infty$ due to Theorem~\ref{theo_fdd_quenched}. The second term, on the other hand, can be further upper bounded using \eqref{eq_est_diff_f_and_F}, such that 
\begin{align*}
    &\sup_{N\in \N}\dP(|
\pi_N^{\ast}\pE^\w_{\beta_N,N}[(F-f\circ\Pi_{t_1,\ldots,t_k} )(X)]
- \pE[(F-f\circ\Pi_{t_1,\ldots,t_k} )(\tfrac{1}{\sqrt{2}}W )]
|>\varepsilon)\\
&\qquad\leq \sup_{N\in \N}
\dP\big(\pi_N^{\ast}\pP^\w_{\beta_N,N}(K^c)>\tfrac{\varepsilon}{8\|F\|_{\infty}}\big)
\leq \frac{8\|F\|_{\infty}}{\varepsilon}\sup_{N\in \N}
\dE[ \pi_N^{\ast}\pP^\w_{\beta_N,N}(K^c)]\leq \delta,
\end{align*}
where we applied Markov's inequality before using \eqref{eq_choice_K_tight} again.
Because $\delta>0$ can be chosen arbitrarily small after taking the large $N$ limit, \eqref{eq_conv_in_prob_path_tested} follows. This concludes the proof.
\end{proof}

It is only left to lift the functional central limit theorem  (in $\dP$-probability) to an invariance principle as stated in Theorem~\ref{theo_main}. 
In fact, we can show more generally the equivalence of functional central limit theorem and invariance principle for random probability measures:

\begin{proposition}\label{prop_equiv_fclt_ivp}
Let $(S,d)$ be a separable, complete metric space and $(\Omega,\mathcal{G}, \mathbb{P})$ a probability space. Moreover, let $(\mathrm{P}_N^\omega)_{N\in \mathbb{N}}$ be random probability measures 
and $\mathrm{P}$ a deterministic probability measure on $(S,\mathcal{F})$, here $\mathcal{F}$ denotes the Borel-$\sigma$-algebra. 
Then the following two statements are equivalent
\begin{enumerate}[label=(\roman*)]
    \item for every $F\in C_b(S)$, $\mathrm{E}_N^\omega[F]\to \mathrm{E}[F]$ in $\mathbb{P}$-probability,
    \item $\mathrm{P}_N^\omega \weakconv \mathrm{P}$ in $\mathbb{P}$-probability,
\end{enumerate}
where $\mathrm{E}_N^\omega$ and $\mathrm{E}$ denote the expectations w.r.t. $\mathrm{P}_N^\omega$ and $\mathrm{P}$, respectively.
The statement remains true when replacing convergence in probability with almost sure convergence.
\end{proposition}

First, we remind the reader that a set of functions $\mathcal{A}\subset C_b(S)$, where $S$ Polish, is called \emph{weak convergence determining}, if for $\nu_n, \nu \in \mathcal{M}_1(S)$
\begin{align*}
    \lim_{n\to\infty}\int F \, d\nu_n  = \int F \, d\nu\, \quad \forall F\in \mathcal{A}\,,
\end{align*}
implies $\nu_n \weakconv \nu$.

The following lemma, which is a Corollary of \cite[Lemma 2]{BK10}, states that we can always find such a family of functions which is countable, provided the probability measures are defined on a Polish space. 

\begin{lemma}\label{lem_BK_sep_comp}
Let $(S,d)$ be a separable and complete metric space, then there exists a countable algebra $\mathcal{A}= \{F_i\}_{i\in \N}\subset C_b(S)$ that is weak convergence determining.
\end{lemma}
\begin{proof}
First, note that the family $\overline{\mathcal{A}}\subset C_b(S)$ of uniformly continuous functions with bounded support strongly separates points.
Then \cite[Lemma 2]{BK10} yields existence of a  countable subfamily $\mathcal{A}:= \{F_i\}_{i\in \N}\subset\overline{\mathcal{A}}$ that strongly separates points and remains closed under multiplication. Lastly applying \cite[Theorem 3.4.5]{EK86}, which states that an algebra in $C_b(S)$ that strongly separates points is convergence determining, finishes the proof.
\end{proof}
\begin{proof}[Proof of Proposition~\ref{prop_equiv_fclt_ivp}]

We begin by proving the direction (i) to (ii).
Let $\{F_i\}_{i\in \N}\subset C_b(S)$ be a countable family of functions that is weak convergence determining, which existence is guaranteed by Lemma~\ref{lem_BK_sep_comp}.
For every $i\in \N$, we have 
\begin{align}\label{eq_fclt_wdcf}
    \lim_{N\to\infty}\mathrm{E}^\w_{N}[F_i]
= \mathrm{E}[F_i],
\quad \text{in }\mathbb{P}\text{-probability}.
\end{align}
Now a diagonal argument, as we performed it in the proof of Lemma~\ref{lem_pi_system_fdd}, yields that for every sequence $(N_j)_j\subset \N$ there exists a further subsequence $(N_{j_m})_m$ and a set $\overline{\Omega}\subset \Omega$, with $\dP(\overline{\Omega})=1$, such that for every $\w\in \overline{\Omega}$
\begin{align*}
    \lim_{m\to\infty}\pi_{N_{j_m}}^{\ast}\mathrm{E}^\w_{{N_{j_m}}}[F_i]
= \mathrm{E}[F_i]\quad \forall i \in \N\,.
\end{align*}
Because $\{F_i\}_{i\in \N}$ is weak convergence determining, this implies $\mathrm{P}^\w_{{N_{j_m}}} \weakconv \mathrm{P}$, as $m\to\infty$, for every $\w \in \overline{\Omega}$. This is equivalent to weak convergence in $\mathbb{P}$-probability as stated in (ii).

The reverse direction, from (ii) to (i), is immediate. 
Moreover, replacing convergence in $\mathbb{P}$-probability with $\mathbb{P}$-almost-sure convergence, we can simply use that 
for every $F_i$ there exists a set $\Omega_i\subset \Omega$ with $\dP(\Omega_i)=1$ such that \eqref{eq_fclt_wdcf} holds pointwise for every $\w \in \Omega_i$. Taking now the countable intersection over all such  $\Omega_i$'s, we conclude $\mathrm{P}_N^\w\weakconv \mathrm{P}$ $\dP$-almost surely.
\end{proof}
\begin{proof}[Proof of Theorem~\ref{theo_main}]
The invariance principle for the polymer measures follows now directly from  the functional CLT in Proposition~\ref{prop_functional_clt} and Proposition~\ref{prop_equiv_fclt_ivp}.
\end{proof}

Instead of taking the detour via the functional central limit theorem first, we could have also argued that $\{\text{Law}(\pi_{N}^{\ast}\pP^\w_{\beta_{N},{N}})\}_N\subset \mathcal{M}_1(\mathcal{M}_1(C[0,1]))$ is tight. Together with convergence of finite-dimensional distributions, Corollary~\ref{cor_conv_fdd_quenched}, this yields a direct argument for the invariance principle. However, we want to put emphasis on the (non-trivial) equivalence of the functional CLT and the invariance principle in the case of random path measures whenever the limit is deterministic.\\

Lastly, we note that Proposition~\ref{prop_equiv_fclt_ivp} also concludes the invariance principle from the functional CLT in higher dimensions \cite[Theorem 1.2]{CY06}, which was -- to the authors' best knowledge -- not yet mentioned in the literature.
\begin{corollary}[Invariance principle for $d\geq 3$, weak disorder]\label{cor_ivp_dgeq3}
Let $\widehat{\beta}\geq 0$ such that weak disorder holds, i.e. $\lim_{N\to \infty}Z_{\widehat{\beta},N}(0,0,N,\star)>0$, then
\begin{align*}
    \pi_N^{*}\pP_{\widehat{\beta}, N}^{\w}
\weakconv \pP\big(\tfrac{1}{\sqrt{d}}W\in \cdot \,\big), \qquad \text{as }N\to\infty,
\quad \text{in }\dP\text{-probability}\,,
\end{align*}
with $\pP$ being the $d$-dimensional Wiener measure.
The statement holds in particular for all $\widehat{\beta}\in [0,\beta_c(d))$.
\end{corollary}


\subsection{Local limit theorem for the polymer marginals}
We want to close this section by proving the local limit theorem for the marginals of the polymer measure, Proposition~\ref{prop_polymer_llt}. 
Recall from \eqref{eq_disc_fdd_as_p2p_product} that the finite-dimensional distributions of the discrete polymer measure can be written as
\begin{align*}
    &\pP_{\beta_N,N}^{\w}(S_{m_1}=z_1,\ldots, S_{m_k}=z_k) \nonumber
    =
    \frac{1}{Z_{\beta_N,N}(0,0,N,\star)}
    \prod_{j=1}^{k+1}
    \ptp_{\beta_N,N}(m_{j-1},z_{j-1}\mid m_j,z_j)\,
q_{m_j-m_{j-1}}(z_{j}-z_{j-1}),
\end{align*}
where $m_0=z_0=0$, $m_{k+1}=N$ and $z_{k+1}=\star$.
Together with Proposition~\ref{prop_fact_part_func} and \cite[Theorem 2.12]{CSZ17}, this suffices to deduce Proposition~\ref{prop_polymer_llt}.


\begin{proof}[Proof of Proposition \ref{prop_polymer_llt}]
By the local limit theorem, we know that $\tfrac{N}{2}\, q_{n_i-n_{i-1}}(z_i-z_{i-1})$ converges to $p_{(t_i-t_{i-1})/2}(x_i-x_{i-1})$. It is only left to show convergence in distribution of the partition functions on the l.h.s. of \eqref{eq_polymer_llt}.
We want to show that, as $N$ diverges, 
\begin{align}\label{eq_llt_poly_repl}
\bigg\|
    \frac{1}{Z_{\beta_N,N}(0,0,N,\star)}
    &\prod_{j=1}^{k+1}
    \ptp_{\beta_N,N}(m_{j-1},z_{j-1}\mid m_j,z_j)
    -
    \prod_{j=1}^{k}
    Z_{\beta_N,N}(m_{j}^{-},\star , m_j,z_j)
    Z_{\beta_N,N}(m_{j},z_{j}, m_{j}^+,\star)
    \bigg\|_1\to 0
    ,
\end{align}
where $(m_j^{\pm})_{j=1}^k$ are non-negative integers such that 
\begin{align*}
    0\leq m_j^-<m_j<m_j^+<m_{j+1}^- <N
    \quad\text{ and }\quad
    \lim_{N\to \infty}\frac{|m_j^{\pm}-m_j|}{N}>0.
\end{align*}
For example, we can choose $m_1^-=0$ and $m^{\pm}_j= \lfloor m_j\pm \tfrac{1}{3}|m_{j\pm 1}-m_j| \rfloor$ for the remaining variables.

Using a chain of triangle inequalities, we will justify the convergence in \eqref{eq_llt_poly_repl}.
First, we note that
\begin{align*}
    &\bigg\|
\frac{1}{Z_{\beta_N,N}(0,0,N,\star)}
\prod_{j=1}^{k+1}
\ptp_{\beta_N,N}(m_{j-1},z_{j-1}\mid m_j,z_j)
 -  
 Z_{\beta_N,N}(0,\star , m_1,z_1)
 \prod_{j=2}^{k+1}
\ptp_{\beta_N,N}(m_{j-1},z_{j-1}\mid m_j,z_j)
\bigg\|_1
\end{align*}
vanishes, as $N$ tends to infinity.
The proof follows the same lines as the one of \eqref{eq_too_tired} in Lemma~\ref{lem_annealed_law}. 
Next, we replace the remaining point-to-point partition functions $(\ptp_{\beta_N,N}(m_{j-1},z_{j-1}\mid m_j,z_j))_{j=2}^{k+1}$ with their point-to-plane counterparts from Proposition~\ref{prop_fact_part_func}. For the sake of brevity, we restrict ourselves to $k=2$ for the remainder of this proof, the general case follows using a telescopic sum argument in the subsequent step. We have
\begin{align*}
        &\bigg\|
 Z_{\beta_N,N}(0,\star , m_1,z_1)
 \prod_{j=2}^{3}
\ptp_{\beta_N,N}(m_{j-1},z_{j-1}\mid m_j,z_j)\\
&\qquad-
Z_{\beta_N,N}(0,\star , m_1,z_1)
\Big(
Z_{\beta_N,N}(m_{1},z_{1}, m_2^-,\star)
Z_{\beta_N,N}(m_{1}^+,\star , m_2,z_2)
\Big)
Z_{\beta_N,N}(m_2,z_2 , N,\star)
\bigg\|_1\\
&=
\|
 Z_{\beta_N,N}(0,\star , m_1,z_1)\|_1\\
&\qquad\times
 \|
\ptp_{\beta_N,N}(m_1,z_1\mid m_2,z_2)
-
Z_{\beta_N,N}(m_{1},z_{1}, m_2^-,\star)
Z_{\beta_N,N}(m_{1}^+,\star , m_2,z_2)
\|_1\\
&\qquad\times
\|
Z_{\beta_N,N}(m_2,z_2 , N,\star)
\|_1.
\end{align*}
Here, we made use of the disorder's independence on the disjoint time intervals $(0,m_1]$, $(m_1,m_2]$ and $(m_2,N]$, and applied the fact that $\ptp_{\beta_N,N}(m_k,z_k \mid N,\star)=Z_{\beta_N,N}(m_k,z_k , N,\star)$. The middle term on the r.h.s. vanishes due to Proposition~\ref{prop_fact_part_func}, whereas the additional terms are all equal to one. This finally yields \eqref{eq_llt_poly_repl} by adding and subtracting the above introduced intermediate terms and applying the triangle inequality.\\

The last step consists of determining the limiting distribution of the partition functions in 
\begin{align*}
    \prod_{j=1}^{k}
    Z_{\beta_N,N}(m_{j}^{-},\star , m_j,z_j)
    Z_{\beta_N,N}(m_{j},z_{j}, m_{j}^+,\star).
\end{align*}
Using \cite[Theorem 2.12]{CSZ17}, see \eqref{eq_csz_weakconv_of_p2l} and the discussion thereafter, we know the limit of such point-to-plane partition functions is given by independent log-normal random variables:
\begin{align*}
    (Z_{\beta_N,N}(m_{j}^{-},\star , m_j,z_j),
    Z_{\beta_N,N}(m_{j},z_{j}, m_{j}^+,\star))_{j=1}^k\weakconv (:e^{Y^-(t_j,x_j)}:,\,:e^{Y^+(t_j,x_j)}:)_{j=1}^k,
\end{align*}
where $Y^{\pm}(t_j,x_j)$ are independent centred Gaussian random variables with variance  $\log (1-\widehat{\beta}^2)^{-1}$.
\end{proof}

%% file: sections/appendix.tex
\appendix 
\section{Transition kernel asymptotics}

In the proof of Proposition~\ref{prop_fact_part_func} we need to handle the ratio of random walk transition probabilities. The following lemma allows us to either ignore such ratios or at least bound them uniformly:

\begin{lemma}\label{lem_rw_tk_properties}
For every $x\in\R^2$ and $r>0$,
\begin{enumerate}[label=(\roman*)]
    \item  let $a_N:=(\log N)^{\gamma-1}$, $\gamma\in (0,1)$, then we have
    \begin{align}
            \sup_{\substack{z\in \sqrt{N}B(x,r)\\ \text{s.t. } q_N(z)>0}}
    \sup_{\substack{|N-n|< 2 N^{1-a_N} \\ |z-y|< 2{N}^{1/2-a_N/4}}}
        \left|
        \frac{q_n(y)}{q_N(z)}-1
        \right|\to 0 \quad \text{as }N\to\infty.
    \end{align}
    
    \item there exists a constant $C>0$ such that for all  $k\in\mathbb{N}$ 
    \begin{align}
        \sup_{\substack{z\in \sqrt{N}B(x,r)\\ \text{s.t. } q_N(z)>0}}
        \sup_{\substack{n\geq N/k\\ y\in\Z^2}}
        \frac{q_n(y)}{q_N(z)}\leq C k
    \end{align}
    for all $N$ large enough.
\end{enumerate}
\end{lemma}

\begin{proof}
(i) 
We want to apply the local limit theorem for simple random walks \cite{LL10}, which is why we write 
\begin{align*}
\sup_{\substack{|N-n|< 2 N^{1-a_N} \\ |z-y|< 2{N}^{1/2-a_N/4}}}
\left|
\frac{q_n(y)}{q_N(z)}-1
\right|
&=
\frac{1}{N\,q_N(z)}
\sup_{\substack{|N-n|< 2 N^{1-a_N} \\ |z-y|< 2{N}^{1/2-a_N/4}}}
\left|
N\,q_n(y)-N\,q_N(z)
\right|.
\end{align*}
Because $N\,q_N(z)$ converges uniformly in $z$ and its limit is lower bounded by $2\, \inf_{\widetilde{x}\in B(x,r)} p_{1/2}(\widetilde{x})$ whenever $q_N(z)>0$, we can ignore the factor in front of the supremum. Adding and subtracting $2\,p_{1/2}(z/\sqrt{N})$ yields 
\begin{align*}
   \sup_{\substack{|N-n|< 2 N^{1-a_N} \\ |z-y|< 2{N}^{1/2-a_N/4}}}
\left|
N\,q_n(y)-N\,q_N(z)
\right|
&\leq 
\sup_{\substack{|N-n|< 2 N^{1-a_N} \\ |z-y|< 2{N}^{1/2-a_N/4}}}\Big(
\left|
N\,q_n(y)-2\,p_{\frac{1}{2}}(z/\sqrt{N})
\right| 
\\&\qquad\qquad 
+
\left|
2\,p_{\frac{1}{2}}(z/\sqrt{N})-N\,q_N(z)
\right|\Big).
\end{align*}
The second term on the r.h.s. vanishes uniformly in $z$ by the local limit theorem. Because $y$ and $z$ are arbitrary close on the macroscopic scale, the first term vanishes for the same reason.\\

(ii) We begin by noting that for $z\in \sqrt{N}B(x,r)$ with $q_N(z)>0$, we have
\begin{align*}
    \sup_{\substack{n\geq N/k\\ y\in\Z^2}}
        \frac{q_n(y)}{q_N(z)}
        \leq k \, \frac{1}{N\,q_N(z)}
        \sup_{\substack{n\geq N/k\\ y\in\Z^2}}
       {n\,q_n(y)},
\end{align*}
where we may again ignore the factor $N\, q_N(z)$ for the same reason as in (i). Hence, it suffices to prove the existence of a constant $C>0$ such that
\begin{align*}
     \sup_{\substack{n\geq N/k\\ y\in\Z^2}}
       {n\,q_n(y)}\leq C\qquad \forall \, k\in \mathbb{N}\,,
\end{align*}
for $N$ large enough.
We make the following choice for $C$:
\begin{align*}
\sup_{{y\in\Z^2}} n\,q_n(y)
&\leq 
\sup_{{y\in\Z^2}}\left(
       |{n\,q_n(y)}-2p_{\frac{1}{2}}(y/\sqrt{n})|
       +|2p_{\frac{1}{2}}(y/\sqrt{n})|\right)\\
       &\leq
       \sup_{{y\in\Z^2}}
       |{n\,q_n(y)}-2p_{\frac{1}{2}}(y/\sqrt{n})|
       +|2p_{\frac{1}{2}}(0)|=:C_n.
\end{align*}
Since the first term on the r.h.s. converges in $n$ by the local limit theorem, $C_n$ is uniformly bounded in $n$.
Thus, we have $$\sup_{n\geq N/k} \sup_{{y\in\Z^2}} n\,q_n(y)\leq \sup_{n\in\mathbb{N}}C_n=:C<\infty.$$
This finishes the proof.
\end{proof}

\section{Decay of remainders in the polynomial chaos expansion}

We still owe the reader a rigorous justification for the exponential decay of second moments of $\widehat{Z}^{(k)}_{\beta_N,N}$, which we use in the proofs of Lemma~\ref{lem_non_box_limit_vanishes} and Corollary~\ref{cor_hypercontract}. 

\begin{lemma}\label{lem_exp_decay_hatZ_k}
For any $k\in \N$, we have
\begin{align*}
    \dE\left[\widehat{Z}^{(k)}_{\beta_N,N}(0,0,N,\star)^2\right]
    \leq C\,
    k^2
    \left(\sigma_N^2 R_N\right)^{\frac{k}{2}}
    \, a_N\,,
\end{align*}
where $C$ is independent of $k$ and $N$. 
\end{lemma}

The proof of this statement can be found in \cite[Section 3.4]{CSZ_KPZ}. 
The original proof makes use of more precise estimates to show that $\dE[\widehat{Z}_{\beta_N,N}(0,0,N,\star)^2]$ decays like $a_N$. One can follow the same steps using less sophisticated estimates to get an uniform estimate on $\dE[\widehat{Z}^{(k)}_{\beta_N,N}(0,0,N,\star)^2]$ in terms of $a_N$ instead, which yields the same qualitative bound.
We sketch the argument for the sake of completeness:\\

Considering $\widehat{Z}^{(k)}_{\beta_N,N}(0,0,N,\star)$ for some $k\leq N$,
there is at least one sample $(n_j,z_j)$ outside the box $A_N^+(0,0)$. Thus, 
\begin{align*}
    \dE\left[\widehat{Z}^{(k)}_{\beta_N,N}(0,0,N,\star)^2\right]
    \leq
    \sigma_N^{2k} 
    \sum_{\substack{1\leq l_1, \ldots, l_k\leq N\\ z_1,\ldots, z_k \in \Z^2}}
    \sum_{j=1}^k
    \left(
    \mathds{1}_{l_j> \frac{1}{k}N^{1-a_N}}
    +
    \mathds{1}_{l_j\leq \frac{1}{k}N^{1-a_N},\, |z_j|\geq \frac{1}{k}N^{1/2-a_N/4}}
    \right)
    \prod_{i=1}^k q_{l_i}^2(z_i)
    \,,
\end{align*}
where we extended the range of time-differences $l_i=n_i-n_{i-1}$ to all of $\{1,\ldots, N\}$. Using once more the identity $\sum_{l=1}^N \sum_{z\in \Z^2}q_l^2(z)=R_N$, we have 
\begin{align*}
    \dE\left[\widehat{Z}^{(k)}_{\beta_N,N}(0,0,N,\star)^2\right]
    \leq k\,
    \sigma_N^{2k} R_N^{k-1} 
    \sum_{\substack{1\leq l \leq N\\ z\in\Z^2}} 
    \left(
    \mathds{1}_{l> \frac{1}{k}N^{1-a_N}}
    +
    \mathds{1}_{l\leq \frac{1}{k}N^{1-a_N}\,, |z|\geq \frac{1}{k}N^{1/2-a_N/4}}
    \right)
    q_{l}^2(z)\,.
\end{align*}
Now, the contribution of the two indicator functions can be considered separately. We follow the exact same steps as in Section 3.4 of \cite{CSZ_KPZ}:
\begin{itemize}
    \item For the contribution of large time-jumps, we have
    \begin{align*}
        \frac{1}{R_N}
        \sum_{\frac{1}{k}N^{1-a_N}\leq l \leq N} q_{2l}(0)
        \leq C \frac{1}{R_N}
        \sum_{\frac{1}{k}N^{1-a_N}\leq l \leq N} \frac{1}{l}
        \leq 
        C'\,  \frac{a_N \log N + \log k}{\log N}
        \leq 2 C'\,  k\, a_N\,,
    \end{align*}
    where we used additionally the crude estimates $\tfrac{\log k}{\log N}\leq k \, a_N$ in the last inequality.
    \item On the other hand, the contribution of the second term is upper bounded by 
    \begin{align*}
        \frac{1}{R_N}
        \sum_{1\leq l \leq \frac{1}{k}N^{1-a_N}}
        \sum_{|z|>\frac{1}{k}N^{1/2-a_N/4}}
        q_l^2(z)
        \leq C
        e^{-\eta \frac{N^{\frac{a_N}{2}}}{k}}\,,
    \end{align*}
    for some uniform $\eta>0$, 
    using Gaussian estimates for the simple random walk. Then, for $N$ large enough
    \begin{align}
        C\, k\, (\sigma_N^{2} R_N)^k
        e^{-\eta \frac{N^{\frac{a_N}{2}}}{k}}
        \leq 
        \begin{cases}
        C\, k (\sigma_N^{2} R_N)^{k}
            e^{-\eta N^{\frac{a_N}{4}}}
            \,,&\quad \text{if } k\leq (N^{a_N/2})^{\frac{1}{2}}\,,\\
            C\, k (\sigma_N^{2} R_N)^{\frac{k}{2}}
            (\widehat{\beta}+\delta)^{N^{\frac{a_N}{4}}}\,,&\quad \text{if } k> (N^{a_N/2})^{\frac{1}{2}} \,,
        \end{cases}
    \end{align}
    where $\delta>0$ small enough such that $\widehat{\beta}+\delta<1$. 
    Note that, $c^{N^{\frac{a_N}{4}}}=o(a_N)$ for any $c\in (0,1)$, because
    \begin{align*}
    (\log N)^\gamma
        {c^{N^{\frac{a_N}{4}}}}
        =(\log N)^\gamma
        c^{\exp(\frac{1}{4}(\log N)^\gamma)}
        \leq (\log N)^\gamma
        c^{\frac{1}{4}(\log N)^\gamma} \to 0
        \,,\quad \text{ as }N\to \infty \,,
    \end{align*}
    which includes in particular the case $c=\max \{\widehat{\beta}+\delta, e^{-\eta}\}$.
\end{itemize}
Adding up all above estimates yields the desired upper bound of $\dE[\widehat{Z}^{(k)}_{\beta_N,N}(0,0,N,\star)^2]$ from Lemma~\ref{lem_exp_decay_hatZ_k}.

%% file: main.bbl
\begin{thebibliography}{AKQ14b}

\bibitem[AKQ14a]{ALKQ13}
Tom Alberts, Konstantin Khanin, and Jeremy Quastel.
\newblock The continuum directed random polymer.
\newblock {\em Journal of Statistical Physics}, 154(1):305--326, 2014.

\bibitem[AKQ14b]{ALKQ14}
Tom Alberts, Konstantin Khanin, and Jeremy Quastel.
\newblock The intermediate disorder regime for directed polymers in dimension
  1+1.
\newblock {\em The Annals of Probability}, 42(3):1212--1256, May 2014.

\bibitem[AZ96]{AZ96}
Sergio Albeverio and Xian~Yin Zhou.
\newblock A martingale approach to directed polymers in a random environment.
\newblock {\em Journal of Theoretical Probability}, 9(1):171--189, January
  1996.

\bibitem[BC95]{BC95}
Lorenzo Bertini and Nicoletta Cancrini.
\newblock The stochastic heat equation: Feynman-kac formula and intermittence.
\newblock {\em Journal of Statistical Physics}, 78(5-6):1377--1401, March 1995.

\bibitem[BC98]{BC98}
Lorenzo Bertini and Nicoletta Cancrini.
\newblock The two-dimensional stochastic heat equation: renormalizing a
  multiplicative noise.
\newblock {\em Journal of Physics A: Mathematical and General}, 31(2):615--622,
  January 1998.

\bibitem[BD00]{Br00}
Eric Brunet and Bernard Derrida.
\newblock Probability distribution of the free energy of a directed polymer in
  a random medium.
\newblock {\em Physical Review E}, 61(6):6789--6801, June 2000.

\bibitem[BG97]{BG97}
Lorenzo Bertini and Giambattista Giacomin.
\newblock Stochastic {Burgers} and {KPZ} equations from particle systems.
\newblock {\em Communications in Mathematical Physics}, 183(3):571--607,
  February 1997.

\bibitem[Bil99]{Bi99}
Patrick Billingsley.
\newblock {\em Convergence of Probability Measures}.
\newblock John Wiley {\&} Sons, Inc., July 1999.

\bibitem[BK10]{BK10}
Douglas Blount and Michael~A. Kouritzin.
\newblock On convergence determining and separating classes of functions.
\newblock {\em Stochastic Processes and their Applications},
  120(10):1898--1907, September 2010.

\bibitem[Bol89]{Bo89}
Erwin Bolthausen.
\newblock A note on the diffusion of directed polymers in a random environment.
\newblock {\em Communications in Mathematical Physics}, 123(4):529--534,
  December 1989.

\bibitem[CDR10]{CDR10}
Pasquale Calabrese, Pierre~Le Doussal, and Alberto Rosso.
\newblock Free-energy distribution of the directed polymer at high temperature.
\newblock {\em {EPL} (Europhysics Letters)}, 90(2):20002, apr 2010.

\bibitem[CNN20]{CNN20}
Cl{\'e}ment Cosco, Shuta Nakajima, and Makoto Nakashima.
\newblock Law of large numbers and fluctuations in the sub-critical and {$L^2$}
  regions for {SHE} and {KPZ} equation in dimension $d\geq 3$.
\newblock {\em arXiv preprint arXiv:2005.12689}, 2020.

\bibitem[Com17]{Co17}
Francis Comets.
\newblock {\em Directed Polymers in Random Environments}.
\newblock Springer International Publishing, 2017.

\bibitem[CSY04]{CSY04}
Francis Comets, Tokuzo Shiga, and Nobuo Yoshida.
\newblock Probabilistic analysis of directed polymers in a random environment:
  a review.
\newblock In {\em Stochastic Analysis on Large Scale Interacting Systems}.
  Mathematical Society of Japan, 2004.

\bibitem[CSZ14]{CSZ_PTRF}
Francesco Caravenna, Rongfeng Sun, and Nikos Zygouras.
\newblock The continuum disordered pinning model.
\newblock {\em Probability Theory and Related Fields}, 164(1-2):17--59,
  December 2014.

\bibitem[CSZ17a]{CSZ_EMS}
Francesco Caravenna, Rongfeng Sun, and Nikos Zygouras.
\newblock Polynomial chaos and scaling limits of disordered systems.
\newblock {\em Journal of the European Mathematical Society}, 19(1):1--65,
  2017.

\bibitem[CSZ17b]{CSZ17}
Francesco Caravenna, Rongfeng Sun, and Nikos Zygouras.
\newblock Universality in marginally relevant disordered systems.
\newblock {\em The Annals of Applied Probability}, 27(5):3050--3112, October
  2017.

\bibitem[CSZ20]{CSZ_KPZ}
Francesco Caravenna, Rongfeng Sun, and Nikos Zygouras.
\newblock The two-dimensional {KPZ} equation in the entire subcritical regime.
\newblock {\em Ann. Probab.}, 48(3):1086--1127, May 2020.

\bibitem[CY06]{CY06}
Francis Comets and Nobuo Yoshida.
\newblock Directed polymers in random environment are diffusive at weak
  disorder.
\newblock {\em The Annals of Probability}, 34(5):1746--1770, September 2006.

\bibitem[CZ21]{CZ21}
Cl{\'e}ment Cosco and Ofer Zeitouni.
\newblock Moments of partition functions of 2d {G}aussian polymers in the weak
  disorder regime.
\newblock {\em arXiv preprint arXiv:2112.03767}, 2021.

\bibitem[EK86]{EK86}
Stewart~N. Ethier and Thomas~G. Kurtz.
\newblock {\em Markov Processes}.
\newblock John Wiley {\&} Sons, Inc., March 1986.

\bibitem[Fen12]{Fe12}
Zi~Sheng Feng.
\newblock Diffusivity of rescaled random polymer in random environment in
  dimensions 1 and 2.
\newblock {\em arXiv preprint arXiv:1201.6215}, 2012.

\bibitem[HH85]{HH85}
David~A. Huse and Christopher~L. Henley.
\newblock Pinning and roughening of domain walls in {I}sing systems due to
  random impurities.
\newblock {\em Physical Review Letters}, 54(25):2708--2711, June 1985.

\bibitem[IS88]{IS88}
John~Z. Imbrie and Thomas Spencer.
\newblock Diffusion of directed polymers in a random environment.
\newblock {\em Journal of Statistical Physics}, 52(3-4):609--626, August 1988.

\bibitem[Jun21a]{Ju21}
Stefan Junk.
\newblock The central limit theorem for directed polymers in weak disorder,
  revisited.
\newblock {\em arXiv preprint arXiv:2105.04082}, 2021.

\bibitem[Jun21b]{Ju21b}
Stefan Junk.
\newblock New characterization of the weak disorder phase of directed polymers
  in bounded random environments.
\newblock {\em Communications in Mathematical Physics}, 389(2):1087--1097,
  November 2021.

\bibitem[Kal02]{Ka02}
Olav Kallenberg.
\newblock {\em Foundations of modern probability}.
\newblock Probability and its Applications (New York). Springer-Verlag, New
  York, second edition, 2002.

\bibitem[Kif97]{Ki97}
Yuri Kifer.
\newblock The {Burgers} equation with a random force and a general model for
  directed polymers in random environments.
\newblock {\em Probability Theory and Related Fields}, 108(1):29--65, May 1997.

\bibitem[Lax02]{La02}
Peter~D. Lax.
\newblock {\em Functional Analysis}.
\newblock Pure and Applied Mathematics: A Wiley Series of Texts, Monographs and
  Tracts. Wiley, 2002.

\bibitem[Led01]{Le01}
Michel Ledoux.
\newblock {\em The Concentration of Measure Phenomenon}.
\newblock Mathematical surveys and monographs. American Mathematical Society,
  2001.

\bibitem[LL10]{LL10}
Greg~F. Lawler and Vlada Limic.
\newblock {\em Random Walk: A Modern Introduction}.
\newblock Cambridge Studies in Advanced Mathematics. Cambridge University
  Press, 2010.

\bibitem[LZ21]{LZ21}
Dimitris Lygkonis and Nikos Zygouras.
\newblock Moments of the 2d directed polymer in the subcritical regime and a
  generalisation of the {E}rd\"os-{T}aylor theorem.
\newblock {\em arXiv preprint arXiv:2109.06115}, 2021.

\bibitem[LZ22]{LZ20}
Dimitris Lygkonis and Nikos Zygouras.
\newblock Edwards{\textendash}{W}ilkinson fluctuations for the directed polymer
  in the full {$L^2$}-regime for dimensions d$\geq$3.
\newblock {\em Annales de l{\textquotesingle}Institut Henri Poincar{\'{e}},
  Probabilit{\'{e}}s et Statistiques}, 58(1), February 2022.

\bibitem[NN21]{NN21}
Shuta Nakajima and Makoto Nakashima.
\newblock Fluctuations of two-dimensional stochastic heat equation and {KPZ}
  equation in subcritical regime for general initial conditions.
\newblock {\em arXiv preprint arXiv:2103.07243}, 2021.

\bibitem[Rud86]{Ru86}
Walter Rudin.
\newblock {\em Real and Complex Analysis}.
\newblock McGraw-Hill series in higher mathematics. McGraw-Hill Professional,
  New York, NY, 3 edition, September 1986.

\bibitem[Sin95]{Si95}
Yakov Sinai.
\newblock A remark concerning random walks with random potentials.
\newblock {\em Fundamenta Mathematicae}, 147(2):173--180, 1995.

\bibitem[SZ96]{SZ96}
Renming Song and Xian~Yin Zhou.
\newblock A remark on diffusion of directed polymers in random environments.
\newblock {\em Journal of Statistical Physics}, 85(1-2):277--289, October 1996.

\bibitem[Var06]{Va06}
Vincent Vargas.
\newblock A local limit theorem for directed polymers in random media: the
  continuous and the discrete case.
\newblock {\em Annales de l{\textquotesingle}Institut Henri Poincare (B)
  Probability and Statistics}, 42(5):521--534, September 2006.

\end{thebibliography}
